\documentclass[10pt]{article}
\usepackage{amsfonts}
\usepackage{amssymb}
\usepackage{amsmath}
\usepackage{amsthm}
\usepackage{latexsym}
\usepackage{graphicx}
\usepackage{url}
\usepackage{color}
\newcommand{\V}[1]{\boldsymbol{#1}}

\newcommand{\Section}[1]{
   \refstepcounter{section}
   \bigskip\noindent
   {\large\bf\hbox{\thesection~~}#1}\par
   \nopagebreak
   \medskip
   \renewcommand{\theequation}{\thesection.\arabic{equation}}
   \setcounter{equation}{0}
   \setcounter{subsection}{0}
}

\newtheorem{thm}{Theorem}[section]
\newtheorem{lem}{Lemma}[section]
\newtheorem{rem}{Remark}[section]
\newtheorem{pro}{Proposition}[section]

\newcommand{\reff}[1]{(\ref{#1})}
\newcommand{\Proof}{\noindent {\bf Proof.~}\ }
\newcommand{\Endproof}{$\hfill\Box$}
\newcommand{\Thm}[1]{Theorem #1}
\newcommand{\Lem}[1]{Lemma #1}
\newcommand{\Rem}[1]{Remark #1}
\newcommand{\Pro}[1]{Proposition #1}

\newcommand{\be}{\begin{equation}}
\newcommand{\ee}{\end{equation}}
\newcommand{\ba}{\begin{array}}
\newcommand{\ea}{\end{array}}
\newcommand{\ben}{\begin{eqnarray}}
\newcommand{\een}{\end{eqnarray}}
\newcommand{\bn}{\begin{eqnarray*}}
\newcommand{\en}{\end{eqnarray*}}
\newcommand{\p}{\partial}
\newcommand{\udl}[1]{\underline{#1}}
\newcommand{\updl}[1]{\hat{\underline{#1}}}
\newcommand{\cTh}{{\cal T}_h}
\newcommand{\cThl}{{\cal T}_{h,1}}
\newcommand{\cThr}{{\cal T}_{h,2}}
\newcommand{\cThm}{{\cal T}_{h,i}}
\newcommand{\cIh}{{\cal I}_h}
\newcommand{\cF}{{\cal F}}
\newcommand{\cFT}{{\cal F}_{\partial T}}
\newcommand{\ife}{\p T\cap \Gamma}
\newcommand{\uife}{\p T\backslash \Gamma}
\newcommand{\Ve}{{\V \epsilon}}
\newcommand{\Vs}{{\V \sigma}}
\newcommand{\mR}{\mathbb{R}}
\newcommand{\mP}{\mathbb{P}}
\newcommand{\bE}{{\bf E}}
\newcommand{\bG}{{\bf G}}
\newcommand{\bt}{{\V \tau}}
\newcommand{\bS}{{\bf S}}
\newcommand{\bK}{{\bf K}}
\newcommand{\bI}{{\bf I}}
\newcommand{\cEh}{{\cal E}_h}
\newcommand{\pk}{\pi_k}
\newcommand{\zt}{\udl{\zeta}}

\newcommand{\norm}[1]{\left\Vert#1\right\Vert}
\newcommand{\norme}[1]{\left\Vert{\hskip -2.6pt}\left\vert #1 \right\vert{\hskip -2.6pt}\right\Vert}
\newcommand{\abs}[1]{\left\vert#1\right\vert}
\newcommand{\pd}[1]{\left\langle #1\right\rangle}
\newcommand{\set}[1]{\left\{#1\right\}}
\newcommand{\av}[1]{\left\{#1\right\}}
\newcommand{\jm}[1]{\left[#1\right]}

\newcommand{\oo}{\Omega_1\cup\Omega_2}
\newcommand{\nn}{\nonumber}
\begin{document}
\centerline{\Large \bf A Hybrid High-Order Method for}

\vskip .2cm
\renewcommand{\thefootnote}{\fnsymbol{footnote}}
\centerline{{\Large \bf the Elasticity Problem with Linear Slip Interface}}

\vskip .5cm
\renewcommand{\thefootnote}{\arabic{footnote}}
\centerline{{\bf Erik Burman}$\ddag$\footnote{Corresponding author E-mail address: e.burman@ucl.ac.uk.} and {\bf Peiqi Huang}$\dagger$\footnote{E-mail address: pqhuang@njfu.edu.cn.}}

\vskip .3cm \centerline{$\ddag$\it
Department of Mathematics, University College London,}
\centerline{\it Gower Street, London WC1E 6BT, UK.}

\vskip .3cm \centerline{$\dagger$\it Department of Applied
Mathematics, Nanjing Forestry University,}
\centerline{\it Nanjing 210037, People's  Republic of China.}

\vskip .5cm \noindent\rule[2mm]{\textwidth}{.1pt}

\noindent {\bf Abstract}
\vskip .5cm

We design and analyse a hybrid high-order (HHO) method for the linear
elasticity problem with a linear slip interface, i.e.\ with interface
conditions of spring type, in which the jump of the displacement across the
interface is proportional to the traction through a compliancy tensor
$\bK=\alpha\bI+(\beta-\alpha)\udl{n}\otimes\udl{n}$.  The mesh is fitted to
the interface, the discrete unknowns are polynomials of degree $k\ge 1$ on the
mesh faces and of degree $k+1$ in the mesh cells, and general polytopal cells
are allowed.  The two novelties of the method lie in the local symmetric
strain reconstruction, which incorporates the interface condition, and in the
interface stabilisation, which is built from a regularised interface stiffness
$\bS_h$ in the spirit of Hansbo and Hansbo \cite{HH04}.  As a consequence, a
single formulation covers the whole range of compliancies, from the perfectly
bonded interface $\alpha=\beta=0$, where the interface stabilisation acts as a
Nitsche-type penalty, to the traction-free interface obtained as
$\alpha,\beta\to+\infty$, and no unknown is attached to the jump.  We prove
that the discrete bilinear form is coercive and that the errors converge as
$h^{k+1}$ in the energy norm and as $h^{k+2}$ in the $L^2$ norm, with
constants that are independent of the compliancy parameters and of the
Lam\'e coefficient $\lambda$, so that the method is also locking free.
Various numerical examples and comparisons are provided to confirm the
theoretical results.



\noindent\rule[2mm]{\textwidth}{.1pt}

\Section{Introduction.}

Composite materials, geological media and glued or welded assemblies are
naturally modelled as elastic bodies made of several components separated by
material interfaces.  If the components are perfectly bonded, the displacement
and the normal traction are both continuous across the interface, and only the
strain is discontinuous.  In many situations of practical interest, however,
the bonding is imperfect: a thin adhesive layer, a damaged zone or a
micro-cracked region between the two materials allows the two sides to slide
or to separate.  A standard way of modelling such an interface, which avoids
resolving the thin layer, is to keep the traction continuous but to let the
displacement jump be proportional to it,
\begin{align*}
  \big[\Vs(\udl{u})\udl{n}\big]=\udl{0},\qquad
  \jm{\udl{u}}=-\bK\Vs(\udl{u})\udl{n}\qquad\mbox{on }\Gamma,
\end{align*}
where the symmetric positive semi-definite tensor $\bK$ measures the
compliancy of the interface.  For an isotropic interface,
$\bK=\alpha\bI+(\beta-\alpha)\udl{n}\otimes\udl{n}$, where $\alpha\geq 0$ and
$\beta\geq 0$ are the tangential and the normal compliancy, respectively; see
\cite{ZM99} and the references therein.  These conditions are known as linear
slip, spring-type or imperfect interface conditions, and they interpolate
between two extreme regimes: for $\bK=\bf{0}$ the interface is perfectly
bonded and the displacement is continuous, whereas in the limit
$\alpha,\beta\to+\infty$ the two subdomains decouple and $\Gamma$ becomes a
traction-free boundary for each of them.

The numerical approximation of such problems raises three distinct
difficulties.  First, the exact solution is discontinuous across $\Gamma$, so
that either the mesh resolves the interface, or the discretisation must be
able to represent the jump inside the mesh cells. Second, the compliance parameters typically span several orders of magnitude, making any discretisation whose stability or accuracy deteriorates in either of the two extreme regimes described above of limited practical value. In the limit of perfect bonding, the interface condition becomes a constraint that must be enforced, for example weakly through a Nitsche-type penalty method, whose penalty parameter must be scaled appropriately with the mesh size. Conversely, for highly compliant interfaces, it is essential to represent the physical interface stiffness accurately.  Third,
if one of the materials is nearly incompressible, the second Lam\'e
coefficient $\lambda$ becomes very large and the discretisation must be
locking free, i.e.\ its accuracy must not deteriorate as
$\lambda\to+\infty$.

The three difficulties have been addressed separately in the literature.
Hansbo and Hansbo \cite{HH04} proposed an unfitted finite element method for
the problem at hand, based on a doubling of the degrees of freedom in the
cells cut by the interface and on a Nitsche-type formulation involving a
regularised interface stiffness
$\bS_h=(h/\delta+\bK)^{-1}$; the resulting method treats the whole range of
compliancies with one and the same formulation and one and the same code, and
was shown to be optimally convergent, uniformly in $\alpha$ and $\beta$, for
piecewise affine approximations.  Locking-free discontinuous Galerkin methods
for nearly incompressible elasticity, based on Nitsche's method, were analysed
in \cite{HL02}.  Hybrid high-order (HHO) methods, introduced in \cite{DEL14}
for diffusion problems and in \cite{DE15} for linear elasticity, provide
arbitrary-order, locking-free discretisations on general polytopal meshes; they
are built from a local reconstruction operator and a local stabilisation
operator, they support hanging nodes and polytopal cells, and the cell
unknowns can be eliminated locally so that the global problem only couples the
face unknowns.  Cardenas and Solano \cite{CS24} developed a high-order unfitted hybridizable discontinuous Galerkin (HDG) method for linear elasticity and established optimal error estimates. In \cite{HWX23} and \cite{WG25}, HDG and weak Galerkin methods, respectively, were employed to solve linear elasticity interface problems under the assumption that the displacement remains continuous across the interface. Carstensen and Tran \cite{CT25} introduced a reconstruction operator for the linearized Green strain tensor, thereby avoiding the classical HHO decomposition into deviatoric and spherical components. To the best of our knowledge, however, no arbitrary-order method capable of handling the full spectrum of interface compliances has yet been rigorously analysed.

The purpose of this paper is to fill this gap.  We design an HHO method of
arbitrary order $k\geq 1$ for the linear slip interface problem on meshes that
are fitted to $\Gamma$, and we prove optimal error estimates with constants
that are independent of the compliancy parameters and of $\lambda$.  The
method has two specific ingredients.  The first one is a local symmetric
strain reconstruction, defined in \reff{gradient-reconstruct}, in which the
interface condition is built into the definition of the operator through the
factor $\bI-\bK\bS_h$; here $\bS_h$ is a regularised interface stiffness which
combines the idea of \cite{HH04} with the mesh-dependent weights that are
natural for HHO stabilisations, and which is defined in
\reff{penalty-matrix} by perturbing the eigenvalues $(\beta,\alpha,\alpha)$ of
$\bK$ by $(h_T/(\mu+\lambda),h_T/\mu,h_T/\mu)$.  The two extreme regimes are
then recovered automatically: when $\alpha=\beta=0$ one has
$\bI-\bK\bS_h=\bI$ and $\bS_h$ behaves like a Nitsche penalty of size
$\mu/h_T$ and $(\mu+\lambda)/h_T$ in the tangential and normal directions,
whereas for $\alpha,\beta\gg h_T/\mu$ one has $\bS_h\simeq \bK^{-1}$, which is
the physical interface stiffness, and $\bI-\bK\bS_h\simeq \bf{0}$.  The
second ingredient is the interface stabilisation in
\reff{stabilization-term}, which penalises the distance between the trace of
the cell unknown and the face unknown in the $\bS_h$-weighted norm.  
Our main results are the energy error estimate of order $h^{k+1}$
(\Pro{\ref{Pro-energy-Err}}) and the $L^2$ error estimate of order $h^{k+2}$
(\Pro{\ref{Pro-L2-Err}}), both with constants independent of $\lambda$,
$\alpha$ and $\beta$.  These estimates are confirmed numerically in
\S\ref{sec-numerics}, where we also verify the robustness of the method with
respect to the compliancy over eight orders of magnitude, including the
degenerate cases $\alpha=0$ or $\beta=0$, and in the quasi-incompressible
limit.

The rest of the paper is organized as follows.  In \S 2, we introduce the
elasticity problem with a linear slip interface and its weak formulation.  The
HHO method is described in \S 3, where the reconstruction operator, the
regularised interface stiffness and the stabilisation are defined.  In \S 4 we
collect the technical tools used in the analysis, namely the properties of the
interface stiffness, the standard inverse and trace inequalities, the local
stability of the reconstruction and the interpolation operators.  The stability
and error analysis are carried out in \S 5.  Numerical experiments are presented
in \S 6, which support our theoretical results.  Conclusions and possible
extensions are presented in the last section.

\begin{figure}[htbp]
\centering
\includegraphics[scale=0.85]{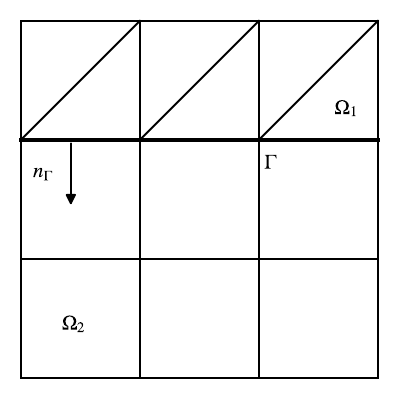}
\caption{A domain $\Omega$ is divided into subdomains $\Omega_1$ and $\Omega_2$ by the interface $\Gamma$.}\label{unfitted_mesh}
\end{figure}

\Section{An elasticity problem with linear slip interface.}
Let $\Omega$ be a connected, convex polygonal domain in ${\mathbb R}^d,d=2,{\rm or}\ 3$, with boundary $\p\Omega$ and an internal smooth interface $\Gamma=\p\Omega_1\cap\p\Omega_2$ dividing $\Omega$ into two subdomains $\Omega_1$ and $\Omega_2$, see Fig. \ref{unfitted_mesh}. Throughout this paper we shall use subscripts to denote the restriction of a function to the subdomain $\Omega_i$. Vectors and tensors are typed with underline and bold face, respectively. Thus, $\udl{u}=[u^i]_{i=1}^d$ may denote a vector valued function in $\Omega$ with components $u^i$, while $\udl{u}_i=\udl{u}|_{\Omega_i}$ denotes its restriction to $\Omega_i$. For any sufficiently regular function $\udl{u}$ in $\Omega_1\cup\Omega_2$, we define the jump of $\udl{u}$ on $\Gamma$ by $[\udl{u}]=\udl{u}_1|_\Gamma-\udl{u}_2|_\Gamma$.

We consider the following elasticity problem having discontinuities in the Lam\'e parameters along the material interface $\Gamma$: Find the displacement $\udl{u}:\Omega\rightarrow{\mathbb R}^d$ such that
\be\label{elasticity-interface}
  \left \{
    \ba {ll}
      {\V \sigma}(\udl{u})=2\mu {\V \epsilon}(\udl{u}) + \lambda \nabla \cdot \udl{u}{\bf I}  \quad & {\rm in} \quad \Omega_1\cup\Omega_2, \\
      -\nabla \cdot {\V \sigma}(\udl{u}) = \udl{f} \quad & {\rm in} \quad \Omega_1\cup\Omega_2, \\
      \big[{\V \sigma}(\udl{u})\udl{n}\big]=\udl{0} \quad & {\rm on} \quad \Gamma, \\
      \jm{\udl{u}}=-\bK \Vs(\udl{u})\udl{n} \quad & {\rm on} \quad \Gamma, \\
      \udl{u}=\udl{0} \quad & {\rm on} \quad \p\Omega,
    \ea
  \right.
\ee
where the external force $\udl{f}\in L^2(\Omega)^d$. Here
\bn
  \lambda=\frac{E\nu}{(1-2\nu)(1+\nu)},\qquad
  \mu=\frac{E}{2(1+\nu)},
\en
are the Lam\'e parameters, satisfying $0<c<\mu <C$ and $\lambda >0$. $E$ is the Young's modulus and $\nu$ is the Poisson's ratio. The strain tensor ${\V \epsilon}(\udl{u})=\frac12(\nabla \udl{u}+\nabla\udl{u}^T)$, the identity tensor ${\bf I}=[\delta^{ij}]_{i,j=1}^d$ with $\delta^{ij}=1$ if $i=j$ and $\delta^{ij}=0$ if $i\neq j$, and $\udl{n}$ is the outward unit normal vector to $\Omega_1$. ${\bf K}$ is a positive semi-definite matrix representing the compliancy of the interface. In this paper, we consider only the case of elastic isotropy, i.e. ${\bf K}$ can be written in the following form
\be\label{Definiteion-K}
  {\bf K}=\alpha {\bf I} +(\beta-\alpha)\udl{n}\otimes\udl{n}\qquad \mbox{or} \qquad K^{ij}=\alpha\delta^{ij}+(\beta-\alpha)n^in^j,
\ee
where $\alpha$ and $\beta$ are nonnegative constants denoting the compliancy in the tangential and the normal directions of the interface, respectively \cite{ZM99}. Since the material properties are different in each region, we set $\mu=\mu_i,\ \lambda=\lambda_i$ in $\Omega_i$ for $i=1,2,$ and assume ${\bf K}$ is constant matrix on $\Gamma$.

Let the interface stiffness ${\bf S}$ be defined by
\be\label{Definition-S}
  {\bf S}=\begin{cases}
    \ba {ll}
      {\bf K}^{-1}  \quad & {\rm for} \quad \alpha>0,\ \beta>0, \\
      \alpha^{-1}({\bf I}-\udl{n}\otimes\udl{n})  \quad & {\rm for} \quad \alpha>0,\ \beta=0, \\
      \beta^{-1}\udl{n}\otimes\udl{n}  \quad & {\rm for} \quad \alpha=0,\ \beta>0, \\
      0  \quad & {\rm for} \quad \alpha=0,\ \beta=0,
    \ea
  \end{cases}
\ee
and define the space $\udl{V}$ in the following
\bn
  \udl{V}=\big\{ \udl{v}\in \udl{V}_1\times\udl{V}_2 \big|\ [\udl{v}]={\bf KS}[\udl{v}] \big\}\quad {\rm where}\quad
  \udl{V}_i=\big\{ \udl{v}_i\in H^1(\Omega_i, {\mathbb R}^d) \big|\ \udl{v}_i|_{\p\Omega}=\udl{0} \big\},i=1,2.
\en

Now we define our variational problem as follows: find $\udl{u}=(\udl{u}_1,\udl{u}_2)\in \udl{V}$ such that
\be\label{weak-form}
  a(\udl{u},\udl{v})=l(\udl{v})\qquad \forall \udl{v}\in \udl{V}.
\ee
Here $a(\udl{u},\udl{v})=( {\V \sigma}(\udl{u}),{\V \epsilon}(\udl{v}) )_{\oo}+<{\bf S}\jm{\udl{u}},\jm{\udl{v}}>_\Gamma$, where $( {\V \sigma}(\udl{u}),{\V \epsilon}(\udl{v}) )_{\Omega_i}=2\mu_i( {\V \epsilon}(\udl{u}),{\V \epsilon}(\udl{v}) )_{\Omega_i} +\lambda_i (\nabla\cdot\udl{u},\nabla\cdot\udl{v})_{\Omega_i}$ and $l(\udl{v})=(\udl{f},\udl{v})_\Omega$. The variational problem \reff{weak-form} is equivalent to the original problem \reff{elasticity-interface}, see Section 2 in \cite{HH04}.

\Section{The HHO method.}

Let $\{\cTh\}_{h>0}$ be a family of regular quasi-uniform triangulations of $\Omega$ with the mesh size $h$, the mesh $\cTh$ is composed of nonempty disjoint cells such that $\overline{\Omega}=\bigcup_{T\in\cTh} \overline{T}.$ The mesh cells are conventionally open subsets in ${\mathbb R}^d$ (not necessarily convex), and they can have a polygonal/polyhedral shape with straight edges (if d = 2) or planar faces (if d = 3). This setting in particular allows for meshes with hanging nodes. The mesh faces are collected in the set ${\cal F}_h$, the set of boundary faces is denoted by ${\cal F}_h^b$. In a nutshell, each mesh $\cTh$ admits a matching simplicial submesh $\cTh'$ such that any cell (or face) of $\cTh'$ is a subset of a cell (or face) of $\cTh$. Moreover, there exists a mesh-regularity parameter $\rho>0$ such that, for all $h>0$, all $T\in\cTh$ and all $S\in\cTh'$ such that $S\subset T,$ we have $\rho h_S\leq r_S$ and $\rho h_T\leq h_S$, where $h_T$ denotes the diameter of the cell $T$ and $r_S$ denotes the inradius of the simplex $S$.

We assume that the triangulation $\cTh$ is fitted with the interface $\Gamma$, i.e. $\cTh=\cThl\cup\cThr$, and the partition of the interface $\Gamma$ induced by $\cThm,\,i=1,2$ is denoted by $\cIh$.
\bn
  \cThm:=\big\{ T\in\cTh|\ T\in \Omega_i \big\}\quad i=1,2,\qquad \cIh:=\big\{ F\subset\Gamma|\ F=T_1\cap T_2,\,
    T_1\in\cThl,T_2\in\cThr \big\}.
\en

Let $k\geq 1$ be the polynomial degree. The discrete unknowns for the displacement are piecewise polynomials of degree $k$ attached to the mesh faces and of degree $k+1$ attached to the mesh cells. For any subset $S\subset{\mathbb R}^d$ consisting of one mesh (sub)cell or one mesh (sub)face, and for all $l\in \mathbb{N}$, we denote $\mathbb{P}^l(S)$ (resp. $\mathbb{P}^l(S)^d$, $\mathbb{P}^l(S)^{d\times d}$) the space of scalar-valued (resp. vector-valued, symmetric matrix-valued) polynomials in $S$ of degree at most $l$. We also denote $(\cdot,\cdot)_S$ the $L^2$-scalar product on $S$ and $\norm{\cdot}_S$ the associated norm, replace $(\cdot,\cdot)_S$ by $\pd{\cdot,\cdot}_S$ when $S\subset{\mathbb R}^{d-1}$. Whenever $S=\emptyset$, we abuse the notation by writing $\mathbb{P}^l(S)={0}$ and $(\cdot,\cdot)_S=0$.

For any $T\in\cTh$, we set $\cFT:=\{F\in{\cal F}_h|\ F\subset\p T\}$, $\cF_{\uife}:=\{F \notin\cIh|\ F\in\cFT\}$ and $\cF_{\ife}:=\{F \in\cIh|\ F\in\cFT\}$. Denote $\mP^k(\cFT):=\times_{F\in\cFT}\mP^k(F)$, $\mP^k(\cF_{\uife})$ and $\mP^k(\cF_{\ife})$ are the restriction of $\mP^k(\cFT)$ on $\uife$ and $\ife$, respectively. We define the local discrete unknowns as
\begin{equation}\label{discrete-variable}
  \updl{v}_T:=(\udl{v}_T, \udl{v}_{\uife}, \udl{v}_{\ife})\in \updl{V}_T
    :=\mP^{k+1}(T)^d \times \mP^k(\cF_{\uife})^d \times \mP^k(\cF_{\ife})^d.
\end{equation}

In order to design a method that is both locking free and robust with respect to the compliancy, some care must be taken when designing the stabilization on the interface.

First we observe that since $\bK$ is a symmetric matrix, there exists ${\bf Q}$ and ${\bf L}$ such that
\begin{align*}
  \bK = {\bf Q} {\bf L}  {\bf Q}^T.
\end{align*}
It is straightforward to verify that if $\udl{n},\udl{t}_1$ and $\udl{t}_2$ form an orthogonal set on the interface, where $\udl{t}_1$ and $\udl{t}_2$ denote the two tangential directions, then
\begin{align*}
  {\bf Q} = [\udl{n} | \udl{t}_1 | \udl{t}_2]\quad\mbox{and}\quad {\bf L} = diag(\beta,\alpha,\alpha).
\end{align*}
Here $diag(x,y,z)$ denotes the diagonal matrix with diagonal elements $x,y$ and $z$. In order to design a method robust with respect to all compliancy parameters $\alpha\geq 0,\, \beta\geq 0$, we perturb the matrix ${\bf L}$ with ${\bf P}$,
\begin{align*}
  \tilde{\bf L} = {\bf L} + {\bf P},\qquad {\bf P} = diag(h_T/(\mu+\lambda),h_T/\mu,h_T/\mu).
\end{align*}
We define the matrix $\bS_h$ (c.f. \cite{HH04}) as follows
\be\label{penalty-matrix}
  \bS_h|_T={\bf Q} \tilde{\bf L}^{-1} {\bf Q}^T, \qquad\forall\ T\in\cTh.
\ee
It is then straightforward to verify that $\udl{n}^T \bS_h \udl{n} = (\beta + h_T/(\mu+\lambda))^{-1}$ and $\udl{t}_i^T \bS_h \udl{t}_i = (\alpha + h_T/\mu)^{-1}$, $i=1,2$, consistently with $\udl{n}^T \bK \udl{n} = \beta$ and $\udl{t}_i^T \bK \udl{t}_i = \alpha$. Furthermore
\begin{align*}
  \bK \bS_h =\bS_h \bK = {\bf Q} {\bf L} \tilde{\bf L}^{-1} {\bf Q}^T,
\end{align*}
and 
\begin{align*}
  \bI - \bK \bS_h = {\bf Q} (\bI-{\bf L} \tilde{\bf L}^{-1}) {\bf Q}^T
    ={\bf Q} {\bf P} \tilde{\bf L}^{-1} {\bf Q}^T
    ={\bf Q} {\bf P} {\bf Q}^T \bS_h.
\end{align*}

The first key ingredient in the devising of the HHO method is a local symmetric strain reconstruction in each mesh cell. For any $T\in\cTh$, we define a local discrete symmetric gradient reconstruction operator $\bE:\ \updl{V}_T\rightarrow \mP^k(T)^{d\times d}$ such that for any $\updl{v}_T=(\udl{v}_T, \udl{v}_{\uife}, \udl{v}_{\ife})$
\begin{align}\label{gradient-reconstruct}
  (\bE(\updl{v}_T), \bt)_T
    &+\frac12\pd{ (\bI-\bK\bS_h)\bK (2\mu \bE +\lambda \bI tr\bE)(\updl{v}_T)\udl{n}, \bt\udl{n} }_{\ife}\nn\\
    =&(\Ve(\udl{v}_T), \bt)_T
    +\pd{ \udl{v}_{\uife}-\udl{v}_T, \bt \udl{n} }_{\uife}\nn\\
    &+\pd{ (\bI-\bK\bS_h)(\udl{v}_{\ife}-\udl{v}_T), \bt\udl{n} }_{\ife}, \quad
      \forall\ \bt\in \mP^k(T)^{d\times d}.
\end{align}

Specially, set $\bt=q \bI$. By the equality $\nabla\cdot \udl{v}=trace\big(\Ve(\udl{v})\big) =\Ve(\udl{v}): \bI$, the local discrete divergence reconstruction operator $D:=tr\bE:\ \updl{V}_T\rightarrow \mP^k(T)$ is defined as follows: For any $\updl{v}_T=(\udl{v}_T, \udl{v}_{\uife}, \udl{v}_{\ife})\in \updl{V}_T$, we have
\begin{align}\label{divergence-reconstruct}
  (tr\bE(\updl{v}_T), q)_T
    &+\frac12\pd{ (\bI-\bK\bS_h)\bK (2\mu \bE +\lambda \bI tr\bE)(\updl{v}_T)\udl{n}, q\udl{n} }_{\ife}\nn\\
    =&(\nabla\cdot\udl{v}_T, q)_T
    +\pd{ \udl{v}_{\uife}-\udl{v}_T, q\udl{n} }_{\uife}\nn\\
    &+\pd{ (\bI-\bK\bS_h)(\udl{v}_{\ife}-\udl{v}_T), q\udl{n} }_{\ife}, \quad
      \forall\ q\in\mP^k(T).
\end{align}

We use the above operators to mimic locally the exact local bilinear form $a(\cdot,\cdot)$ defined
in \reff{weak-form} by means of the following local bilinear form defined on $\updl{V}_T\times \updl{V}_T$
\begin{align}\label{Local-aT}
  a_T(\updl{v}_T,\updl{w}_T)=&2\mu\big(\bE(\updl{v}_T), \bE(\updl{w}_T)\big)_T
    +\lambda\big(D(\updl{v}_T), D(\updl{w}_T)\big)_T\nn\\
    &+\frac12\pd{ (\bI-\bK\bS_h)\bK \bG(\updl{v}_T)\udl{n},\bG(\updl{w}_T)\udl{n} }_{\ife},
\end{align}
where
\begin{equation*}
  \bG(\updl{v}_T) = \mathbb{C}\bE(\updl{v}_T) = (2\mu \bE +\lambda \bI tr\bE)(\updl{v}_T)
    = 2\mu \bE(\updl{v}_T) +\lambda D(\updl{v}_T)\bI.
\end{equation*}

The second key ingredient is the local stabilization operator used to penalize in a least-squares sense the difference between the face unknown $\udl{v}_{\p T}$ and the trace of the cell unknown $\udl{v}_T|_{\p T}$. Let $\Pi_{\p T}$ be the $L^2$-orthogonal projection onto $\mP^k(\cFT)^d$, and $\Pi_F$ will be the $L^2$-orthogonal projection onto $\mP^k(F)^d$ when $F\in\cFT$. For any $\updl{v}_T, \updl{w}_T\in \updl{V}_T$, the stabilization bilinear form is defined as
\begin{align}\label{stabilization-term}
  s_T(\updl{v}_T,\updl{w}_T):=& s_{T^{\backslash \Gamma}}(\updl{v}_T,\updl{w}_T) +s_{T^{\Gamma}}(\updl{v}_T,\updl{w}_T)\nn\\
    :=&\pd{ h_T^{-1} (2\mu\bI +\lambda \udl{n}\otimes\udl{n})
     (\udl{v}_{\uife}-\Pi_{\uife}\udl{v}_{T}), \udl{w}_{\uife}-\Pi_{\uife}\udl{w}_{T} }_{\uife}\nn\\
    &+2\pd{ \bS_h(\udl{v}_{\ife}-\Pi_{\ife}\udl{v}_{T}), \udl{w}_{\ife}-\Pi_{\ife}\udl{w}_{T} }_{\ife}.
\end{align}

Going from local to global bilinear forms proceeds, as in standard finite element methods, by a cell-wise assembly. The global space of DOFs is obtained by patching local DOFs at interfaces, yielding the discrete spaces
\bn
  \updl{V}_h=\big(\times_{T\in\cTh} \mP^{k+1}(T)^d \big)\times \big(\times_{F\in{\cal F}_h}\mP^k(F)^d \big).
\en
We enforce strongly the homogeneous Dirichlet boundary condition on $\p\Omega$ by considering the subspace
\bn
  \updl{V}_{h0}=\big\{ \updl{v}_h =\big( (\udl{v}_T)_{T\in\cTh}, (\udl{v}_F)_{F\in{\cal F}_h} \big) \in\updl{V}_h\big|\ \udl{v}_F=\udl{0},\ \forall \ F\in{\cal F}_h^b\big\}.
\en
The HHO discrete problem is: Find $\updl{u}_h\in \updl{V}_{h0}$ such that
\be\label{discrete-form}
  a_h(\updl{u}_h,\updl{w}_h)=l(\updl{w}_h),\qquad\forall \ \updl{w}_h\in \updl{V}_{h0},
\ee
where
\begin{align}\label{ah}
  a_h(\updl{u}_h,\updl{w}_h)=\sum_{T\in\cTh}  \big( a_T(\updl{u}_T,\updl{w}_T)
    +s_T(\updl{u}_T,\updl{w}_T) \big).
\end{align}


Two comments on the design of the method are in order.  The first one explains
why no unknown is attached to the jump of the displacement, and fixes the
weight of the interface stabilization in \reff{stabilization-term}.  The
second one describes the behaviour of the method in the two extreme regimes of
the compliancy.

\begin{rem}[Interface unknowns and the compliancy term]\label{Rem-interface-unknown}
Since $\cIh\subset{\cal F}_h$, the unknowns attached to the interface faces are
single valued, and they approximate the average $\av{\udl{u}}$ rather than
either of the two traces; the jump $\jm{\udl{u}}$ is not an unknown of the
method and is recovered a posteriori as the difference of the traces of the
two adjacent cell unknowns.  The discrete counterpart of the compliancy term
$\pd{\bS\jm{\udl{u}},\jm{\udl{v}}}_\Gamma$ of \reff{weak-form} is produced by
the interface stabilization $s_{T^\Gamma}$, as the following computation
shows.  Let $F=T_1\cap T_2\in\cIh$ and assume for simplicity that
$h_{T_1}=h_{T_2}$ and that the two cells carry the same matrix $\bS_h$.
Eliminating the sole interface unknown $\udl{v}_F$ from
$s_{T_1^\Gamma}+s_{T_2^\Gamma}$, i.e.\ minimising
\begin{align*}
  2\pd{ \bS_h(\udl{v}_F-\Pi_F\udl{v}_{T_1}), \udl{v}_F-\Pi_F\udl{v}_{T_1} }_F
  +2\pd{ \bS_h(\udl{v}_F-\Pi_F\udl{v}_{T_2}), \udl{v}_F-\Pi_F\udl{v}_{T_2} }_F
\end{align*}
over $\udl{v}_F\in\mP^k(F)^d$, gives
$\udl{v}_F=\frac12\Pi_F(\udl{v}_{T_1}+\udl{v}_{T_2})$ and the value
$\pd{\bS_h\Pi_F\jm{\udl{v}},\Pi_F\jm{\udl{v}}}_F$.  The factor $2$ in
\reff{stabilization-term} is thus precisely the one for which the elimination
of the interface unknown reproduces the compliancy term of \reff{weak-form},
with $\bS$ replaced by its regularisation $\bS_h$.
\end{rem}

\begin{rem}[The two extreme regimes]\label{Rem-limits}
The definition \reff{penalty-matrix} of $\bS_h$ is designed so that one single
formulation covers the whole range of compliancies.  When $\alpha=\beta=0$,
i.e.\ for a perfectly bonded interface, we have $\bK=\bf{0}$ and therefore
$\bI-\bK\bS_h=\bI$ and
$\bS_h={\bf Q}{\bf P}^{-1}{\bf Q}^T =\frac{\mu+\lambda}{h_T}\udl{n}\otimes\udl{n}
+\frac{\mu}{h_T}\sum_{i=1,2}\udl{t}_i\otimes\udl{t}_i$; the interface faces
are then treated as ordinary faces and $s_{T^\Gamma}$ is a Nitsche-type
penalty enforcing the continuity of the displacement, the weights being the
same as those of $s_{T^{\backslash\Gamma}}$ up to the factor $2$ discussed in
\Rem{\ref{Rem-interface-unknown}}.  In the opposite regime
$\alpha,\beta\gg h_T/\mu$, we have $\bS_h\simeq\bK^{-1}$, so that
$s_{T^\Gamma}$ represents the physical interface stiffness, and
$\bI-\bK\bS_h\simeq\bf{0}$, so that the interface contributions to
\reff{gradient-reconstruct} and \reff{Local-aT} disappear; in the limit
$\alpha,\beta\to+\infty$ the interface unknowns become inactive and the
discrete problem decouples into two independent problems with a traction-free
boundary condition on $\Gamma$, which is the correct behaviour of
\reff{elasticity-interface} in that limit.  The perturbation $\bf P$ is what
makes both limits, and the degenerate cases $\alpha=0$ or $\beta=0$, uniformly
accessible; this is confirmed numerically in \S\ref{sec-numerics}.
\end{rem}

\Section{Some preliminaries.}

We have some properties of the penalty matrix $\bS_h$ (c.f. Lemma 2 in \cite{HH04}).
\begin{lem}\label{Lem-Sh} Define $\bS_h$ as in \reff{penalty-matrix}, and let $|\cdot|$ denote the matrix norm induced by the Euclidean norm on $\mR^d$. Set $\bS_{h,n}=\bS_h \udl{n}\otimes\udl{n}$ and $\bS_{h,t}=\bS_h (\bI -\udl{n}\otimes\udl{n})$, then
\be\label{Lem-Sh-1}
  \bI-\bK\bS_h={\bf Q}{\bf P}{\bf Q}^T\bS_h,
\ee
\be\label{Lem-Sh-2}
  \udl{n}^T\bS_h\udl{n}\leq \frac{\mu+\lambda}{h}, \quad \udl{t}_i^T\bS_h\udl{t}_i\leq \frac{\mu}{h},i=1,2,\quad
  \mbox{and}\ \abs{\bS_{h,t}}\leq \frac{\mu}{h},
\ee
\be\label{Lem-Sh-3}
  \abs{\bI-\bK\bS_h}\leq 1, \qquad \abs{\bS_h\bK}\leq 1,
\ee
\be\label{Lem-Sh-4}
  \udl{n}^T{\bf P}\udl{n}= \frac{h}{\mu+\lambda}, \quad 
  \udl{t}_i^T{\bf P}\udl{t}_i= \frac{h}{\mu},i=1,2, \quad 
  \mbox{and}\ \abs{(\bI-\bK\bS_h)\bK}\leq \frac{h}{\mu}.
\ee
\end{lem}

In this section we recall some basic results on admissible mesh sequences (c.f. Lemma 1.46 and 1.49 in \cite{DE12}, Lemma B.66 and 3.78 in \cite{EG04}).
\begin{lem}\label{Lem-Basic-Ineq}
There are $C_{dtr}>0$ and $C_{mtr}>0$, such that the following inequalities hold true:
\begin{equation}\label{Lem-Discrete-trace-Ineq}
  \norm{v}_F\leq C_{dtr} h_T^{-\frac12}\norm{v}_T,\quad
    \forall \, v\in \mP^l(T), \, \forall\, F\in\cFT,\, T\in\cTh,
\end{equation}
and
\begin{equation}\label{Lem-Multip-trace-Ineq}
  \norm{v}_{\p T}\leq  C_{mtr}\big( h_T^{-\frac12}\norm{v}_T +h_T^{\frac12}\norm{\nabla v}_T \big), \quad
    \forall\, v\in H^1(T),\, T\in\cTh.
\end{equation}
For any $T\in\cTh$, the following Poincar\'e inequality is valid
\begin{equation}\label{Lem-Poincare-Ineq}
  \norm{v}_T\leq  C_Ph_T \norm{\nabla v}_T, \quad
    \forall\, v\in H^1(T),\, \int_T v=0.
\end{equation}
Let $RM=\{\udl{v} \in \mP^1(T)^d |\ \Ve (\udl{v}) ={\bf 0}\}$, the following second Korn's inequality holds for any $T\in\cTh, \, \udl{v}\in H^1(T)^d$ satisfying $(\udl{v},\udl{r})_T=0$, $\forall\, \udl{r}\in RM$,
\begin{equation}\label{Lem-Korn-Ineq}
  \norm{\nabla \udl{v}}_T\leq  C_K \norm{\Ve (\udl{v})}_T.
\end{equation}
\end{lem}

In what follows, we use the convention $A\lesssim B$ to abbreviate the inequality $A\leq CB$ for positive real numbers $A$ and $B$, where the constant $C$ only depends on the polynomial degree $k\geq 1$ used in the HHO method, the mesh parameters $\rho$ and the above constants $C_{dtr}$, $C_P$ and $C_K$, but does not depend neither on the Lam\'e coefficients $\lambda_i,\mu_i,i=1,2$ nor on the mesh size $h>0$.

We equip the space $\updl{V}_T$ with the following local energy norm:
\begin{align}\label{Energy-norm-local}
  \norme{\updl{v}_T}_E^2= 2\mu\norm{\Ve(\udl{v}_T)}^2_T
    +\frac{2\mu}{h_T}\norm{\udl{v}_{\uife}-\udl{v}_T}^2_{\uife}
    +\norm{\bS_h^\frac12(\udl{v}_{\ife}-\udl{v}_T)}^2_{\ife}.
\end{align}
Observe that the three contributions to \reff{Energy-norm-local} are the ones
that survive in the two extreme regimes discussed in \Rem{\ref{Rem-limits}}:
the last one behaves like $\mu h_T^{-1}\norm{\udl{v}_{\ife}-\udl{v}_T}^2_{\ife}$
when $\alpha=\beta=0$, in which case the interface faces are treated as
ordinary faces, and it degenerates to
$\norm{\bK^{-\frac12}(\udl{v}_{\ife}-\udl{v}_T)}^2_{\ife}$ when the compliancy
dominates, in which case it mimics the compliancy term of \reff{weak-form}.

The following elementary estimate will be used twice; it expresses that the
distance between a cell polynomial and the $L^2$-projection of its trace onto
the face polynomial spaces is controlled by its symmetric gradient alone.
Note that no $\lambda$ appears in it.
\begin{lem}\label{Lem-Fce-Korn} Let $k\geq1$. For any $T\in\cTh$ and any
$\udl{v}_T\in\mP^{k+1}(T)^d$, we have
\begin{equation}\label{Lem-Fce-Korn-0}
    h_T^{-\frac12}\norm{\udl{v}_T -\Pi_{\p T}\udl{v}_T}_{\p T}
    \leq C_{dtr}C_PC_K \norm{\Ve (\udl{v}_T)}_T .
\end{equation}
\end{lem}
\Proof Let $\Pi_{RM}$ denote the $L^2(T)$-orthogonal projection onto the space
$RM$ of rigid-body motions defined in \Lem{\ref{Lem-Basic-Ineq}}.  Since
$RM\subset \mP^1(T)^d$ and $k\geq1$, the trace on any face $F\in\cFT$ of a
rigid-body motion belongs to $\mP^k(F)^d$ so that
$\Pi_{\p T}\Pi_{RM}(\udl{v}_T)=\Pi_{RM}(\udl{v}_T)$ and the left-hand side of
\reff{Lem-Fce-Korn-0} is unchanged if $\udl{v}_T$ is replaced by
$\udl{v}_T-\Pi_{RM}(\udl{v}_T)$.  Setting
$\udl{\psi}:=\udl{v}_T-\Pi_{RM}(\udl{v}_T)$ and using that $\Pi_{\p T}$ is an
$L^2(\p T)$-orthogonal projection, we get
\begin{align*}
  h_T^{-1}\norm{\udl{v}_T -\Pi_{\p T}\udl{v}_T }^2_{\p T}
    = h_T^{-1}\norm{(I-\Pi_{\p T} )\udl{\psi} }^2_{\p T}
    \leq h_T^{-1}\norm{\udl{\psi} }^2_{\p T}.
\end{align*}
The discrete trace inequality \reff{Lem-Discrete-trace-Ineq}, the
Poincar\'e inequality \reff{Lem-Poincare-Ineq} applied component-wise, which
is licit since $\int_T\udl{\psi}=\udl{0}$, and Korn's inequality
\reff{Lem-Korn-Ineq}, which is licit since $(\udl{\psi},\udl{r})_T=0$ for all
$\udl{r}\in RM$, then give in turn
\begin{align*}
  h_T^{-1}\norm{\udl{\psi} }^2_{\p T}
    &\leq C_{dtr}^2 h_T^{-2}\norm{\udl{\psi} }^2_T,\\
  h_T^{-2}\norm{\udl{\psi} }^2_T
    &\leq C_P^2 \norm{\nabla\udl{\psi} }^2_T,\\
  \norm{\nabla\udl{\psi} }^2_T
    &\leq C_K^2\norm{\Ve (\udl{\psi}) }^2_T .
\end{align*}
Finally, every element of $RM$ has a vanishing symmetric gradient, so that
$\Ve (\udl{\psi})=\Ve (\udl{v}_T)$, and collecting the three inequalities
above yields \reff{Lem-Fce-Korn-0}. \Endproof

The following local stability property of the reconstruction operator and
stabilization terms is established similar to Lemma 4 of \cite{DE15}.
\begin{lem}\label{Lem-equivalent-norm} Let $\bE$ be defined by \reff{gradient-reconstruct}. There is $0<\alpha_{\flat}<+\infty$ such that, for any $T\in\cTh$, any $\updl{v}_T\in\updl{V}_T$, we have
\begin{align}\label{Lem-equivalent-norm-0}
  \alpha_{\flat}\norme{\updl{v}_T}_E \leq& \bigg( 2\mu\norm{\bE(\updl{v}_T)}^2_T
    +\norm{[(\bI-\bK\bS_h)\bK]^\frac12 \bG(\updl{v}_T)\udl{n}}^2_{\ife} \nn\\
      &+\frac{2\mu}{h_T}\norm{\udl{v}_{\uife} -\Pi_{\uife}\udl{v}_T}^2_{\uife}
      +s_{T^{\Gamma}}(\updl{v}_T,\updl{v}_T) \bigg)^{\frac12}.
\end{align}
\end{lem}
\Proof The right-hand side of \reff{Lem-equivalent-norm-0} is, up to the
factor $\lambda\norm{D(\updl{v}_T)}^2_T$ and up to the $\lambda$-weighted part
of $s_{T^{\backslash\Gamma}}$, the square of $\norme{\updl{v}_T}_{a_T}$; the point
of the lemma is therefore that the first contribution
$2\mu\norm{\Ve(\udl{v}_T)}^2_T$ to the energy norm
\reff{Energy-norm-local}, which involves the {\em exact} symmetric gradient of
the cell unknown, is controlled by the {\em reconstructed} one.  The proof
proceeds in two steps: we first bound $\norm{\sqrt{2\mu}\Ve(\udl{v}_T)}_T$,
and we then bound the two remaining contributions to
\reff{Energy-norm-local}.

\smallskip\noindent{\em Step 1: bound on $\norm{\sqrt{2\mu}\Ve(\udl{v}_T)}_T$.}
We test the definition \reff{gradient-reconstruct} of $\bE$ with the admissible
choice $\bt=2\mu\Ve(\udl{v}_T)\in\mP^{k}(T)^{d\times d}$, and we move the
interface term of the left-hand side to the right-hand side.  Since
$\big( \Ve(\udl{v}_T) ,2\mu\Ve(\udl{v}_T) \big)_T
=2\mu\norm{\Ve(\udl{v}_T)}^2_T$, this yields the identity
\begin{align}\label{Lem-equivalent-norm-1}
  2\mu\norm{\Ve(\udl{v}_T)}^2_T
    =& \big( \bE(\updl{v}_T) ,2\mu\Ve(\udl{v}_T) \big)_T
      +\frac12 \pd{ (\bI-\bK\bS_h)\bK \bG(\updl{v}_T) \udl{n}, 2\mu\Ve(\udl{v}_T) \udl{n} }_{\ife} \nn\\
    &-\pd{ \udl{v}_{\uife} -\Pi_{\uife}\udl{v}_T, 2\mu\Ve(\udl{v}_T) \udl{n} }_{\uife}  \nn\\
    &-\pd{ (\bI-\bK\bS_h) (\udl{v}_{\ife} -\Pi_{\ife}\udl{v}_T), 2\mu\Ve(\udl{v}_T) \udl{n} }_{\ife},
\end{align}
where we have also used that $2\mu\Ve(\udl{v}_T)\udl{n}$ has its components in
$\mP^{k}(F)^d$ on every face $F\in\cFT$, so that the traces
$\udl{v}_{\uife}-\udl{v}_T$ and $\udl{v}_{\ife}-\udl{v}_T$ appearing in
\reff{gradient-reconstruct} may be replaced by
$\udl{v}_{\uife}-\Pi_{\uife}\udl{v}_T$ and
$\udl{v}_{\ife}-\Pi_{\ife}\udl{v}_T$, respectively.

We bound the four terms of \reff{Lem-equivalent-norm-1} separately.  The first
one is handled by the Cauchy--Schwarz inequality.  For the second one we use
that $(\bI-\bK\bS_h)\bK$ is symmetric positive semi-definite together with the
discrete trace inequality \reff{Lem-Discrete-trace-Ineq}, which gives
\begin{align*}
  \frac12&\pd{ (\bI-\bK\bS_h)\bK \bG(\updl{v}_T) \udl{n}, 2\mu\Ve(\udl{v}_T) \udl{n} }_{\ife}\\
  &\leq \norm{\big[ (\bI-\bK\bS_h)\bK \big]^\frac12\bG(\updl{v}_T)\udl{n}}_{\ife}
        \norm{\big[ (\bI-\bK\bS_h)\bK \big]^\frac12 \mu \Ve(\udl{v}_T)\udl{n}}_{\ife}\\
  &\leq C_{dtr} \norm{\big[ (\bI-\bK\bS_h)\bK \big]^\frac12\bG(\updl{v}_T)\udl{n}}_{\ife}
        \norm{\sqrt{2\mu}\Ve(\udl{v}_T)}_{T},
\end{align*}
where the last inequality follows from
$\abs{(\bI-\bK\bS_h)\bK}\leq h_T/\mu$, see \reff{Lem-Sh-4}.  The third term is
bounded in the same way, using $\abs{2\mu\Ve(\udl{v}_T)\udl{n}}\leq
2\mu\abs{\Ve(\udl{v}_T)}$ and again \reff{Lem-Discrete-trace-Ineq}.  For the
fourth term we use the factorisation
\begin{align}\label{Lem-equivalent-norm-2}
  \abs{\bI-\bK\bS_h}
    =\abs{ ({\bf Q} {\bf P} {\bf Q}^T)^\frac12\bS_h^\frac12
      [\bI-\bK\bS_h]^\frac12}
    \leq \Big(\frac{h_T}{\mu}\Big)^\frac12\abs{\bS_h^\frac12},
\end{align}
which follows from \reff{Lem-Sh-1} and \reff{Lem-Sh-3}--\reff{Lem-Sh-4}.
Collecting the four bounds and dividing by
$\norm{\sqrt{2\mu}\Ve(\udl{v}_T)}_T$, we obtain
\begin{align}\label{Lem-equivalent-norm-3}
  \norm{\sqrt{2\mu}\Ve(\udl{v}_T)}_T
  \leq& \sqrt{2\mu}\norm{\bE(\updl{v}_T)}_T
    +C_{dtr} \norm{\big[ (\bI-\bK\bS_h)\bK \big]^\frac12\bG(\updl{v}_T)\udl{n}}_{\ife} \nn\\
    &+C_{dtr} \Big( \sqrt{\frac{2\mu}{h_T}}\norm{\udl{v}_{\uife} -\Pi_{\uife}\udl{v}_T}_{\uife}
     +s_{T^{\Gamma}}(\updl{v}_T,\updl{v}_T)^\frac12 \Big),
\end{align}
which is the announced bound on the first contribution to
\reff{Energy-norm-local}.

\smallskip\noindent{\em Step 2: bound on the face contributions.}
It remains to bound $h_T^{-1}\norm{\udl{v}_{\p T} -\udl{v}_T}^2_{\p T}$, in
which the {\em trace} of the cell unknown appears, whereas the right-hand side
of \reff{Lem-equivalent-norm-0} only involves its $L^2$-projection onto the
face spaces.  Writing
$\udl{v}_{\p T} -\udl{v}_T=\big(\udl{v}_{\p T} -\Pi_{\p T}\udl{v}_T\big)
-\big(I-\Pi_{\p T}\big)\udl{v}_T$ and recalling the definition
\reff{stabilization-term} of the stabilization bilinear form, we only need to
bound $h_T^{-1}\norm{(I-\Pi_{\p T})\udl{v}_T}^2_{\p T}$, and this is precisely
the content of \Lem{\ref{Lem-Fce-Korn}}:
\begin{align*}
    h_T^{-1}\norm{(I-\Pi_{\p T})\udl{v}_T }^2_{\p T}
    \leq C_{dtr}^2C_P^2C_K^2\norm{\Ve (\udl{v}_T) }^2_T .
\end{align*}
Together with \reff{Lem-equivalent-norm-3} and with the bound
$\abs{\bS_h^{\frac12}}^2\leq\mu/h_T$ of \reff{Lem-Sh-2}, which allows us to
compare the $\bS_h$-weighted face term of \reff{Energy-norm-local} with
$s_{T^{\Gamma}}$, this gives the inequality
\reff{Lem-equivalent-norm-0}.\Endproof

The key operator in the HHO error analysis is the local interpolation operator $\updl{I}_T:\ H^1(T)^d\rightarrow \updl{V}_T$. For any $\udl{v}\in H^1(\Omega_1\cup \Omega_2)^d$, we define
\be\label{reduction-map}
  \updl{I}_T(\udl{v})=\big( \Pi_T(\udl{v}), \Pi_{\uife}(\udl{v}), \Pi_{\ife}(\udl{v}) \big)\in \updl{V}_T,\quad
    \forall\ T\in\cTh.
\ee
Here $\Pi_T$ denotes the $L^2$-orthogonal projection onto $\mP^{k+1}(T)^d$ and $\Pi_F$ is the $L^2$-orthogonal projection onto $\mP^k(F)^d$ defined in \reff{stabilization-term} for any $F\in\cFT$, when $\updl{I}_T(\udl{v})$ is used for the reconstruction operator $\bE$. While in the definition of $D\updl{I}_T(\udl{v})$, $\Pi_T$ and $\Pi_{\p T}$ are the $H(div)$-type BDM interpolation operator defined onto $BDM(k+1)$ (c.f. section III.3.3 in \cite{BF91}). We also define $\updl{I}_h(\udl{v})\in\updl{V}_h$ such that, for any $T\in\cTh$, the local components of $\updl{I}_h(\udl{v})$ in $T$ are $\updl{I}_T(\udl{v})\in\updl{V}_T$. Note that $\updl{I}_h(\udl{v})\in\updl{V}_{h0}$ whenever $\udl{v}\in H_0^1(\Omega_1\cup \Omega_2)^d:=\set{\udl{v}\in H^1(\Omega_1\cup \Omega_2)^d\big|\ \udl{v}|_{\p\Omega}=\udl{0}}$.

\begin{rem}[The two interpolation operators]\label{Rem-interpolants}
The reason for using two different sets of projections in \reff{reduction-map}
is the following. The $L^2$-orthogonal projections are the ones for which the
reconstruction $\bE$ is consistent, in the sense that
$\bE(\updl{I}_T\udl{v})$ approximates $\pk\Ve(\udl{v})$ optimally, where $\pk$ is the elementwise projection on polynomials of order $k$; this is the
content of \Thm{\ref{Thm-reconstruct-Approx}} below.  They do not, however,
preserve normal traces, and the resulting bound on the {\em divergence} part
$D(\updl{I}_T\udl{v})-\pk(\nabla\cdot\udl{v})$ would then carry a factor
$\lambda$ instead of $\mu$, which would spoil the robustness in the
quasi-incompressible limit.  The BDM interpolation operator, on the contrary,
satisfies $\big(\udl{v}-\Pi_T\udl{v}\big)\cdot\udl{n}=0$ on every face of $T$
and commutes with the divergence in the sense that
$\nabla\cdot\Pi_T\udl{v}=\pk(\nabla\cdot\udl{v})$; the terms in which the
normal component of the interpolation error appears therefore drop out.  Both
operators enjoy the approximation properties of
\Lem{\ref{Lem-Interpolation-properties}}, which is all that is used below, and
we keep the same symbol $\updl{I}_T$ for both, indicating explicitly which one
is meant whenever the distinction matters.  This device is the same as the one
introduced in \cite{DE15} to obtain locking-free estimates; see also
\Rem{\ref{Rem-lambda-norm}}.
\end{rem}

\begin{lem}(Local approximation c.f. Proposition 1.134 in \cite{EG04})\label{Lem-Interpolation-properties} Let $\udl{v}\in H^1(\Omega_1\cup \Omega_2)^d$. We assume that $\udl{v}|_{\Omega_i}\in H^{m+1}(\Omega_i)^d$ for some $m\geq 0,i=1,2$. Set $l=\min\set{k+1,m}$. For any $T\in\cTh$, we have
\begin{align}\label{Lem-Interpolation-properties-0}
  \norm{\udl{v}-\Pi_T(\udl{v})}_T
    &+h_T^{\frac12}\norm{\udl{v}-\Pi_T(\udl{v})}_{\p T}
     +h_T\abs{\udl{v}-\Pi_T(\udl{v})}_{H^1(T)^d}\nn\\
     &+h_T^{\frac{3}{2}}\abs{\udl{v}-\Pi_T(\udl{v})}_{H^1(\p T)^d}
     \lesssim h_T^{l+1}\abs{\udl{v}}_{H^{l+1}(T)^d}.
\end{align}
\end{lem}

\Section{Stability and error analysis.}

For any $\updl{v}_h\in\updl{V}_h$, define the energy norms
\be\label{Energy-norm}
  \norme{\updl{v}_h}^2_e:=\sum_{T\in\cTh}\norme{\updl{v}_T}^2_E,
\ee
and
\begin{align}\label{Energy-norm-h}
  \norme{\updl{v}_h}^2_h:=&a_h(\updl{v}_h,\updl{v}_h) \nn\\
    =&\sum_{T\in\cTh}\norme{\updl{v}_T}^2_{a_T},
\end{align}
with
\begin{align}\label{Energy-norm-T}
  \norme{\updl{v}_T}_{a_T}^2:=&a_T(\updl{v}_T,\updl{v}_T) +s_T(\updl{v}_T,\updl{v}_T) \nn\\
        =&2\mu\norm{\bE(\updl{v}_T)}^2_T
    +\lambda\norm{D(\updl{v}_T)}^2_T
    +\frac{2\mu}{h_T}\norm{\udl{v}_{\uife}-\Pi_{\uife}\udl{v}_T}^2_{\uife} \nn\\
    &+\frac{\lambda}{h_T}\norm{(\udl{v}_{\uife}-\Pi_{\uife}\udl{v}_T)\cdot\udl{n}}^2_{\uife}
    +2\norm{\bS_h^\frac12(\udl{v}_{\ife}-\Pi_{\ife}\udl{v}_T)}^2_{\ife} \nn\\
    &+\frac12\norm{[(\bI-\bK\bS_h)\bK]^\frac12 \bG(\updl{v}_T)\udl{n}}^2_{\ife}.
\end{align}
The six contributions to \reff{Energy-norm-T} are, in this order, the three
terms of the local bilinear form $a_T$ of \reff{Local-aT} and the three
terms of the stabilization $s_T$ of \reff{stabilization-term}.  We can prove
that $\norme{\cdot}_e$ is a norm on $\updl{V}_{h0}$ using the discrete Korn's
inequality, and we get the stability of the HHO method from
\Lem{\ref{Lem-equivalent-norm}}.
\begin{thm}\label{Thm-Stability} The discrete form $a_h(\cdot,\cdot)$ is continuous, i.e.
\be\label{Thm-Stability-continuity}
  a_h(\updl{v}_h,\updl{w}_h)\leq \norme{\updl{v}_h}_h\norme{\updl{w}_h}_h,\qquad\forall \ \updl{v}_h,\updl{w}_h\in\updl{V}_h.
\ee
Furthermore, it is coercive on $\updl{V}_h$, i.e.
\be\label{Thm-Stability-coercivity}
  a_h(\updl{v}_h,\updl{v}_h) = \norme{\updl{v}_h}^2_h
    \geq \alpha^2_\flat \norme{\updl{v}_h}^2_e, \qquad\forall \ \updl{v}_h\in\updl{V}_h.
\ee
\end{thm}
\Proof The bilinear form $a_h$ is symmetric and, by \reff{Energy-norm-h},
positive semi-definite, so that \reff{Thm-Stability-continuity} is nothing but
the Cauchy--Schwarz inequality for $a_h$.  The coercivity
\reff{Thm-Stability-coercivity} follows by summing the local bound
\reff{Lem-equivalent-norm-0} of \Lem{\ref{Lem-equivalent-norm}} over
$T\in\cTh$ and comparing the right-hand side with \reff{Energy-norm-T}, which
is licit since each of the four terms in \reff{Lem-equivalent-norm-0} is one
of the six terms of \reff{Energy-norm-T}, up to the factors $2$ and
$\frac12$.\Endproof

The next result is the cornerstone of the error analysis: it states that,
applied to the interpolate of a smooth function, the reconstruction operator
$\bE$ is an optimal approximation of the symmetric gradient, and that the same
holds for the discrete divergence and for the reconstructed traction on the
interface.  In the last term of \reff{Thm-reconstruct-Approx-0}, and only
there, $\Pi_{\ife}$ denotes the $L^2(\ife)$-orthogonal projection onto
$\mP^k(\ife)^{d\times d}$; see \Rem{\ref{Rem-projections}} below for the
variant in which the cell projection is used instead.
\begin{thm}\label{Thm-reconstruct-Approx} Let $\udl{u}$ be the solution of \reff{elasticity-interface} and assume $\udl{u}\in H^{k+2}(\Omega_1\cup \Omega_2)^d$. For any $T\in\cTh$, we have
\begin{align}\label{Thm-reconstruct-Approx-0}
  \sqrt{2\mu}\norm{\bE\big( \updl{I}_T(\udl{u}) \big)-\pk\Ve(\udl{u})}_T
    &+\sqrt{\lambda}\norm{D\big( \updl{I}_T(\udl{u}) \big)-\pk(\nabla\cdot\udl{u})}_T \nn\\
    &+\norm{[(\bI-\bK\bS_h)\bK]^\frac12 \big( \bG\big( \updl{I}_T(\udl{u}) \big) -\Pi_{\ife}\Vs(\udl{u}) \big)\udl{n}}_{\ife}\nn\\
   \lesssim& \sqrt{\mu} h_T^{k+1}\abs{\udl{u}}_{H^{k+2}(T)^d},
\end{align}
where $\pk$ is the $L^2(T)$-orthogonal projection on the polynomial function space of degree $k$.
\end{thm}
\Proof Throughout the proof, $\bt$ is an arbitrary element of
$\mP^k(T)^{d\times d}$ and we abbreviate $\updl{I}_T\udl{u}$ for
$\updl{I}_T(\udl{u})$.  Since $\udl{u}$ is discontinuous across $\Gamma$, we
set $\Pi_{\ife}\udl{u}:=\Pi_{\ife}\av{\udl{u}}$, i.e.\ the interface
component of the interpolate approximates the average of the two traces, and
we shall repeatedly use the two interface conditions of
\reff{elasticity-interface} in the form
\be\label{Thm-reconstruct-Approx-jump}
  \av{\udl{u}}=\udl{u}-\frac12\jm{\udl{u}},\qquad
  \jm{\udl{u}}=-\bK\Vs(\udl{u})\udl{n} \qquad \mbox{on }\ \Gamma.
\ee
The proof consists in identifying the residual of the defining relation
\reff{gradient-reconstruct} when $\updl{v}_T=\updl{I}_T\udl{u}$ (Steps 1--3),
and in choosing the test function $\bt$ so that this residual is the quantity
to be estimated (Step 4).

\smallskip\noindent{\em Step 1: integration by parts.}
By the definition \reff{gradient-reconstruct} of the reconstruction operator,
\begin{align}\label{Thm-reconstruct-Approx-2}
  (\bE(\updl{I}_T\udl{u}), \bt)_T &+\frac12\pd{(\bI-\bK\bS_h)\bK\bG(\updl{I}_T\udl{u})\udl{n} ,\bt\udl{n}}_{\ife} \nn\\
    =& (\Ve(\Pi_T\udl{u}), \bt)_T +\pd{\Pi_{\uife}\udl{u}-\Pi_T\udl{u}, \bt\udl{n}}_{\uife}
       +\pd{(\bI-\bK\bS_h)(\Pi_{\ife}\av{\udl{u}}-\Pi_T\udl{u}), \bt\udl{n}}_{\ife}.
\end{align}
We integrate by parts the volume term of the right-hand side, which produces
the boundary contribution $\pd{\Pi_T\udl{u},\bt\udl{n}}_{\p T}$; the latter
cancels the two terms $-\pd{\Pi_T\udl{u},\bt\udl{n}}$ on $\uife$ and on $\ife$
up to the factor $\bI-\bK\bS_h$ in the second one.  Since
$\bt\udl{n}|_F\in\mP^k(F)^d$ for every $F\in\cFT$, the projections
$\Pi_{\uife}$ and $\Pi_{\ife}$ may be removed from the first argument of the
remaining boundary pairings, and we obtain
\begin{align}\label{Thm-reconstruct-Approx-3}
  (\bE(\updl{I}_T\udl{u}), \bt)_T &+\frac12\pd{(\bI-\bK\bS_h)\bK\bG(\updl{I}_T\udl{u})\udl{n} ,\bt\udl{n}}_{\ife} \nn\\
    =& -(\Pi_T\udl{u}, \nabla\cdot\bt)_T +\pd{\udl{u}, \bt\udl{n}}_{\uife}
       +\pd{\av{\udl{u}}, \bt\udl{n}}_{\ife}
      - \pd{\bK\bS_h (\Pi_{\ife}\av{\udl{u}}-\Pi_T\udl{u}), \bt\udl{n}}_{\ife}.
\end{align}

\smallskip\noindent{\em Step 2: back to the exact strain.}
In \reff{Thm-reconstruct-Approx-3} we may replace $\Pi_T\udl{u}$ by $\udl{u}$
in the volume term, because $\nabla\cdot\bt\in\mP^{k-1}(T)^d\subset
\mP^{k+1}(T)^d$, and integrate by parts backwards.  Using the first identity in
\reff{Thm-reconstruct-Approx-jump} to write
$\av{\udl{u}}=\udl{u}-\frac12\jm{\udl{u}}$ in the interface pairing, and
removing once more the projection $\Pi_{\ife}$ in the last term, this yields
\begin{align}\label{Thm-reconstruct-Approx-4}
  (\bE(\updl{I}_T\udl{u}), \bt)_T &+\frac12\pd{(\bI-\bK\bS_h)\bK\bG(\updl{I}_T\udl{u})\udl{n} ,\bt\udl{n}}_{\ife} \nn\\
    =& (\Ve(\udl{u}), \bt)_T -\frac12\pd{\jm{\udl{u}}, \bt\udl{n}}_{\ife}
       -\pd{\bK\bS_h (\av{\udl{u}}-\Pi_T\udl{u}), \bt\udl{n}}_{\ife}.
\end{align}
The first term of the right-hand side is $(\pk\Ve(\udl{u}),\bt)_T$, since
$\bt\in\mP^k(T)^{d\times d}$.  In the last one we insert again
$\av{\udl{u}}-\Pi_T\udl{u}=(\udl{u}-\Pi_T\udl{u})-\frac12\jm{\udl{u}}$, so
that the two terms involving the jump combine into
$-\frac12\pd{(\bI-\bK\bS_h)\jm{\udl{u}},\bt\udl{n}}_{\ife}$, and we use the
second identity in \reff{Thm-reconstruct-Approx-jump} to replace
$-\jm{\udl{u}}$ by $\bK\Vs(\udl{u})\udl{n}$.  This gives
\begin{align}\label{Thm-reconstruct-Approx-5}
  (\bE(\updl{I}_T\udl{u}), \bt)_T &+\frac12\pd{(\bI-\bK\bS_h)\bK\bG(\updl{I}_T\udl{u})\udl{n} ,\bt\udl{n}}_{\ife} \nn\\
    =& (\pk\Ve(\udl{u}), \bt)_T +\frac12\pd{(\bI-\bK\bS_h)\bK\Vs(\udl{u})\udl{n}, \bt\udl{n}}_{\ife} \nn\\
    &-\pd{\bK\bS_h (\udl{u}-\Pi_T\udl{u}), \bt\udl{n}}_{\ife}.
\end{align}

\smallskip\noindent{\em Step 3: the residual.}
Two simplifications are now available in \reff{Thm-reconstruct-Approx-5}.
First, in the interface pairing we may replace $\Vs(\udl{u})$ by
$\Pi_{\ife}\Vs(\udl{u})$: indeed, for any symmetric tensor field $\bt$ with
$\bt\udl{n}|_F\in\mP^k(F)^d$ one has
$\pd{\big(\Vs(\udl{u})-\Pi_{\ife}\Vs(\udl{u})\big)\udl{n},\bt\udl{n}}_{\ife}=0$,
because $\mbox{sym}\big((\bt\udl{n})\otimes\udl{n}\big)$ belongs to
$\mP^k(\ife)^{d\times d}$ and $\udl{n}$ is constant on each interface face.
Second, we decompose $\bK\bS_h(\udl{u}-\Pi_T\udl{u})$ along the normal and the
two tangential directions,
\begin{align*}
  \bK\bS_h (\udl{u}-\Pi_T\udl{u})
    = \bK\bS_h\udl{n}\, (\udl{u}-\Pi_T\udl{u})\cdot\udl{n}
      +\sum_{i=1,2} \bK\bS_h\udl{t}_i\, (\udl{u}-\Pi_T\udl{u})\cdot\udl{t}_i,
\end{align*}
and we observe that the normal contribution vanishes.  Indeed, as recalled in
\Rem{\ref{Rem-interpolants}}, the normal component of the cell interpolate is
the one of the BDM interpolation operator, for which
$(\udl{u}-\Pi_T\udl{u})\cdot\udl{n}=0$ on every face of $T$.  Collecting these
two observations, and writing the tangential contribution in a form that is
adapted to the interface term of the left-hand side by means of the
factorisation $\bI-\bK\bS_h = {\bf QPQ}^T\bS_h$ of \reff{Lem-Sh-1}, we arrive
at the residual identity
\begin{align}\label{Thm-reconstruct-Approx-6}
  &(\bE(\updl{I}_T\udl{u})-\pk\Ve(\udl{u}), \bt)_T
     +\frac12\pd{(\bI-\bK\bS_h)\bK\big(\bG(\updl{I}_T\udl{u})-\Pi_{\ife}\Vs(\udl{u})\big)\udl{n} ,\bt\udl{n}}_{\ife} \nn\\
    &\hspace{1cm}= -\sum_{i=1,2}\pd{\udl{t}_i^T({\bf QPQ}^T)^{-\frac12} \bK^\frac12\bS_h^\frac12\udl{t}_i\,
        (\udl{u}-\Pi_T\udl{u}),
        (\bI-\bK\bS_h)^{\frac12} \bK^\frac12\bt\udl{n}}_{\ife},
\end{align}
valid for every $\bt\in\mP^k(T)^{d\times d}$.

\smallskip\noindent{\em Step 4: choice of the test function and conclusion.}
We now take
$\bt=\bG(\updl{I}_T\udl{u}) -\Pi_{\ife}\Vs(\udl{u})
=\mathbb{C}\big(\bE(\updl{I}_T\udl{u}) -\pk\Ve(\udl{u})\big)$, which is an
admissible test function, and we use that $\mathbb{C}$ is symmetric positive
definite with
\begin{align*}
  \big(\bE ,\mathbb{C}\bE\big)_T = 2\mu\norm{\bE}^2_T +\lambda\norm{tr\bE}^2_T .
\end{align*}
The left-hand side of \reff{Thm-reconstruct-Approx-6} then becomes
\begin{align*}
  2\mu\norm{\bE(\updl{I}_T\udl{u})-\pk\Ve(\udl{u})}^2_T
    &+\lambda\norm{D(\updl{I}_T\udl{u})-\pk(\nabla\cdot\udl{u})}^2_T \\
    &+\frac12\norm{[(\bI-\bK\bS_h)\bK]^\frac12\big(\bG(\updl{I}_T\udl{u})-\Pi_{\ife}\Vs(\udl{u}) \big)\udl{n}}^2_{\ife},
\end{align*}
while the right-hand side is bounded by the Cauchy--Schwarz inequality.  For
the first factor we use that, in the tangential directions,
$\abs{({\bf QPQ}^T)^{-\frac12} \bK^\frac12\bS_h^\frac12}
\leq \big(\frac{\mu}{h_T}\big)^{\frac12}$, which follows from
\reff{Lem-Sh-2}--\reff{Lem-Sh-4}; the second factor is one half of the third
term of the left-hand side.  Absorbing the latter, we obtain
\begin{align}\label{Thm-reconstruct-Approx-1}
  2\mu\norm{\bE\big( \updl{I}_T(\udl{u}) \big)-\pk\Ve(\udl{u})}^2_T
    &+\lambda\norm{D\big( \updl{I}_T(\udl{u}) \big)-\pk(\nabla\cdot\udl{u})}^2_T \nn\\
    &+\frac{1}{4}\norm{[(\bI-\bK\bS_h)\bK]^\frac12\big(\bG\big(\updl{I}_T(\udl{u})\big)
      -\Pi_{\ife}\Vs(\udl{u}) \big)\udl{n}}^2_{\ife} \nn\\
    \leq& \frac{\mu}{h_T}\norm{\udl{u}-\Pi_T\udl{u}}^2_{\ife},
\end{align}
and \Lem{\ref{Lem-Interpolation-properties}} bounds the right-hand side by
$\mu h_T^{2(k+1)}\abs{\udl{u}}^2_{H^{k+2}(T)^d}$, which proves
\reff{Thm-reconstruct-Approx-0}.\Endproof

\begin{rem}[On the projections in \reff{Thm-reconstruct-Approx-0}]\label{Rem-projections}
The interface term of \reff{Thm-reconstruct-Approx-0} involves the face
projection $\Pi_{\ife}$, because this is the projection for which Step 3 of
the proof is exact.  The variant with the cell projection $\pk$ is also
available and is the one used in the sequel; it follows from
\reff{Thm-reconstruct-Approx-0} and from the triangle inequality, at the price
of the additional term $\lambda\abs{\nabla\cdot\udl{u}}_{H^{k+1}(T)}$ coming
from the approximation of $\Vs(\udl{u})$ on $\ife$; see
\reff{Lem-traction-2}--\reff{Lem-traction-3} below.
\end{rem}

The following three bounds on the interface traction are used repeatedly in
the sequel.  The first one says that the reconstructed traction converges at
the optimal rate in the $\bS_h$-weighted trace norm that is natural for the
interface terms, the second one is the corresponding bound for the exact
traction, and the third one is their combination.
\begin{lem}\label{Lem-traction} Under the assumptions of \Thm{\ref{Thm-reconstruct-Approx}}, for any $T\in\cTh$ we have
\begin{align}
  \sqrt{\frac{h_T}{\mu+\lambda}} \norm{\bG(\updl{I}_T\udl{u})-\pk\Vs(\udl{u})}_{\ife}
    &\lesssim \sqrt{\mu} h_T^{k+1}\abs{\udl{u}}_{H^{k+2}(T)^d},\label{Lem-traction-1}\\
  \sqrt{\frac{h_T}{\mu}} \norm{\Vs(\udl{u})-\pk\Vs(\udl{u})}_{\p T}
    &\lesssim \mu^{-\frac12}h_T^{k+1}\big( \mu\abs{\udl{u}}_{H^{k+2}(T)^d}
      +\lambda\abs{\nabla\cdot\udl{u}}_{H^{k+1}(T)}\big),\label{Lem-traction-2}\\
  \sqrt{\frac{h_T}{\mu+\lambda}} \norm{\bG(\updl{I}_T\udl{u})-\Vs(\udl{u})}_{\ife}
    &\lesssim \mu^{-\frac12}h_T^{k+1}\big( \mu\abs{\udl{u}}_{H^{k+2}(T)^d}
      +\lambda\abs{\nabla\cdot\udl{u}}_{H^{k+1}(T)}\big).\label{Lem-traction-3}
\end{align}
\end{lem}
\Proof We prove the three bounds in turn.

For \reff{Lem-traction-1} we observe that
$\bG(\updl{I}_T\udl{u})-\pk\Vs(\udl{u})
=2\mu\big(\bE(\updl{I}_T\udl{u})-\pk\Ve(\udl{u})\big)
+\lambda\big(D(\updl{I}_T\udl{u})-\pk(\nabla\cdot\udl{u})\big)\bI$ is a
polynomial of degree at most $k$ on $T$, so that the discrete trace inequality
\reff{Lem-Discrete-trace-Ineq} applies to each of the two contributions.
Using $\lambda\leq\mu+\lambda$ and $2\mu\leq2(\mu+\lambda)$, we get
\begin{align*}
  \sqrt{\frac{h_T}{\mu+\lambda}} \norm{\bG(\updl{I}_T\udl{u})-\pk\Vs(\udl{u})}_{\ife}
    \lesssim& \norm{\sqrt{\mu}\big( \bE(\updl{I}_T\udl{u})-\pk\Ve(\udl{u}) \big)}_T \\
    &+    \sqrt{\frac{\lambda}{\mu+\lambda}}\norm{\sqrt{\lambda}\big( D(\updl{I}_T\udl{u})-\pk(\nabla\cdot\udl{u}) \big)}_T ,
\end{align*}
and \reff{Thm-reconstruct-Approx-0} bounds both terms by
$\sqrt{\mu} h_T^{k+1}\abs{\udl{u}}_{H^{k+2}(T)^d}$.

The bound \reff{Lem-traction-2} is a pure approximation result.  Since
$\Vs(\udl{u})-\pk\Vs(\udl{u})$ is not a polynomial, we use the multiplicative
trace inequality \reff{Lem-Multip-trace-Ineq} and then
$\Vs(\udl{u})=2\mu\Ve(\udl{u})+\lambda(\nabla\cdot\udl{u})\bI$ together with
\Lem{\ref{Lem-Interpolation-properties}}:
\begin{align*}
  \sqrt{\frac{h_T}{\mu}} \norm{\Vs(\udl{u})-\pk\Vs(\udl{u})}_{\p T}
    \lesssim& \mu^{-\frac12}\Big(\mu\big(\norm{\Ve(\udl{u})-\pk\Ve(\udl{u})}_T
         +h_T\norm{\nabla(\Ve(\udl{u})-\pk\Ve(\udl{u}))}_T\big)\\
    &+    \lambda\big( \norm{\nabla\cdot\udl{u}-\pk(\nabla\cdot\udl{u})}_T
         +h_T\norm{ \nabla\big(\nabla\cdot\udl{u}-\pk(\nabla\cdot\udl{u})\big) }_T \big)\Big) \\
    \lesssim& \mu^{-\frac12}h_T^{k+1}\big( \mu\abs{\udl{u}}_{H^{k+2}(T)^d}
    +\lambda\abs{\nabla\cdot\udl{u}}_{H^{k+1}(T)}\big).
\end{align*}
Finally, \reff{Lem-traction-3} follows from
$\bG(\updl{I}_T\udl{u})-\Vs(\udl{u})
=\big(\bG(\updl{I}_T\udl{u})-\pk\Vs(\udl{u})\big)
-\big(\Vs(\udl{u})-\pk\Vs(\udl{u})\big)$, the triangle inequality,
\reff{Lem-traction-1}--\reff{Lem-traction-2} and $\mu\leq\mu+\lambda$.\Endproof

In what follows, we formulate the regularity of the exact solution using broken Sobolev space $H^m(\cTh)$ for some positive interger $m$ equipped with the norm $\norm{\cdot}_{H^m(\cTh)}^2=\sum_{T\in\cTh}\norm{\cdot}_{H^m(T)}^2$.
\begin{pro}\label{Pro-energy-Err} Let $k\geq 1$. Assume that the exact solution $\udl{u}$ of \reff{weak-form} satisfies $\udl{u}\in H^{k+2}(\cTh)^d, \nabla\cdot\udl{u}\in H^{k+1}(\cTh)$. Let $\updl{u}_h$ be the solution of the discrete problem \reff{discrete-form}, then we have
\be\label{Pro-energy-Err-0}
  \norme{\udl{u}-\updl{u}_h}_e\lesssim \mu^{-\frac12}h^{k+1}\big( \mu\abs{\udl{u}}_{H^{k+2}(\cTh)^d}
    +\lambda\abs{\nabla\cdot\udl{u}}_{H^{k+1}(\cTh)}\big).
\ee
\end{pro}
\Proof As in \Thm{\ref{Thm-reconstruct-Approx}}, we set
$\Pi_{\ife}\udl{u}:= \Pi_{\ife}\av{\udl{u}}$, since $\udl{u}$ is discontinuous
across $\Gamma$, and we write $\updl{e}_h=\updl{u}_h-\updl{I}_h(\udl{u})$.
The proof is organised as follows.  We first reduce the estimate to the
consistency error of the method (Step 1), which we then rewrite in a form in
which only differences between the reconstructed and the exact stress appear
(Step 2).  The five resulting contributions are estimated one by one in
Steps 3--6, and Step 7 concludes.

\smallskip\noindent{\em Step 1: reduction to the consistency error.}
Since $a_h$ is symmetric positive semi-definite and
$\norme{\cdot}^2_h=a_h(\cdot,\cdot)$, the Cauchy--Schwarz inequality gives
\begin{align}\label{Pro-energy-Err-1}
  \norme{\updl{e}_h}_h
    \leq&\sup_{\updl{w}_h\in\updl{V}_{h0},|||\updl{w}_h|||_h=1} a_h(\updl{e}_h,\updl{w}_h)\nn\\
    =& \sup_{\updl{w}_h\in\updl{V}_{h0},|||\updl{w}_h|||_h=1} \cEh(\updl{w}_h),
\end{align}
where, by the discrete problem \reff{discrete-form},
$\cEh(\updl{w}_h)= l(\updl{w}_h) -a_h(\updl{I}_h(\udl{u}),\updl{w}_h)$ is the
consistency error of the method.  Splitting $a_h$ into its local
contributions and $s_T$ into its two parts, we restate it as
\begin{align}\label{Pro-energy-Err-2}
  \cEh(\updl{w}_h) =& \sum_{T\in \cTh}\big[ l_T(\updl{w}_h) -a_T(\updl{I}_T\udl{u},\updl{w}_h)
      -s_T(\updl{I}_T\udl{u},\updl{w}_h)\big] \nn\\
    =&\sum_{T\in \cTh}\big[ l_T(\updl{w}) -a_T(\updl{I}_T\udl{u},\updl{w})
      -s_{T^{\Gamma}}(\updl{I}_T\udl{u},\updl{w}) -s_{T^{\backslash \Gamma}}(\updl{I}_T\udl{u},\updl{w})\big],
\end{align}
in which we drop the subscript $h$ on $\updl{w}$ for brevity.  The last
contribution is estimated separately in Step 6; we now concentrate on the
first three.

\smallskip\noindent{\em Step 2: rewriting the consistency error.}
On every element $T\in\cTh$, we use $\bG=\mathbb{C}\bE$, the symmetry of
$\mathbb{C}$ and of $(\bI-\bK\bS_h)\bK$, and the equation
$-\nabla\cdot\Vs(\udl{u})=\udl{f}$ to write $l_T$ as a volume term, which
gives
\begin{align*}
  a_T(\updl{I}_T\udl{u},\updl{w})
      &+s_{T^{\Gamma}}(\updl{I}_T\udl{u},\updl{w}) -l_T(\updl{w}) \\
    =& 2\mu\big( \bE(\updl{I}_T\udl{u}) ,\bE(\updl{w}) \big)_T
      +\lambda \big( tr\bE(\updl{I}_T\udl{u}) ,tr\bE(\updl{w}) \big)_T \\
    &+\frac12\pd{ (\bI-\bK\bS_h)\bK\bG(\updl{I}_T\udl{u})\udl{n} ,\bG(\updl{w})\udl{n} }_{\ife}
      +\big(\nabla\cdot \Vs(\udl{u}), \udl{w}\big)_T +s_{T^{\Gamma}}(\updl{I}_T\udl{u},\updl{w}) \\
    =& \big( \bG(\updl{I}_T\udl{u}) ,\bE(\updl{w}) \big)_T
      +\frac12\pd{ \bG(\updl{I}_T\udl{u})\udl{n} ,(\bI-\bK\bS_h)\bK\bG(\updl{w})\udl{n} }_{\ife} \\
    &+\big(\nabla\cdot \Vs(\udl{u}), \udl{w}\big)_T +s_{T^{\Gamma}}(\updl{I}_T\udl{u},\updl{w}).
\end{align*}
We now apply the definition \reff{gradient-reconstruct} of the reconstruction
operator to $\updl{w}$, with the admissible test function
$\bt=\bG(\updl{I}_T\udl{u})$, which transforms the first two terms into a
volume term involving $\Ve(\udl{w})$ and boundary terms involving the
differences $\udl{w}_{\uife}-\udl{w}_T$ and $\udl{w}_{\ife}-\udl{w}_T$; we
integrate by parts the term $\big(\nabla\cdot \Vs(\udl{u}), \udl{w}\big)_T$,
and we use \reff{Thm-reconstruct-Approx-jump} in the form
$\av{\udl{u}}-\udl{u}_i=\frac12\bK\Vs(\udl{u})\udl{n}_i$, $i=1,2$, to rewrite
$s_{T^{\Gamma}}(\updl{I}_T\udl{u},\updl{w})$.  Summing over the cells, in
which process the traction terms of two neighbouring cells cancel on the
interior faces not lying on $\Gamma$, we deduce
\begin{align}\label{Pro-energy-Err-3}
  \sum_{T\in\cTh} \big[a_T(\updl{I}_T\udl{u},\updl{w})
      &+s_{T^{\Gamma}}(\updl{I}_T\udl{u},\updl{w}) -l_T(\updl{w})\big] \nn\\
    =& \sum_{T\in\cTh} \Big[ \big( \bG(\updl{I}_T\udl{u}) ,\Ve(\udl{w}) \big)_T
      +\pd{ \bG(\updl{I}_T\udl{u})\udl{n} ,\udl{w}_{\uife}-\udl{w}_T }_{\uife} \nn\\
    &+\pd{ (\bI-\bS_h\bK) \bG(\updl{I}_T\udl{u})\udl{n} ,\udl{w}_{\ife}-\udl{w}_T }_{\ife} \nn\\
    &-\big(\Vs(\udl{u}), \Ve(\udl{w})\big)_T
      -\pd{ \Vs(\udl{u})\udl{n} ,\udl{w}_{\uife}-\udl{w}_T }_{\uife}
      -\pd{ \Vs(\udl{u})\udl{n} ,\udl{w}_{\ife}-\udl{w}_T }_{\ife} \nn\\
    &+2\pd{ \bS_h\Pi_{\ife}\big( \av{\udl{u}} -\Pi_T\udl{u} \big), \udl{w}_{\ife}-\Pi_{\ife}\udl{w}_T }_{\ife}  \Big]\nn\\
    =& \sum_{T\in\cTh} \Big[ \big( \bG(\updl{I}_T\udl{u})-\Vs(\udl{u}) ,\Ve(\udl{w}) \big)_T
      +\pd{ \big( \bG(\updl{I}_T\udl{u})-\Vs(\udl{u}) \big)\udl{n} ,\udl{w}_{\uife}-\udl{w}_T }_{\uife} \nn\\
    &+\pd{ \big( \bG(\updl{I}_T\udl{u})-\Vs(\udl{u}) \big)\udl{n} ,\udl{w}_{\ife}-\udl{w}_T }_{\ife} \nn\\
    &-\pd{ \bS_h\bK\big( \bG(\updl{I}_T\udl{u})-\Vs(\udl{u}) \big)\udl{n},\udl{w}_{\ife}-\Pi_{\ife}\udl{w}_T }_{\ife} \nn\\
    &+2\pd{ \bS_h( \udl{u} -\Pi_T\udl{u} ), \udl{w}_{\ife}-\Pi_{\ife}\udl{w}_T }_{\ife}  \Big].
\end{align}
The gain of this rewriting is that every term of the last expression is a
duality pairing between an interpolation error of the stress and a quantity
that is controlled by $\norme{\updl{w}}_{a_T}$.  We estimate the five terms
successively.

\smallskip\noindent{\em Step 3: the volume term.}
We first split the volume term according to
$\bG(\updl{I}_T\udl{u})-\Vs(\udl{u})
=\big(\bG(\updl{I}_T\udl{u})-\pi_k\Vs(\udl{u})\big)
-\big(\Vs(\udl{u})-\pi_k\Vs(\udl{u})\big)$ and we expand the first part
according to the constitutive law, which gives
\begin{align*}
  \big( \bG(\updl{I}_T\udl{u})-\Vs(\udl{u}) ,\Ve(\udl{w}) \big)_T
    =& 2\mu\big( \bE(\updl{I}_T\udl{u})-\pi_k\Ve(\udl{u}) ,\Ve(\udl{w}) \big)_T
      +\lambda\big( D(\updl{I}_T\udl{u})-\pi_k(\nabla\cdot\udl{u}) ,D(\updl{w}) \big)_T \\
    &+ \lambda\big( D(\updl{I}_T\udl{u})-\pi_k(\nabla\cdot\udl{u}) ,\nabla\cdot\udl{w}-D(\udl{w}) \big)_T
      -\big( \Vs(\udl{u})-\pi_k\Vs(\udl{u}) ,\Ve(\udl{w}) \big)_T .
\end{align*}
The first, second and fourth terms are directly bounded by the
Cauchy--Schwarz inequality, since
$\norm{\sqrt{2\mu}\Ve(\udl{w})}_T\leq\norme{\updl{w}}_E$ and
$\norm{\sqrt{\lambda}D(\udl{w})}_T\leq\norme{\updl{w}}_{a_T}$.  The third one
requires an argument, because $\nabla\cdot\udl{w}-D(\udl{w})$ is not
controlled by the energy norm: we use the definition
\reff{divergence-reconstruct} of the discrete divergence with the test
function $q=\lambda\big(D(\updl{I}_T\udl{u})-\pi_k(\nabla\cdot\udl{u})\big)$,
which turns it into three boundary terms,
\begin{align*}
  \lambda\big( D(\updl{I}_T\udl{u})-\pi_k(\nabla\cdot\udl{u}) &,\nabla\cdot\udl{w}-D(\updl{w}) \big)_T
    = \frac{\lambda}{2}\pd{ \big(D(\updl{I}_T\udl{u})-\pi_k(\nabla\cdot\udl{u})\big)\udl{n},
        (\bI-\bK\bS_h)\bK \bG(\updl{w})\udl{n} }_{\ife} \\
    &- \pd{ \lambda^{\frac12}\big(D(\updl{I}_T\udl{u})-\pi_k(\nabla\cdot\udl{u})\big)\udl{n},
       \lambda^{\frac12}\udl{n}\otimes\udl{n}(\udl{w}_{\uife}-\Pi_{\uife}\udl{w}_T) }_{\uife} \\
    &- \lambda\pd{ \big(D(\updl{I}_T\udl{u})-\pi_k(\nabla\cdot\udl{u})\big)\udl{n},
        (\bI-\bK\bS_h) (\udl{w}_{\ife}-\udl{w}_T) }_{\ife}.
\end{align*}
Each of the three is now bounded by the Cauchy--Schwarz inequality, the
weights being $\big(\frac{h_T}{\mu+\lambda}\big)^{\frac12}$ for the first and
the third one, by \reff{Lem-Sh-4} and \reff{Lem-equivalent-norm-2}, so that
\begin{align*}
  \lambda\big( D(\updl{I}_T\udl{u})-\pi_k(\nabla\cdot\udl{u}) &,\nabla\cdot\udl{w}-D(\updl{w}) \big)_T \\
    \leq& \lambda\sqrt{\frac{h_T}{\mu+\lambda}} \norm{D(\updl{I}_T\udl{u})-\pi_k(\nabla\cdot\udl{u})}_{\ife}
           \norm{[(\bI-\bK\bS_h)\bK]^{\frac12} \bG(\udl{w})\udl{n} }_{\ife} \\
    &+  \norm{\lambda^{\frac12}\big(D(\updl{I}_T\udl{u})-\pi_k(\nabla\cdot\udl{u})\big)}_T \norme{\updl{w}}_{a_T} \\
    &+  \sqrt{\frac{h_T}{\mu+\lambda}}\norm{ \lambda\big(D(\updl{I}_T\udl{u})-\pi_k(\nabla\cdot\udl{u})\big) }_{\ife}
          \norm{\bS_h^{\frac12} (\udl{w}_{\ife}-\Pi_{\ife}\udl{w}_T)}_{\ife} \\
    \leq& \norm{\lambda^{\frac12}\big(D(\updl{I}_T\udl{u})-\pi_k(\nabla\cdot\udl{u})\big)}_T \norme{\updl{w}}_{a_T},
\end{align*}
where the last inequality uses the discrete trace inequality
\reff{Lem-Discrete-trace-Ineq} to convert the interface norms of the
polynomial $D(\updl{I}_T\udl{u})-\pi_k(\nabla\cdot\udl{u})$ into cell norms.
Collecting the four contributions and applying
\reff{Thm-reconstruct-Approx-0} to the first two,
\Lem{\ref{Lem-Interpolation-properties}} to the fourth, we obtain
\begin{align}\label{Pro-energy-Err-4}
  \big( \bG(\updl{I}_T\udl{u})-\Vs(\udl{u}) ,\Ve(\udl{w}) \big)_T
    \leq& \sqrt{2\mu}\norm{\bE(\updl{I}_T\udl{u})-\pi_k\Ve(\udl{u})}_T \norm{\sqrt{2\mu}\Ve(\udl{w})}_T \nn\\
    &+  \sqrt{\lambda}\norm{D(\updl{I}_T\udl{u})-\pi_k(\nabla\cdot\udl{u})}_T
       \Big( \norm{\sqrt{\lambda}D(\udl{w})}_T +\norme{\updl{w}}_{a_T}\Big) \nn\\
    &+ \mu^{-\frac12}\big( 2\mu\norm{\Ve(\udl{u})-\pi_k\Ve(\udl{u})}_T
           +\lambda\norm{\nabla\cdot\udl{u}-\pi_k(\nabla\cdot\udl{u})}_T\big)\norm{\sqrt{2\mu}\Ve(\udl{w})}_T \nn\\
    \lesssim& \mu^{-\frac12}h^{k+1}\big( \mu\abs{\udl{u}}_{H^{k+2}(T)^d}
               +\lambda\abs{\nabla\cdot\udl{u}}_{H^{k+1}(T)}\big)\norme{\updl{w}}_{a_T}.
\end{align}

\smallskip\noindent{\em Step 4: the faces away from the interface.}
On $\uife$ we split the pairing in the same way; in the first part we may
insert the projection $\Pi_{\uife}$ in front of $\udl{w}_T$, because
$\big(\bG(\updl{I}_T\udl{u})-\pi_k\Vs(\udl{u})\big)\udl{n}$ has its components
in $\mP^k(F)^d$ on every face.  \Thm{\ref{Thm-reconstruct-Approx}}, the
discrete trace inequality \reff{Lem-Discrete-trace-Ineq} for the first part
and the continuous trace inequality \reff{Lem-Multip-trace-Ineq} for the
second one then give
\begin{align}\label{Pro-energy-Err-5}
  \pd{ \big( \bG(\updl{I}_T\udl{u})-\Vs(\udl{u}) \big)\udl{n} ,\udl{w}_{\uife}-\udl{w}_T }_{\uife}
    =& \pd{ \big( \bG(\updl{I}_T\udl{u})-\pi_k\Vs(\udl{u}) \big)\udl{n} ,\udl{w}_{\uife}-\Pi_{\uife}\udl{w}_T }_{\uife} \nn\\
    &\hspace{-2.5cm} - \pd{ \big( \Vs(\udl{u})-\pi_k\Vs(\udl{u}) \big)\udl{n} ,\udl{w}_{\uife}-\udl{w}_T }_{\uife} \nn\\
    &\hspace{-3cm} \leq \big(\sqrt{2\mu}\norm{\bE(\updl{I}_T\udl{u})-\pi_k\Ve(\udl{u})}_T
           +\sqrt{\lambda}\norm{D(\updl{I}_T\udl{u})-\pi_k(\nabla\cdot\udl{u})}_T\big)\norme{\updl{w}}_{a_T} \nn\\
    &\hspace{-2.5cm} +\mu^{-\frac12}h^{\frac12}\big( 2\mu\norm{\Ve(\udl{u})-\pi_k\Ve(\udl{u})}_{\uife}
           +\lambda\norm{\nabla\cdot\udl{u}-\pi_k(\nabla\cdot\udl{u})}_{\uife}\big) \norme{\updl{w}}_E \nn\\
    &\hspace{-3cm} \lesssim \mu^{-\frac12}h^{k+1}\big( \mu\abs{\udl{u}}_{H^{k+2}(T)^d}
           +\lambda\abs{\nabla\cdot\udl{u}}_{H^{k+1}(T)}\big)\norme{\updl{w}}_{a_T}.
\end{align}

\smallskip\noindent{\em Step 5: the interface terms.}
The third and fourth terms in \reff{Pro-energy-Err-3} have to be treated
together, because neither of them is separately controlled by the energy norm:
the third one involves the trace $\udl{w}_T$ instead of its projection, and
the fourth one carries the factor $\bS_h\bK$.  Adding and subtracting
$\Pi_{\ife}\udl{w}_T$, and using
$\bI-\bS_h\bK={\bf QPQ}^T\bS_h$ from \reff{Lem-Sh-1}, we obtain
\begin{align}\label{Pro-energy-Err-6}
  \pd{ \big( \bG(\updl{I}_T\udl{u})-\Vs(\udl{u}) \big)\udl{n} ,\udl{w}_{\ife}-\udl{w}_T }_{\ife}
    &-\pd{ \bS_h\bK\big( \bG(\updl{I}_T\udl{u})-\Vs(\udl{u}) \big)\udl{n},\udl{w}_{\ife}-\Pi_{\ife}\udl{w}_T }_{\ife} \nn\\
    &\hspace{-5.5cm}= \pd{ \big( \bG(\updl{I}_T\udl{u})-\Vs(\udl{u}) \big)\udl{n} ,\Pi_{\ife}\udl{w}_T-\udl{w}_T }_{\ife}\nn\\
    &\hspace{-4cm}+\pd{ (\bI-\bS_h\bK)\big( \bG(\updl{I}_T\udl{u})-\Vs(\udl{u}) \big)\udl{n},
     \udl{w}_{\ife}-\Pi_{\ife}\udl{w}_T }_{\ife} .
\end{align}
In the first term of \reff{Pro-energy-Err-6} we split again
$\bG(\updl{I}_T\udl{u})-\Vs(\udl{u})$ around $\pi_k\Vs(\udl{u})$ and expand
$\bG-\pi_k\Vs$ according to the constitutive law, which produces the two terms
$\pd{ 2\mu\big( \bE(\updl{I}_T\udl{u})-\pi_k\Ve(\udl{u}) \big)\udl{n}
,\Pi_{\ife}\udl{w}_T-\udl{w}_T }_{\ife}$ and
$\pd{ \lambda\udl{n}^T\big( D(\updl{I}_T\udl{u})-\pi_k(\nabla\cdot\udl{u})
\big)\udl{n} ,\big(\Pi_{\ife}\udl{w}_T-\udl{w}_T\big)\cdot\udl{n} }_{\ife}$.
The second of these vanishes: by \Rem{\ref{Rem-interpolants}}, the normal
component of $\Pi_{\ife}\udl{w}_T-\udl{w}_T$ is orthogonal to
$\mP^k(\ife)$ when the BDM interpolation is used for the divergence part, so
that only the $\mu$-weighted contribution survives.  Hence
\begin{align}\label{Pro-energy-Err-6b}
  \pd{ \big( \bG(\updl{I}_T\udl{u})-\Vs(\udl{u}) \big)\udl{n} ,\Pi_{\ife}\udl{w}_T-\udl{w}_T }_{\ife}
    =& \pd{ 2\mu\big( \bE(\updl{I}_T\udl{u})-\pi_k\Ve(\udl{u}) \big)\udl{n} ,\Pi_{\ife}\udl{w}_T-\udl{w}_T }_{\ife} \nn\\
    &-\pd{ \big( \Vs(\udl{u})-\pi_k\Vs(\udl{u}) \big)\udl{n} ,\Pi_{\ife}\udl{w}_T-\udl{w}_T }_{\ife}.
\end{align}
The two remaining factors are now estimated by \Lem{\ref{Lem-Fce-Korn}}, which
gives
\begin{align*}
  \sqrt{\frac{\mu}{h_T}} \norm{\udl{w}_T -\Pi_{\ife}\udl{w}_T}_{\ife}
    \lesssim \sqrt{\mu}\norm{\Ve(\udl{w}_T) }_T
    \lesssim \norme{\updl{w}}_E \lesssim \norme{\updl{w}}_{a_T},
\end{align*}
and by \reff{Lem-traction-1}--\reff{Lem-traction-2} of
\Lem{\ref{Lem-traction}}, so that
\begin{align}\label{Pro-energy-Err-8}
  &\pd{ \big( \bG(\updl{I}_T\udl{u})-\Vs(\udl{u}) \big)\udl{n} ,\Pi_{\ife}\udl{w}_T-\udl{w}_T }_{\ife} \nn\\
    &\hspace{1cm}\leq \sqrt{\frac{h_T}{\mu}} \Big( 2\mu\norm{\bE(\updl{I}_T\udl{u})-\pi_k\Ve(\udl{u})}_{\ife}
                                    +\norm{\Vs(\udl{u})-\pi_k\Vs(\udl{u})}_{\ife} \Big)
         \sqrt{\frac{\mu}{h_T}} \norm{\udl{w}_T -\Pi_{\ife}\udl{w}_T}_{\ife} \nn\\
    &\hspace{1cm}\lesssim \mu^{-\frac12}h^{k+1}\big( \mu\abs{\udl{u}}_{H^{k+2}(T)^d}
             +\lambda\abs{\nabla\cdot\udl{u}}_{H^{k+1}(T)}\big) \norme{\updl{w}}_{a_T}.
\end{align}
For the second term of \reff{Pro-energy-Err-6} we use the factorisation
$\bI-\bS_h\bK={\bf QPQ}^T\bS_h={\bf QPQ}^T\bS_h^{\frac12}\bS_h^{\frac12}$
together with $\abs{{\bf QPQ}^T\bS_h^{\frac12}}\leq
\big(\frac{h_T}{\mu+\lambda}\big)^{\frac12}$ in the normal direction and
$\big(\frac{h_T}{\mu}\big)^{\frac12}$ in the tangential ones, see
\reff{Lem-Sh-2}--\reff{Lem-Sh-4}, which moves one half of the weight onto the
test function and produces exactly the interface contribution to
$\norme{\updl{w}}_{a_T}$; with \reff{Lem-traction-3} we get
\begin{align}\label{Pro-energy-Err-9}
  &\pd{ (\bI-\bS_h\bK)\big( \bG(\updl{I}_T\udl{u})-\Vs(\udl{u}) \big)\udl{n},
     \udl{w}_{\ife}-\Pi_{\ife}\udl{w}_T }_{\ife} \nn\\
    &\hspace{1cm}= \pd{ {\bf QPQ}^T\bS_h^{\frac12}\big( \bG(\updl{I}_T\udl{u})-\Vs(\udl{u}) \big)\udl{n},
     \bS_h^{\frac12}\big(\udl{w}_{\ife}-\Pi_{\ife}\udl{w}_T\big) }_{\ife} \nn\\
    &\hspace{1cm}\leq \sqrt{\frac{h_T}{\mu+\lambda}} \norm{\bG(\updl{I}_T\udl{u})-\Vs(\udl{u})}_{\ife} \norme{\updl{w}}_{a_T} \nn\\
    &\hspace{1cm}\lesssim \mu^{-\frac12}h^{k+1}\big( \mu\abs{\udl{u}}_{H^{k+2}(T)^d}
             +\lambda\abs{\nabla\cdot\udl{u}}_{H^{k+1}(T)}\big) \norme{\updl{w}}_{a_T}.
\end{align}
Combining \reff{Pro-energy-Err-6}, \reff{Pro-energy-Err-8} and
\reff{Pro-energy-Err-9}, it is derived that
\begin{align}\label{Pro-energy-Err-10}
  &\pd{ \big( \bG(\updl{I}_T\udl{u})-\Vs(\udl{u}) \big)\udl{n} ,\udl{w}_{\ife}-\udl{w}_T }_{\ife}
    -\pd{ \bS_h\bK\big( \bG(\updl{I}_T\udl{u})-\Vs(\udl{u}) \big)\udl{n},\udl{w}_{\ife}-\Pi_{\ife}\udl{w}_T }_{\ife}\nn\\
    &\hspace{1cm}\lesssim \mu^{-\frac12}h^{k+1}\big( \mu\abs{\udl{u}}_{H^{k+2}(T)^d}
             +\lambda\abs{\nabla\cdot\udl{u}}_{H^{k+1}(T)}\big) \norme{\updl{w}}_{a_T}.
\end{align}

\smallskip\noindent{\em Step 6: the two stabilization terms.}
The last term in \reff{Pro-energy-Err-3} only involves the tangential part of
$\bS_h$, since the normal component of $\udl{u}-\Pi_T\udl{u}$ vanishes on the
faces by \Rem{\ref{Rem-interpolants}}; recalling that
$\abs{\bS_{h,t}}\leq\mu/h_T$ by \reff{Lem-Sh-2} and using the continuous trace
inequality \reff{Lem-Multip-trace-Ineq}, we get
\begin{align}\label{Pro-energy-Err-11}
  \pd{ \bS_h( \udl{u} -\Pi_T\udl{u} ), \udl{w}_{\ife}-\Pi_{\ife}\udl{w}_T }_{\ife}
    =& \pd{ \bS_{h,t}( \udl{u} -\Pi_T\udl{u} ), \udl{w}_{\ife}-\Pi_{\ife}\udl{w}_T }_{\ife} \nn\\
    \leq& \norm{\bS_{h,t}^{\frac12}( \udl{u} -\Pi_T\udl{u} )}_{\ife} \norme{\updl{w}}_{a_T} \nn\\
    \leq& \mu^{\frac12}h^{-\frac12} \big( h^{-\frac12}\norm{\udl{u} -\Pi_T\udl{u}}_T
          +h^{\frac12}\norm{\nabla( \udl{u} -\Pi_T\udl{u} )}_T \big) \norme{\updl{w}}_{a_T} \nn\\
    \lesssim& \mu^{-\frac12}h^{k+1}\mu\abs{\udl{u}}_{H^{k+2}(T)^d} \norme{\updl{w}}_{a_T}.
\end{align}
The contribution $s_{T^{\backslash \Gamma}}(\updl{I}_T\udl{u},\updl{w})$ left
aside in Step 1 is bounded in the same spirit; observing (c.f.\ Lemma 7 in
\cite{HL02}) that $s_{T^{\backslash \Gamma}}$ is a positive semi-definite
bilinear form, the Cauchy--Schwarz inequality and
\Lem{\ref{Lem-Interpolation-properties}} give
\begin{align}\label{Pro-energy-Err-13}
  s_{T^{\backslash \Gamma}}(\updl{I}_T\udl{u},\updl{w})
    \leq& s_{T^{\backslash \Gamma}}^{\frac12}(\updl{I}_T\udl{u},\updl{I}_T\udl{u}) \norme{\updl{w}}_{a_T} \nn\\
    \leq& \mu^{\frac12}h^{-\frac12} \norm{\Pi_{\uife}(\udl{u}-\Pi_T\udl{u})}_{\uife} \norme{\updl{w}}_{a_T} \nn\\
    \leq& \mu^{\frac12}h^{-\frac12} \norm{\udl{u}-\Pi_T\udl{u}}_{\uife} \norme{\updl{w}}_{a_T} \nn\\
    \lesssim& \mu^{-\frac12}h^{k+1}\mu\abs{\udl{u}}_{H^{k+2}(T)^d} \norme{\updl{w}}_{a_T},
\end{align}
where the $\lambda$-weighted part of $s_{T^{\backslash \Gamma}}$ has been
discarded, again because the normal component of the interpolation error
vanishes on the faces.

\smallskip\noindent{\em Step 7: conclusion.}
Applying \reff{Pro-energy-Err-4}, \reff{Pro-energy-Err-5},
\reff{Pro-energy-Err-10} and \reff{Pro-energy-Err-11} in
\reff{Pro-energy-Err-3}, and summing over the cells with the Cauchy--Schwarz
inequality in $\ell^2$, we obtain
\begin{align}\label{Pro-energy-Err-12}
  \sum_{T\in\cTh} \big[a_T(\updl{I}_T\udl{u},\updl{w})
      &+s_{T^{\Gamma}}(\updl{I}_T\udl{u},\updl{w}) -l_T(\updl{w})\big] \nn\\
    \lesssim& \mu^{-\frac12}h^{k+1}\big( \mu\abs{\udl{u}}_{H^{k+2}(\cTh)^d}
             +\lambda\abs{\nabla\cdot\udl{u}}_{H^{k+1}(\cTh)}\big) \norme{\updl{w}}_h.
\end{align}
By \reff{Pro-energy-Err-1}--\reff{Pro-energy-Err-2},
\reff{Pro-energy-Err-12}--\reff{Pro-energy-Err-13} and
\Lem{\ref{Lem-equivalent-norm}}, we arrive at
\begin{align}\label{Pro-energy-Err-14}
  \norme{\updl{u}_h-\updl{I}_h(\udl{u})}_e
    \lesssim \norme{\updl{u}_h-\updl{I}_h(\udl{u})}_h
    \lesssim \mu^{-\frac12}h^{k+1}\big( \mu\abs{\udl{u}}_{H^{k+2}(\cTh)^d}
             +\lambda\abs{\nabla\cdot\udl{u}}_{H^{k+1}(\cTh)}\big) .
\end{align}
Then \reff{Pro-energy-Err-0} follows from \reff{Pro-energy-Err-14}, the
approximation properties of the interpolation operator $\updl{I}_h$ collected
in \Lem{\ref{Lem-Interpolation-properties}} and the triangle
inequality. \Endproof

In order to give the $L^2$-error estimate, we start by proving three
elementary Lemmas.  The first one measures the price to pay when the exact
solution is replaced by its interpolate in the first argument of the local
bilinear form; the only remainder is an interface term.
\begin{lem}\label{Lem-reconstruct-aT} Let $\udl{u}$ be the solution of the original problem \reff{elasticity-interface}. For any $\updl{w}\in \updl{V}_T$, we have
\be\label{Lem-reconstruct-aT-0}
  a_T(\updl{I}_T (\udl{u}), \updl{w}) =a_T( \udl{u}, \updl{w})
    -\pd{\bS_h(\udl{u}-\Pi_T\udl{u}), \bK\big( 2\mu \bt(\updl{w}) +\lambda \bI tr\bt(\updl{w}) \big)\udl{n} }_{\ife},
\ee
where the operator $\bt$ can be the symmetric gradient $\Ve$ and the symmetric gradient reconstruction $\bE$ according to vector $\udl{w}$ and $\updl{w}$, respectively.
\end{lem}

\Proof Throughout the proof we abbreviate
$\bt_{\mathbb{C}}(\updl{w}):=2\mu \bt(\updl{w}) +\lambda \bI tr\bt(\updl{w})$,
so that the assertion reads
$a_T(\updl{I}_T\udl{u},\updl{w})=a_T(\udl{u},\updl{w})
-\pd{\bS_h(\udl{u}-\Pi_T\udl{u}),\bK\bt_{\mathbb{C}}(\updl{w})\udl{n}}_{\ife}$.

\smallskip\noindent{\em Step 1: the reconstruction is tested with
$\bt_{\mathbb{C}}(\updl{w})$.}
By the definition \reff{Local-aT} of $a_T$ and since
$D=tr\bE$ implies
\begin{align*}
  \lambda\big(D(\updl{I}_T\udl{u}),tr\bt(\updl{w})\big)_T
    =\big(\bE(\updl{I}_T\udl{u}),\lambda\,\bI\, tr\bt(\updl{w})\big)_T,
\end{align*}
the first two terms of $a_T(\updl{I}_T\udl{u},\updl{w})$ can be gathered into
a single inner product,
\begin{align}\label{Lem-reconstruct-aT-1}
  a_T( \updl{I}_T (\udl{u}), \updl{w} )
    =& \big(\bE(\updl{I}_T\udl{u}), \bt_{\mathbb{C}}(\updl{w})\big)_T
      +\frac12\pd{ (\bI-\bK\bS_h)\bK \bG(\updl{I}_T\udl{u})\udl{n}
      ,\bt_{\mathbb{C}}(\updl{w})\udl{n} }_{\ife}.
\end{align}
The right-hand side of \reff{Lem-reconstruct-aT-1} is exactly the left-hand
side of the definition \reff{gradient-reconstruct} of the reconstruction
operator, tested with $\bt=\bt_{\mathbb{C}}(\updl{w})\in\mP^k(T)^{d\times d}$.
Hence
\begin{align}\label{Lem-reconstruct-aT-2}
  a_T( \updl{I}_T (\udl{u}), \updl{w} )
    =& \big(\Ve(\Pi_T\udl{u}), \bt_{\mathbb{C}}(\updl{w})\big)_T
      +\pd{ \Pi_{\uife}\udl{u} -\Pi_T\udl{u} ,\bt_{\mathbb{C}}(\updl{w})\udl{n} }_{\uife} \nn\\
    &+\pd{ (\bI-\bK\bS_h) \big( \Pi_{\ife}\av{\udl{u}} -\Pi_T\udl{u}\big)
      ,\bt_{\mathbb{C}}(\updl{w})\udl{n} }_{\ife}.
\end{align}

\smallskip\noindent{\em Step 2: removing the projections.}
We integrate by parts the volume term of \reff{Lem-reconstruct-aT-2} and use
that $\bt_{\mathbb{C}}(\updl{w})\udl{n}|_F\in\mP^k(F)^d$ for every
$F\in\cFT$, so that the face projections may be dropped, and that
$\nabla\cdot\bt_{\mathbb{C}}(\updl{w})\in\mP^{k+1}(T)^d$, so that $\Pi_T$ may
be dropped in the volume term.  This gives
\begin{align}\label{Lem-reconstruct-aT-3}
  a_T( \updl{I}_T (\udl{u}), \updl{w} )
    =& -\big(\udl{u}, \nabla\cdot\bt_{\mathbb{C}}(\updl{w}) \big)_T
      +\pd{ \udl{u} ,\bt_{\mathbb{C}}(\updl{w})\udl{n} }_{\uife}
      +\pd{ \av{\udl{u}} ,\bt_{\mathbb{C}}(\updl{w})\udl{n} }_{\ife} \nn\\
    &-\pd{ \bS_h \big( \av{\udl{u}} -\Pi_T\udl{u}\big)
      ,\bK\bt_{\mathbb{C}}(\updl{w})\udl{n} }_{\ife},
\end{align}
where we have also used the symmetry of $\bK$ and $\bS_h$ in the last term.

\smallskip\noindent{\em Step 3: back to $a_T(\udl{u},\updl{w})$.}
Integrating by parts backwards, the first three terms of
\reff{Lem-reconstruct-aT-3} produce
$\big(\Ve(\udl{u}), \bt_{\mathbb{C}}(\updl{w}) \big)_T
-\pd{\udl{u}-\av{\udl{u}},\bt_{\mathbb{C}}(\updl{w})\udl{n}}_{\ife}$.  We now
use the interface identities $\udl{u}-\av{\udl{u}}=\frac12 \jm{\udl{u}}$ and
$-\jm{\udl{u}}=\bK \Vs(\udl{u})\udl{n}$ of
\reff{Thm-reconstruct-Approx-jump}, and we insert
$\av{\udl{u}}-\Pi_T\udl{u}=(\udl{u}-\Pi_T\udl{u})-\frac12\jm{\udl{u}}$ in the
last term of \reff{Lem-reconstruct-aT-3}, so that the two contributions
involving $\jm{\udl{u}}$ combine into
$\frac12\pd{(\bI-\bK\bS_h)\bK\Vs(\udl{u})\udl{n},
\bt_{\mathbb{C}}(\updl{w})\udl{n}}_{\ife}$.  We obtain
\begin{align*}
  a_T( \updl{I}_T (\udl{u}), \updl{w} )
    =& 2\mu\big(\Ve(\udl{u}), \bt(\updl{w}) \big)_T +\lambda\big(\nabla\cdot\udl{u}, tr\bt(\updl{w}) \big)_T\\
    &+\frac12\pd{(\bI-\bK\bS_h)\bK\Vs(\udl{u})\udl{n},\bt_{\mathbb{C}}(\updl{w})\udl{n} }_{\ife}
    -\pd{ \bS_h \big( \udl{u} -\Pi_T\udl{u}\big) ,\bK\bt_{\mathbb{C}}(\updl{w})\udl{n} }_{\ife},
\end{align*}
and the first three terms of the right-hand side are precisely
$a_T( \udl{u}, \updl{w} )$, which completes the proof.\Endproof

We assume the following auxiliary problem has a unique solution
\be\label{Auxiliary-problem}
  \left \{
    \ba {ll}
      -\nabla \cdot \Vs(\zt) = \udl{g} \quad & {\rm in} \quad \Omega_1\cup\Omega_2, \\
      \jm{\zt}+\bK\Vs(\udl{\zeta})\udl{n}=\udl{0} \quad & {\rm on} \quad \Gamma, \\
      \big[\Vs(\zt)\udl{n}\big]=\udl{0} \quad & {\rm on} \quad \Gamma, \\
      \zt=\udl{0} \quad & {\rm on} \quad \p\Omega,
    \ea
  \right.
\ee
and satisfies the regularity
\begin{equation}\label{Regularity-Asm}
  \mu\norm{\zt}_{H^2(\oo)} +\lambda\norm{\nabla\cdot\zt}_{H^1(\oo)} \lesssim \norm{\udl{g}}_\Omega.
\end{equation}
Note that \reff{Auxiliary-problem} is the adjoint of
\reff{elasticity-interface}, which is self-adjoint, so that the same interface
conditions appear.

\begin{lem}\label{Lem-reconstruct-Difference} Let $\udl{u}$ and $\zt$ be the solutions of \reff{elasticity-interface} and \reff{Auxiliary-problem}, respectively. We have the following difference of the reconstruction operator in the local inner product
\begin{align}\label{Lem-reconstruct-Difference-0}
  a_T( \udl{u}, \zt) -a_T(\updl{I}_T (\udl{u}), \updl{I}_T (\zt))
    =&a_T(\udl{u}-\updl{I}_T (\udl{u}), \zt-\updl{I}_T (\zt))
      +\pd{\bS_h(\udl{u}-\Pi_T\udl{u}), \bK\bG( \updl{I}_T \zt )\udl{n} }_{\ife}\nn\\
    &+\pd{\bS_h \bK\bG( \updl{I}_T \udl{u} )\udl{n}, \zt-\Pi_T\zt }_{\ife}.
\end{align}
\end{lem}

\Proof Expanding the quadratic term, which is licit because $a_T$ is bilinear,
we have
\begin{align*}
  a_T( \udl{u}-\updl{I}_T \udl{u}, \zt-\updl{I}_T \zt) =a_T( \udl{u}, \zt)
    -a_T( \udl{u}, \updl{I}_T \zt)
    -a_T(\updl{I}_T \udl{u}, \zt)
    +a_T(\updl{I}_T \udl{u}, \updl{I}_T \zt).
\end{align*}
The two mixed terms are computed with \Lem{\ref{Lem-reconstruct-aT}}, applied
to $\udl{u}$ with $\updl{w}=\updl{I}_T\zt$ and, by symmetry of $a_T$, to $\zt$
with $\updl{w}=\updl{I}_T\udl{u}$, which gives
\begin{align*}
  a_T( \udl{u}, \updl{I}_T \zt)
    &=a_T(\updl{I}_T \udl{u}, \updl{I}_T \zt)
    +\pd{\bS_h(\udl{u}-\Pi_T\udl{u}), \bK\bG( \updl{I}_T \zt )\udl{n} }_{\ife},\\
  a_T(\updl{I}_T \udl{u}, \zt)
    &=a_T(\updl{I}_T \udl{u}, \updl{I}_T \zt)
    +\pd{\bS_h \bK\bG( \updl{I}_T \udl{u} )\udl{n}, \zt-\Pi_T\zt }_{\ife}.
\end{align*}
Combining these three equalities gives
\reff{Lem-reconstruct-Difference-0}.\Endproof

The last lemma is the consistency of the local bilinear form itself.  Here and
in the sequel, we use for a smooth function $\udl{v}$ the notation
\begin{align}\label{aT-extended}
  a_T( \udl{v}-\updl{I}_T\udl{v}, \udl{z}-\updl{I}_T\udl{z})
    :=& 2\mu\big( \Ve(\udl{v})-\bE(\updl{I}_T\udl{v})
       ,\Ve(\udl{z})-\bE(\updl{I}_T\udl{z}) \big)_T \nn\\
    &+\lambda\big( \nabla\cdot\udl{v}-D(\updl{I}_T\udl{v})
       ,\nabla\cdot\udl{z}-D(\updl{I}_T\udl{z}) \big)_T \nn\\
    &+ \frac12\pd{(\bI-\bK\bS_h)\bK\big(\Vs(\udl{v})-\bG(\updl{I}_T\udl{v})\big)\udl{n}
         ,\big(\Vs(\udl{z})-\bG(\updl{I}_T\udl{z})\big)\udl{n} }_{\ife},
\end{align}
which is consistent with \reff{Local-aT} and with \Lem{\ref{Lem-reconstruct-aT}}.
Since $\mathbb{C}$ and $(\bI-\bK\bS_h)\bK$ are symmetric positive
semi-definite, so is the bilinear form \reff{aT-extended}, and the
Cauchy--Schwarz inequality is therefore available for it.
\begin{lem}\label{Lem-aT-consistency} Let $m\in\{1,k+1\}$ and let $\udl{v}\in H^{m+1}(T)^d$ with $\nabla\cdot\udl{v}\in H^{m}(T)$. Then
\begin{align}\label{Lem-aT-consistency-0}
  a_T( \udl{v}-\updl{I}_T\udl{v}, \udl{v}-\updl{I}_T\udl{v})^{\frac12}
    \lesssim \mu^{-\frac12}h_T^{m}\big( \mu\abs{\udl{v}}_{H^{m+1}(T)^d}
      +\lambda\abs{\nabla\cdot\udl{v}}_{H^{m}(T)}\big).
\end{align}
\end{lem}
\Proof We bound the three contributions to \reff{aT-extended} separately.  We
shall use that the proof of \Thm{\ref{Thm-reconstruct-Approx}} only relies on
\Lem{\ref{Lem-Interpolation-properties}}, so that
\reff{Thm-reconstruct-Approx-0} holds with $h_T^{k+1}$ and
$\abs{\cdot}_{H^{k+2}}$ replaced by $h_T^{m}$ and $\abs{\cdot}_{H^{m+1}}$ for
$m\in\{1,k+1\}$.

For the first contribution, the triangle inequality,
\Lem{\ref{Lem-Interpolation-properties}} and
\reff{Thm-reconstruct-Approx-0} give
\begin{align*}
  \sqrt{2\mu}\norm{\Ve(\udl{v})-\bE(\updl{I}_T\udl{v})}_T
    \leq& \sqrt{2\mu}\norm{\Ve(\udl{v})-\pk\Ve(\udl{v})}_T
      +\sqrt{2\mu}\norm{\pk\Ve(\udl{v})-\bE(\updl{I}_T\udl{v})}_T \\
    \lesssim& \sqrt{\mu}h_T^{m}\abs{\udl{v}}_{H^{m+1}(T)^d}.
\end{align*}
For the second one, we argue in the same way and then use
$\sqrt{\lambda\mu}\leq\frac12(\mu+\lambda)$ and
$\abs{\nabla\cdot\udl{v}}_{H^{m}(T)}\lesssim\abs{\udl{v}}_{H^{m+1}(T)^d}$,
which gives
\begin{align*}
  \sqrt{\lambda}\norm{\nabla\cdot\udl{v}-D(\updl{I}_T\udl{v})}_T
    \lesssim& \sqrt{\lambda}h_T^{m}\abs{\nabla\cdot\udl{v}}_{H^{m}(T)}
      +\sqrt{\mu}h_T^{m}\abs{\udl{v}}_{H^{m+1}(T)^d} \\
    \lesssim& \mu^{-\frac12}h_T^{m}\big( \mu\abs{\udl{v}}_{H^{m+1}(T)^d}
      +\lambda\abs{\nabla\cdot\udl{v}}_{H^{m}(T)}\big).
\end{align*}
For the interface contribution we split the traction into its normal and
tangential parts.  Since $\bI$ has no tangential contribution to the traction,
the tangential part of
$\big(\Vs(\udl{v})-\bG(\updl{I}_T\udl{v})\big)\udl{n}$ is
$2\mu\big(\Ve(\udl{v})-\bE(\updl{I}_T\udl{v})\big)\udl{n}\cdot\udl{t}_i$,
which involves $\mu$ only, and by \reff{Lem-Sh-4} the corresponding weight is
$\udl{t}_i^T(\bI-\bK\bS_h)\bK\udl{t}_i\leq h_T/\mu$, whereas the weight of the
normal part is $\udl{n}^T(\bI-\bK\bS_h)\bK\udl{n}\leq h_T/(\mu+\lambda)$.
Using the multiplicative trace inequality \reff{Lem-Multip-trace-Ineq} for the
non-polynomial parts and the discrete trace inequality
\reff{Lem-Discrete-trace-Ineq} for the polynomial ones, as in the proof of
\Lem{\ref{Lem-traction}}, we obtain
\begin{align*}
  \norm{[(\bI-\bK\bS_h)\bK]^\frac12\big(\Vs(\udl{v})-\bG(\updl{I}_T\udl{v})\big)\udl{n}}_{\ife}
    \lesssim& \sqrt{\frac{h_T}{\mu}}\,
      2\mu\norm{\Ve(\udl{v})-\bE(\updl{I}_T\udl{v})}_{\ife} \\
    &+\sqrt{\frac{h_T}{\mu+\lambda}}
      \norm{\big(\Vs(\udl{v})-\bG(\updl{I}_T\udl{v})\big)\udl{n}\cdot\udl{n}}_{\ife} \\
    \lesssim& \mu^{-\frac12}h_T^{m}\big( \mu\abs{\udl{v}}_{H^{m+1}(T)^d}
      +\lambda\abs{\nabla\cdot\udl{v}}_{H^{m}(T)}\big),
\end{align*}
where the last inequality also uses $\mu\leq\mu+\lambda$.  Collecting the
three bounds gives \reff{Lem-aT-consistency-0}.\Endproof

\begin{pro}\label{Pro-L2-Err} Let $k\geq 1$. Assume that the exact solution $\udl{u}$ of \reff{weak-form} satisfies $\udl{u}\in H^{k+2}(\cTh)^d, \nabla\cdot\udl{u}\in H^{k+1}(\cTh)$. Let $\updl{u}_h$ be the solution of the discrete problem \reff{discrete-form}, then we have
\be\label{Pro-L2-Err-0}
  \norm{\udl{u}_h-\udl{I}_h(\udl{u})}_\Omega\lesssim \mu^{-1}h^{k+2}\big( \mu\abs{\udl{u}}_{H^{k+2}(\cTh)^d}
    +\lambda\abs{\nabla\cdot\udl{u}}_{H^{k+1}(\cTh)}\big).
\ee
\end{pro}
\Proof We use a duality argument.  Set
$\udl{g}=\udl{e}_h=\udl{u}_h-\udl{I}_h(\udl{u})$ in the auxiliary problem
\reff{Auxiliary-problem} and let $\zt$ be its solution, so that
$\norm{\udl{e}_h}^2_\Omega=(\udl{e}_h,-\nabla\cdot\Vs(\zt))_\Omega$.  The
proof splits this quantity into two parts (Step 1), which are then estimated
in Steps 2--3 and Steps 4--5, respectively; the extra power of $h$ with
respect to \Pro{\ref{Pro-energy-Err}} comes from the fact that only the
$H^2$-regularity \reff{Regularity-Asm} of the dual solution is used.

\smallskip\noindent{\em Step 1: splitting of the $L^2$-norm.}
Integrating by parts cell by cell and using the continuity of
$\Vs(\zt)\udl{n}$ across the interior faces, together with the fact that the
face components of $\updl{e}_h$ are single valued and vanish on $\p\Omega$, we
obtain
\begin{align}\label{Pro-L2-Err-1}
  \norm{\udl{e}_h}^2_\Omega =& \sum_{T\in\cTh}(\udl{e}_T, -\nabla\cdot\Vs(\zt))_T \nn\\
    =& \sum_{T\in\cTh}\Big[\big(\Ve(\udl{e}_T), \Vs(\zt)\big)_T +\pd{ \udl{e}_{\uife}-\udl{e}_T, \Vs(\zt)\udl{n}}_{\uife}
       +\pd{ \udl{e}_{\ife}-\udl{e}_T, \Vs(\zt)\udl{n}}_{\ife} \Big].
\end{align}
We now insert $\bG(\updl{I}_T\zt)$, i.e.\ the reconstructed dual stress, in the
volume term, and we treat the resulting term
$\big(\Ve(\udl{e}_T), \bG(\updl{I}_T\zt) \big)_T$ with the definition
\reff{gradient-reconstruct} of the reconstruction operator applied to
$\updl{e}_T$ with the test function $\bt=\bG(\updl{I}_T\zt)$, which yields
\begin{align*}
  \big(\Ve(\udl{e}_T), \bG(\updl{I}_T\zt) \big)_T
    =& a_T( \updl{e}_T, \updl{I}_T\zt )
      -\pd{ \udl{e}_{\uife}-\udl{e}_T, \bG(\updl{I}_T\zt)\udl{n} }_{\uife} \\
    &- \pd{ (\bI-\bK\bS_h) (\udl{e}_{\ife}-\udl{e}_T), \bG(\updl{I}_T\zt)\udl{n} }_{\ife}.
\end{align*}
Substituting this identity into \reff{Pro-L2-Err-1}, adding and subtracting
$s_T(\updl{e}_T,\updl{I}_T\zt)$ so as to reconstruct $a_h$, and using the
definition \reff{ah}, we deduce that
\begin{align}\label{Pro-L2-Err-2}
  \norm{\udl{e}_h}^2_\Omega =& \sum_{T\in\cTh}\Big[\big(\Ve(\udl{e}_T), \Vs(\zt)-\bG(\updl{I}_T\zt) \big)_T
       +\pd{ \udl{e}_{\uife}-\udl{e}_T, \big( \Vs(\zt)-\bG(\updl{I}_T\zt) \big)\udl{n}}_{\uife} \nn\\
    &+ \pd{ \udl{e}_{\ife}-\udl{e}_T, \Vs(\zt)\udl{n}}_{\ife}
       -\pd{ (\bI-\bK\bS_h) (\udl{e}_{\ife}-\udl{e}_T), \bG(\updl{I}_T\zt)\udl{n} }_{\ife} \nn\\
    &- s_T( \updl{e}_T, \updl{I}_T\zt ) \Big]
       +a_h( \updl{e}_h, \updl{I}_h\zt ).
\end{align}
The last term is estimated in Steps 4--5; the sum is estimated in Steps 2--3.

\smallskip\noindent{\em Step 2: rewriting the interface contributions of the sum.}
The third and fourth terms of \reff{Pro-L2-Err-2} and the interface
stabilization do not individually have the right structure, and we first
combine them.  Using $\av{\zt}=\zt-\frac12\jm{\zt}$ and
$-\jm{\zt}=\bK \Vs(\zt)\udl{n}$ on $\Gamma$, and the fact that
$\Pi_{\ife}$ is self-adjoint, the interface stabilization can be written as
\begin{align}\label{Pro-L2-Err-sTG}
  s_{T^{\Gamma}}(\updl{e}_T, \updl{I}_T\zt)
    =&\pd{ \bS_h (\udl{e}_{\ife}-\udl{e}_T) ,\bK \Vs(\zt)\udl{n} }_{\ife}
      +\pd{ \bS_h(I-\Pi_{\ife}) \udl{e}_T ,\bK \Vs(\zt)\udl{n} }_{\ife} \nn\\
    &+2\pd{ \bS_h\Pi_{\ife} (\udl{e}_{\ife}-\udl{e}_T) ,\zt-\Pi_T\zt }_{\ife}.
\end{align}
The first term of \reff{Pro-L2-Err-sTG} is precisely what is needed to
transform the third and fourth terms of \reff{Pro-L2-Err-2} into a single
pairing against the dual consistency error
$\big(\Vs(\zt)-\bG(\updl{I}_T\zt)\big)\udl{n}$, since
$\bI-\bK\bS_h$ appears in the fourth one.  We are thus left with
\begin{align}\label{Pro-L2-Err-3}
  \sum_{T\in\cTh}\Big[ & \big(\Ve(\udl{e}_T), \Vs(\zt)-\bG(\updl{I}_T\zt) \big)_T
       +\pd{ \udl{e}_{\uife}-\udl{e}_T, \big( \Vs(\zt)-\bG(\updl{I}_T\zt) \big)\udl{n}}_{\uife} \nn\\
    &+ \pd{ \udl{e}_{\ife}-\udl{e}_T, \Vs(\zt)\udl{n}}_{\ife}
       -\pd{(\bI-\bK\bS_h)(\udl{e}_{\ife}-\udl{e}_T),\bG(\updl{I}_T\zt)\udl{n}}_{\ife}
       -s_T( \updl{e}_T, \updl{I}_T\zt) \Big] \nn\\
    =&\sum_{T\in\cTh}\Big[ \big(\Ve(\udl{e}_T), \Vs(\zt)-\bG(\updl{I}_T\zt) \big)_T
       +\pd{ \udl{e}_{\uife}-\udl{e}_T, \big( \Vs(\zt)-\bG(\updl{I}_T\zt) \big)\udl{n}}_{\uife} \nn\\
    &+ \pd{(\bI-\bK\bS_h)(\udl{e}_{\ife}-\udl{e}_T), \big(\Vs(\zt)-\bG(\updl{I}_T\zt) \big)\udl{n}}_{\ife}
       -\pd{ \bS_h(I-\Pi_{\ife}) \udl{e}_T ,\bK \Vs(\zt)\udl{n} }_{\ife} \nn\\
    &-2\pd{ \bS_h\Pi_{\ife} (\udl{e}_{\ife}-\udl{e}_T) ,\zt-\Pi_T\zt }_{\ife}
       -s_{T^{\backslash \Gamma}}(\updl{e}_T, \updl{I}_T\zt) \Big] \nn\\
    :=& \diamondsuit.
\end{align}

\smallskip\noindent{\em Step 3: estimate of $\diamondsuit$.}
Every term of \reff{Pro-L2-Err-3} is now a pairing between a component of
$\updl{e}_T$, controlled by $\norme{\updl{e}_T}_{a_T}$, and an interpolation error
of the dual solution, which is of order $h_T$ by
\Lem{\ref{Lem-Interpolation-properties}}, \Thm{\ref{Thm-reconstruct-Approx}}
and \Lem{\ref{Lem-traction}} used with $m=1$.  The five estimates are the
exact analogues of \reff{Pro-energy-Err-4}, \reff{Pro-energy-Err-5},
\reff{Pro-energy-Err-9}, \reff{Pro-energy-Err-8} and
\reff{Pro-energy-Err-11}, with the roles of the primal and the dual solution
interchanged.  For the first one, splitting around $\pk\Vs(\zt)$ and
expanding according to the constitutive law as in Step 3 of the proof of
\Pro{\ref{Pro-energy-Err}},
\begin{align*}
  \big(\Ve(\udl{e}_T), \Vs(\zt)-\bG(\updl{I}_T\zt) \big)_T
    \leq& \norme{\updl{e}_T}_{a_T} \mu^{-\frac12}h\big(\mu\abs{\zt}_{H^2(T)^d}+\lambda\abs{\nabla\cdot\zt}_{H^1(T)}\big).
\end{align*}
For the second one, inserting $\Pi_{\uife}\udl{e}_T$ in the polynomial part
and using the trace inequalities,
\begin{align*}
  \pd{ \udl{e}_{\uife}-\udl{e}_T, \big( \Vs(\zt)-\bG(\updl{I}_T\zt) \big)\udl{n}}_{\uife}
    \leq& \norme{\updl{e}_T}_{a_T} \mu^{-\frac12}h\big(\mu\abs{\zt}_{H^2(T)^d}+\lambda\abs{\nabla\cdot\zt}_{H^1(T)}\big).
\end{align*}
The third one is bounded by moving half of the weight
${\bf QPQ}^T$ onto the test function, as in \reff{Pro-energy-Err-9}, and by
\reff{Lem-traction-3},
\begin{align*}
  &\pd{(\bI-\bK\bS_h)(\udl{e}_{\ife}-\udl{e}_T),
      \big(\Vs(\zt)-\bG(\updl{I}_T\zt) \big)\udl{n}}_{\ife} \\
  &\hspace{1cm}\leq  \norm{\bS_h^{\frac12}(\udl{e}_{\ife}-\udl{e}_T)}_{\ife}
      \sqrt{\frac{h}{\mu+\lambda}}\norm{\Vs(\zt)-\bG(\updl{I}_T\zt)}_{\ife} \\
  &\hspace{1cm}\leq \norme{\updl{e}_T}_{a_T} \mu^{-\frac12}h\big(\mu\abs{\zt}_{H^2(T)^d}
      +\lambda\abs{\nabla\cdot\zt}_{H^1(T)}\big).
\end{align*}
For the fourth one we use $\abs{\bK\bS_h}\leq 1$ from \reff{Lem-Sh-3}, the
orthogonality of $I-\Pi_{\ife}$ to $\mP^k(\ife)^d$,
\Lem{\ref{Lem-Fce-Korn}} and the trace inequality
\reff{Lem-Multip-trace-Ineq}, which give
\begin{align*}
  \pd{ \bS_h(I-\Pi_{\ife}) \udl{e}_T ,\bK \Vs(\zt)\udl{n} }_{\ife}
    =& \pd{ \bK\bS_h(I-\Pi_{\ife}) \udl{e}_T ,(I-\Pi_{\ife}) \Vs(\zt)\udl{n} }_{\ife} \\
    \lesssim& \sqrt{2\mu}\norm{\Ve(\udl{e}_T)}_T \mu^{-\frac12}h^{\frac12}\norm{(I-\Pi_{\ife}) \Vs(\zt)}_{\ife}\\
    \lesssim& \norme{\updl{e}_T}_{a_T} \mu^{-\frac12}h\big(\mu\abs{\zt}_{H^2(T)^d}+\lambda\abs{\nabla\cdot\zt}_{H^1(T)}\big).
\end{align*}
Finally, in the last two terms only the tangential part of $\bS_h$ survives,
because the normal component of $\zt-\Pi_T\zt$ vanishes on the faces by
\Rem{\ref{Rem-interpolants}}, and we conclude as in
\reff{Pro-energy-Err-11}--\reff{Pro-energy-Err-13} that
\begin{align*}
  2\pd{ \bS_h\Pi_{\ife} (\udl{e}_{\ife}-\udl{e}_T) ,\zt-\Pi_T\zt }_{\ife}
       +s_{T^{\backslash \Gamma}}(\updl{e}_T, \updl{I}_T\zt)
    \lesssim& \norme{\updl{e}_T}_{a_T} h\mu^{\frac12}\abs{\zt}_{H^2(T)^d}.
\end{align*}
Summing the five bounds over $T\in\cTh$, using the energy estimate
\reff{Pro-energy-Err-14} for $\norme{\updl{e}_h}_h$ and the regularity
\reff{Regularity-Asm} with $\udl{g}=\udl{e}_h$, we arrive at
\begin{align}\label{Pro-L2-Err-4}
  \diamondsuit\lesssim& \norme{\updl{e}_h}_h \mu^{-\frac12}h\big(\mu\abs{\zt}_{H^2(\cTh)^d}+\lambda\abs{\nabla\cdot\zt}_{H^1(\cTh)}\big) \nn\\
  \lesssim& \mu^{-1}h^{k+2}\big(\mu\abs{\udl{u}}_{H^{k+2}(\cTh)^d}+\lambda\abs{\nabla\cdot\udl{u}}_{H^{k+1}(\cTh) }\big)
          \norm{\udl{e}_h}_\Omega.
\end{align}

\smallskip\noindent{\em Step 4: the term $a_h(\updl{e}_h,\updl{I}_h\zt)$.}
Using the discrete problem \reff{discrete-form} with the test function
$\updl{I}_h\zt\in\updl{V}_{h0}$ and then the equation
$-\nabla\cdot\Vs(\udl{u})=\udl{f}$ together with an integration by parts, we
obtain
\begin{align}\label{Pro-L2-Err-6}
  a_h( \updl{e}_h, \updl{I}_h\zt ) =& \sum_{T\in \cTh} \big[ l_T(\Pi_T\zt) -a_T(\updl{I}_T\udl{u}, \updl{I}_T\zt)
       -s_T(\updl{I}_T\udl{u}, \updl{I}_T\zt)\big] \nn\\
    =& \sum_{T\in \cTh} \Big[ (\udl{f},\Pi_T\zt-\zt)_T +(\Vs(\udl{u}), \Ve(\zt))_T
       -\frac12\pd{ \Vs(\udl{u})\udl{n}, \jm{\zt} }_{\ife}
       -a_T(\updl{I}_T\udl{u}, \updl{I}_T\zt) \nn\\
    &- s_{T^{\Gamma}}(\updl{I}_T\udl{u}, \updl{I}_T\zt)
       -s_{T^{\backslash \Gamma}}(\updl{I}_T\udl{u}, \updl{I}_T\zt) \Big].
\end{align}
The point is now that the leading parts of the four middle terms cancel.
Indeed, expanding the interface stabilization as in Step 2, this time for the
pair $(\updl{I}_T\udl{u},\updl{I}_T\zt)$, gives
\begin{align}\label{Pro-L2-Err-5}
  s_{T^{\Gamma}}(\updl{I}_T\udl{u}, \updl{I}_T\zt)
    =& \pd{ \bS_h\Pi_{\ife} (\udl{u}-\Pi_T\udl{u}) ,\bK \Vs(\zt)\udl{n} }_{\ife}
       +\frac12 \pd{ \bS_h(I-\Pi_{\ife}) \jm{\udl{u}} ,\bK \Vs(\zt)\udl{n} }_{\ife} \nn\\
    &+ \frac12 \pd{ \bS_h\bK \Vs(\udl{u})\udl{n} ,\bK \Vs(\zt)\udl{n} }_{\ife}
       +2\pd{ \bS_h\Pi_{\ife} (\udl{u}-\Pi_T\udl{u}) ,\zt-\Pi_T\zt }_{\ife} \nn\\
    &+ \pd{ \bS_h(I-\Pi_{\ife}) \jm{\udl{u}} ,\zt-\Pi_T\zt }_{\ife}
       +\pd{ \bS_h\bK \Vs(\udl{u})\udl{n} ,\zt-\Pi_T\zt }_{\ife},
\end{align}
while, using $-\jm{\zt}=\bK\Vs(\zt)\udl{n}$ and the splitting
$\bI=(\bI-\bK\bS_h)+\bK\bS_h$,
\begin{align*}
  (\Vs(\udl{u}), \Ve(\zt))_T -\frac12\pd{ \Vs(\udl{u})\udl{n}, \jm{\zt} }_{\ife}
    =& (\Vs(\udl{u}), \Ve(\zt))_T
      +\frac12\pd{ (\bI-\bK\bS_h)\bK \Vs(\udl{u})\udl{n}, \Vs(\zt)\udl{n} }_{\ife} \\
    &+\frac12\pd{ \bS_h\bK\Vs(\udl{u})\udl{n}, \bK \Vs(\zt)\udl{n} }_{\ife} \\
    =& a_T(\udl{u}, \zt)
      +\frac12\pd{ \bS_h\bK\Vs(\udl{u})\udl{n}, \bK \Vs(\zt)\udl{n} }_{\ife}.
\end{align*}
Inserting these two identities in \reff{Pro-L2-Err-6} and using
\Lem{\ref{Lem-reconstruct-Difference}} to replace
$a_T(\udl{u},\zt)-a_T(\updl{I}_T\udl{u},\updl{I}_T\zt)$, all the terms
containing $\Vs(\udl{u})\udl{n}$ paired with $\bK\Vs(\zt)\udl{n}$ cancel and
we are left with
\begin{align}\label{Pro-L2-Err-6b}
  (\Vs(\udl{u}), \Ve(\zt))_T &-\frac12\pd{ \Vs(\udl{u})\udl{n}, \jm{\zt} }_{\ife}
       -a_T(\updl{I}_T\udl{u}, \updl{I}_T\zt) -s_{T^{\Gamma}}(\updl{I}_T\udl{u}, \updl{I}_T\zt)  \nn\\
    =& a_T(\udl{u}-\updl{I}_T\udl{u}, \zt-\updl{I}_T\zt)
       -2\pd{ \bS_h\Pi_{\ife} (\udl{u}-\Pi_T\udl{u}) ,\zt-\Pi_T\zt }_{\ife}   \nn\\
    &+ \pd{ \bS_h\bK(\udl{u}-\Pi_T\udl{u}), \big( \bG(\updl{I}_T\zt)-\Pi_{\ife}\Vs(\zt) \big)\udl{n} }_{\ife} \nn\\
    &+ \pd{\bS_h\bK \big( \bG(\updl{I}_T\udl{u})-\Vs(\udl{u}) \big)\udl{n}, \zt-\Pi_T\zt }_{\ife}
       -\frac12 \pd{ \bS_h(I-\Pi_{\ife}) \jm{\udl{u}} ,\bK \Vs(\zt)\udl{n} }_{\ife}  \nn\\
    &- \pd{ \bS_h(I-\Pi_{\ife}) \jm{\udl{u}} ,\zt-\Pi_T\zt }_{\ife}.
\end{align}

\smallskip\noindent{\em Step 5: estimate of the six remainders.}
The first term of \reff{Pro-L2-Err-6b} is estimated by the Cauchy--Schwarz
inequality for the bilinear form \reff{aT-extended} and by
\Lem{\ref{Lem-aT-consistency}}, applied with $m=k+1$ to $\udl{u}$ and with
$m=1$ to $\zt$:
\begin{align*}
  a_T(\udl{u}-\updl{I}_T\udl{u}, \zt-\updl{I}_T\zt)
    \leq& a_T(\udl{u}-\updl{I}_T\udl{u},\udl{u}-\updl{I}_T\udl{u})^{\frac12}\,
          a_T(\zt-\updl{I}_T\zt,\zt-\updl{I}_T\zt)^{\frac12} \\
    \lesssim& \mu^{-1}h^{k+2} \big(\mu \abs{\udl{u}}_{H^{k+2}(T)^d}+\lambda\abs{\nabla\cdot\udl{u}}_{H^{k+1}(T)}\big)
             \big(\mu \abs{\zt}_{H^2(T)^d}+\lambda\abs{\nabla\cdot\zt}_{H^1(T)}\big).
\end{align*}
The second term only involves the tangential part of $\bS_h$, as in
\reff{Pro-energy-Err-11}, and the trace inequality
\reff{Lem-Multip-trace-Ineq} gives
\begin{align*}
  \pd{ \bS_h\Pi_{\ife} (\udl{u}-\Pi_T\udl{u}) ,\zt-\Pi_T\zt }_{\ife}
    =& \pd{ \bS_{h,t}\Pi_{\ife} (\udl{u}-\Pi_T\udl{u}) ,\zt-\Pi_T\zt }_{\ife} \\
    \leq& \mu\big(h^{-1}\norm{\udl{u}-\Pi_T\udl{u}}_T +\norm{\nabla(\udl{u}-\Pi_T\udl{u})}_T\big) \\
    &     \big(h^{-1}\norm{\zt-\Pi_T\zt}_T +\norm{\nabla(\zt-\Pi_T\zt)}_T\big) \\
    \leq& \mu^{-1}h^{k+2} \mu \abs{\udl{u}}_{H^{k+2}(T)^d}\, \mu \abs{\zt}_{H^2(T)^d}.
\end{align*}
The third and fourth terms are estimated with $\abs{\bS_h\bK}\leq1$, the trace
inequality \reff{Lem-Multip-trace-Ineq} and
\reff{Lem-traction-1}, \reff{Lem-traction-3} of \Lem{\ref{Lem-traction}},
\begin{align*}
  \pd{ \bS_h\bK(\udl{u}-\Pi_T\udl{u}), \big( \bG(\updl{I}_T\zt)-\Pi_{\ife}\Vs(\zt) \big)\udl{n} }_{\ife}
    &+ \pd{\bS_h\bK \big( \bG(\updl{I}_T\udl{u})-\Vs(\udl{u}) \big)\udl{n}, \zt-\Pi_T\zt }_{\ife} \\
    &\hspace{-5cm}\leq \mu^{\frac12}h^{-\frac12}\norm{\udl{u}-\Pi_T\udl{u}}_{\ife}
          \sqrt{\frac{h}{\mu+\lambda}}\norm{\bG(\updl{I}_T\zt)-\Pi_{\ife}\Vs(\zt)}_{\ife}\\
    &\hspace{-4.5cm}+ \sqrt{\frac{h}{\mu+\lambda}} \norm{\bG(\updl{I}_T\udl{u})-\Vs(\udl{u})}_{\ife}
         \mu^{\frac12}h^{-\frac12}\norm{\zt-\Pi_T\zt}_{\ife} \\
    &\hspace{-5cm}\leq \mu^{-1}h^{k+2}\big(\mu\abs{\udl{u}}_{H^{k+2}(T)^d}+\lambda\abs{\nabla\cdot\udl{u}}_{H^{k+1}(T)}\big)
      \big(\mu \abs{\zt}_{H^2(T)^d}+\lambda\abs{\nabla\cdot\zt}_{H^1(T)}\big).
\end{align*}
The two last terms are the ones in which the jump of the exact solution
appears; they are treated with the orthogonality
$\pd{(I-\Pi_{\ife})\xi,v_h}_{\ife}=0$ for all $v_h\in \mP^k(\ife)^d$, with
$\abs{\bK\bS_h}\leq 1$ and, for the very last one, with
$\jm{\udl{u}}=-\bK\Vs(\udl{u})\udl{n}$, which gives
\begin{align*}
  \pd{ \bS_h(I-\Pi_{\ife}) \jm{\udl{u}} ,\bK \Vs(\zt)\udl{n} }_{\ife}
    =& \pd{ \bK\bS_h(I-\Pi_{\ife}) \jm{\udl{u}} ,(I-\Pi_{\ife})\Vs(\zt)\udl{n} }_{\ife} \\
    \leq& \norm{(I-\Pi_{\ife})\udl{u}}_{\ife} \norm{(I-\Pi_{\ife})\Vs(\zt)}_{\ife} \\
    \lesssim& \mu^{-1}h^{k+2} \mu\abs{\udl{u}}_{H^{k+2}(T)^d}
          \big(\mu \abs{\zt}_{H^2(T)^d}+\lambda\abs{\nabla\cdot\zt}_{H^1(T)}\big),
\end{align*}
and
\begin{align*}
  \pd{ \bS_h(I-\Pi_{\ife}) \jm{\udl{u}} ,\zt-\Pi_T\zt }_{\ife}
    =&\pd{ \bS_h\bK(I-\Pi_{\ife}) \Vs(\udl{u})\udl{n} ,\zt-\Pi_T\zt }_{\ife} \\
    \leq& h^{\frac12}\norm{(I-\Pi_{\ife})\Vs(\udl{u})}_{\ife} h^{-\frac12}\norm{\zt-\Pi_T\zt}_{\ife} \\
    \lesssim& \mu^{-1}h^{k+2}\big(\mu\abs{\udl{u}}_{H^{k+2}(T)^d}+\lambda\abs{\nabla\cdot\udl{u}}_{H^{k+1}(T)}\big)
         \mu \abs{\zt}_{H^2(T)^d}.
\end{align*}

\smallskip\noindent{\em Step 6: conclusion.}
Inserting the six bounds of Step 5 into \reff{Pro-L2-Err-6}, and estimating
the two remaining terms of \reff{Pro-L2-Err-6} by the orthogonality of
$\Pi_T$, by the Cauchy--Schwarz inequality for the positive semi-definite form
$s_{T^{\backslash \Gamma}}$ and by \Lem{\ref{Lem-Interpolation-properties}},
we obtain with the regularity assumption \reff{Regularity-Asm} that
\begin{align}\label{Pro-L2-Err-7}
  a_h( \updl{e}_h, \updl{I}_h\zt ) \leq& \sum_{T\in \cTh} \Big[ (\udl{f},\Pi_T\zt-\zt)_T
      +s_{T^{\backslash \Gamma}}(\updl{I}_T\udl{u}, \updl{I}_T\zt) \Big]\nn\\
    &+ \mu^{-1}h^{k+2}\big(\mu\abs{\udl{u}}_{H^{k+2}(\cTh)^d}+\lambda\abs{\nabla\cdot\udl{u}}_{H^{k+1}(\cTh)}\big)
          \big(\mu\abs{\zt}_{H^2(\cTh)^d}+\lambda\abs{\nabla\cdot\zt}_{H^1(\cTh)}\big) \nn\\
    \lesssim& \sum_{T\in \cTh} \Big[ (\udl{f}-\Pi_T\udl{f},\Pi_T\zt-\zt)_T
      +s_{T^{\backslash \Gamma}}^{\frac12}(\updl{I}_T\udl{u}, \updl{I}_T\udl{u})
       s_{T^{\backslash \Gamma}}^{\frac12}(\updl{I}_T\zt, \updl{I}_T\zt) \Big]\nn\\
    &+ \mu^{-1}h^{k+2}\big(\mu\abs{\udl{u}}_{H^{k+2}(\cTh)^d}+\lambda\abs{\nabla\cdot\udl{u}}_{H^{k+1}(\cTh)}\big)
          \norm{\udl{e}_h}_\Omega \nn\\
    \lesssim& \sum_{T\in \cTh} \Big[ h^{k+2}\abs{\udl{f}}_{H^k(T)^d} \abs{\zt}_{H^2(T)^d}
          +h^{k+2} \abs{\udl{u}}_{H^{k+2}(T)^d} \mu\abs{\zt}_{H^2(T)^d} \Big]\nn\\
    &+ \mu^{-1}h^{k+2}\big(\mu\abs{\udl{u}}_{H^{k+2}(\cTh)^d}+\lambda\abs{\nabla\cdot\udl{u}}_{H^{k+1}(\cTh)}\big)
          \norm{\udl{e}_h}_\Omega \nn\\
    \lesssim& \mu^{-1}h^{k+2}\big(\mu\abs{\udl{u}}_{H^{k+2}(\cTh)^d}+\lambda\abs{\nabla\cdot\udl{u}}_{H^{k+1}(\cTh) }\big)
          \norm{\udl{e}_h}_\Omega.
\end{align}
Combining \reff{Pro-L2-Err-2}, \reff{Pro-L2-Err-4} and \reff{Pro-L2-Err-7},
and dividing by $\norm{\udl{e}_h}_\Omega$, deduces \reff{Pro-L2-Err-0}.
\Endproof

\Section{Numerical experiments.}
\label{sec-numerics}

In this section we discuss the implementation of the method and we report
two-dimensional numerical experiments that assess the error estimates of
\Pro{\ref{Pro-energy-Err}} and \Pro{\ref{Pro-L2-Err}}, i.e.\ the convergence
rates $h^{k+1}$ and $h^{k+2}$ of the energy- and $L^2$-errors, their
robustness with respect to the compliancy parameters $\alpha\ge0$,
$\beta\ge0$, and their robustness in the quasi-incompressible limit
$\lambda\to+\infty$.  Two test cases are considered, in increasing order of
difficulty.  The first one, in \S\ref{sec-num1}, is a manufactured solution
across a {\em straight} interface; it satisfies the interface conditions
exactly for every $\alpha\geq0$, $\beta\geq0$, it is set on general polytopal
meshes, and it is the configuration for which the analysis of \S 4--\S 5
applies verbatim.  It is used to assess the convergence rates and the
behaviour in the quasi-incompressible limit.  The
second one, in \S\ref{sec-num2}, is the elastic inclusion problem of
\cite[Sect.~5.1]{HH04}, whose exact solution we extend here to a compliant
interface.  It is used to assess the robustness with respect to the compliancy
parameters over sixteen orders of magnitude.

\subsection{Implementation.}
\label{sec-implementation}

As is customary for hybrid high-order methods \cite{CDE18}, an important step
in the implementation is the choice of a basis for each of the polynomial
spaces appearing in the construction.  For all $T\in\cTh$ and all
$l\in\{k,k+1\}$ we use the hierarchical basis of $\mP^l(T)$ obtained by
$L^2(T)$-orthonormalisation (Cholesky factorisation of the mass matrix) of the
monomials in the scaled variables
$\xi_{T}=(\udl{x}-\udl{x}_T)/h_T$, where $\udl{x}_T$ is the barycentre of $T$;
a basis of $\mP^l(T)^d$, resp.\ of the space $\mP^l(T)^{d\times d}$ of
symmetric-matrix-valued polynomials, is then obtained by taking the Cartesian
product with the canonical basis of $\mR^d$, resp.\ with an orthonormal basis
of the symmetric $d\times d$ matrices.  Faces are handled in the same way,
with respect to the arclength parameter.  Two remarks are in order.  First,
the bases being $L^2$-orthonormal, the mass matrices in
\reff{gradient-reconstruct} and in the local $L^2$-projections are identity
matrices, and the $L^2$-norm of a cell unknown is the Euclidean norm of its
coefficients.  Second, the basis being hierarchical, the basis of $\mP^k(T)$
needed in \reff{gradient-reconstruct} and in \reff{stabilization-term} is
obtained from that of $\mP^{k+1}(T)$ by discarding the highest-order
functions.

For every mesh cell, the reconstruction \reff{gradient-reconstruct} amounts to
solving a linear system of size $\frac12 d(d+1)\dim\mP^{k}(T)$, whose matrix
is the identity for a cell such that $\ife=\emptyset$ and the identity plus
the contribution of the second term in the left-hand side of
\reff{gradient-reconstruct} otherwise; the discrete divergence is recovered as
$D(\updl{v}_T)=tr\bE(\updl{v}_T)$, which is \reff{divergence-reconstruct}.
The local matrices of $a_T+s_T$ are symmetric positive semi-definite with a
three-dimensional kernel, that of the rigid-body motions (checked numerically,
cell by cell, on the coarsest mesh of each of the three families used below
and for $k\in\{1,2,3\}$).  The cell unknowns are eliminated by static
condensation, so that the global linear system is symmetric positive definite
and only involves the face unknowns; its size is
$N_{\mathrm{dof}}=d\dim\mP^k(F)\,\mathrm{card}(\cF_h\setminus\cF_h^b)$, that
is $2(k+1)\,\mathrm{card}(\cF_h\setminus\cF_h^b)$ for $d=2$.  It is solved by
a sparse direct method.  We recall from
\Rem{\ref{Rem-interface-unknown}} that the unknowns attached to the interface
faces are single valued, exactly as the other face unknowns, so that the
implementation of the interface conditions is entirely local: it only affects
the local matrices of the cells that touch $\Gamma$.

Concerning numerical integration, all the integrals over straight (polygonal)
cells are evaluated by decomposing the cell into triangles and using a
collapsed Gauss--Jacobi rule with $k+6$ nodes per direction, which is exact
for polynomials of total degree $2k+11$; face integrals use Gauss--Legendre
rules with $k+6$ nodes, of the same accuracy.  This is well beyond what the
products appearing in \reff{gradient-reconstruct}--\reff{stabilization-term}
require.  On the curved meshes of \S\ref{sec-num2}, the integrals over the
annular sectors and over the circular faces are computed with tensor-product
Gauss rules in the polar variables $(r,\theta)$ with $k+6$ nodes per
direction; enriching these rules further modifies the errors reported below by
less than $10^{-6}$ in relative value, so that the quadrature error is
negligible.

The Dirichlet condition is non-homogeneous in both test cases; it is enforced
strongly by prescribing $\udl{v}_F=\Pi_F(\udl{u})$ for all $F\in\cF_h^b$.
Since the interpolate $\updl{I}_h(\udl{u})$ satisfies the same constraint, the
error analysis of \S 5 applies verbatim.  Setting
$\updl{e}_h:=\updl{u}_h-\updl{I}_h(\udl{u})$, we monitor the three error
measures
\begin{equation}\label{errors}
  E_{e}:=\norme{\updl{e}_h}_{e},
  \qquad
  E_{h}:=\norme{\updl{e}_h}_{h}=a_h(\updl{e}_h,\updl{e}_h)^{1/2},
  \qquad
  E_{0}:=\norm{\udl{e}_h}_{\Omega},
\end{equation}
cf.\ \reff{Energy-norm}, \reff{Energy-norm-h} and \Pro{\ref{Pro-L2-Err}};
recall that $E_{e}\lesssim E_{h}$ by \reff{Thm-Stability-coercivity}.
\Pro{\ref{Pro-energy-Err}} predicts $E_{e}=O(h^{k+1})$ and
\Pro{\ref{Pro-L2-Err}} predicts $E_{0}=O(h^{k+2})$.  The error
$\norm{\udl{u}-\udl{u}_h}_{\Omega}$ differs from $E_{0}$ by at most
$Ch^{k+2}\abs{\udl{u}}_{H^{k+2}(\cTh)^d}$ and was found to behave identically
(in all the computations reported below, the ratio of the two quantities lies
between $1.01$ and $1.11$); it is not reported.  As a first, elementary
validation of the implementation, we checked that, for $k\in\{1,2,3\}$, the
discrete problem reproduces to machine precision every rigid-body translation,
on all three mesh families used below and for any value of $\alpha$ and
$\beta$, as well as, on the two families with straight faces, every piecewise
affine displacement field satisfying the interface conditions of
\reff{elasticity-interface} with $\bK=\bf{0}$.  Both restrictions are
intrinsic: on a circular face the trace of an affine field is not a polynomial
in the arclength, and for $\bK\neq\bf{0}$ one has
$s_{T^{\Gamma}}(\updl{I}_T\udl{u},\cdot)\neq 0$, which is precisely the
consistency error estimated in Step 6 of the proof of
\Pro{\ref{Pro-energy-Err}}.  All the computations reported below were
performed with a plain {\tt Python}/{\tt NumPy} implementation written for
verification purposes; since no attempt was made at optimising it, we do not
report computational times.

\subsection{Test case 1: straight interface and general meshes.}
\label{sec-num1}

We first let $\Omega=(0,1)^2$, $\Gamma=\{x_1=\frac12\}$,
$\Omega_1=\{x_1<\frac12\}$ and $\udl{n}=\udl{e}_1$, and we manufacture an
exact solution which satisfies the interface conditions of
\reff{elasticity-interface} exactly for every $\alpha\geq0$, $\beta\geq0$.
Let $\psi(\udl{x})=\sin(\pi x_1+\frac{\pi}{4})\sin(2\pi x_2)$ and set
\begin{equation}\label{num-u2}
  \udl{u}_2=\underline{\mathrm{curl}}\,\psi +\frac{1}{2\lambda_2}\udl{x}
    =\Big(\p_2\psi+\frac{x_1}{2\lambda_2},\
      -\p_1\psi+\frac{x_2}{2\lambda_2}\Big),
\end{equation}
so that $\nabla\cdot\udl{u}_2=\lambda_2^{-1}$ and
$\udl{f}_2=5\pi^2\mu_2\,\underline{\mathrm{curl}}\,\psi$; the two different
frequencies in $\psi$ make the tangential traction on $\Gamma$ non-zero, so
that $\alpha$ is active, and the affine term is the two-dimensional analogue
of the one used in \cite[Eq.~(71)]{DE15}, which keeps
$\lambda\abs{\nabla\cdot\udl{u}}$ bounded, and non-zero, as
$\nu\to\frac12$.  Let then
$\udl{t}:=\Vs_2(\udl{u}_2)\udl{n}|_{\Gamma}$ and define, on $\Gamma$,
\begin{equation}\label{num-VW}
  \udl{V}:=-\bK\udl{t},
  \qquad
  \udl{R}:=\udl{t}-\Vs_1(\udl{u}_2)\udl{n}|_{\Gamma},
  \qquad
  W_1:=\frac{R_1-\lambda_1\p_2V_2}{2\mu_1+\lambda_1},
  \qquad
  W_2:=\frac{R_2}{\mu_1}-\p_2V_1,
\end{equation}
where $\Vs_i$ denotes the constitutive law of $\Omega_i$, and finally
\begin{equation}\label{num-u1}
  \udl{u}_1:=\udl{u}_2+\udl{V} +\big(x_1-\tfrac12\big)\udl{W}
    -\tfrac12\big(x_1-\tfrac12\big)^2(\p_2W_2)\,\udl{e}_1,
  \qquad
  \udl{f}_1:=-\nabla\cdot\Vs_1(\udl{u}_1).
\end{equation}
By construction $\jm{\udl{u}}=\udl{V}=-\bK\Vs(\udl{u})\udl{n}$ and
$\big[\Vs(\udl{u})\udl{n}\big]=\udl{0}$ on $\Gamma$: the quadratic term in
\reff{num-u1} vanishes together with its normal derivative on $\Gamma$, and is
included so that $\nabla\cdot\udl{u}_1=O(\lambda_1^{-1})$ irrespective of
$\alpha$ and $\beta$.  The functions $\udl{u}_i$ and $\udl{f}_i$ are computed
by symbolic differentiation.  Unless stated otherwise we take $E_1=1$,
$\nu_1=0.25$, $E_2=10$, $\nu_2=0.3$ and $\alpha=\beta=5\cdot10^{-2}$.

Two mesh families, both fitted to $\Gamma$, are considered, see
Fig.~\ref{fig-meshes}(b)--(c): a matching triangular family, and a polygonal
family obtained by shifting every other column of a Cartesian mesh by half a
cell in the $x_2$-direction.  The cells of the latter have up to six faces and
the mesh has hanging nodes, a situation which the framework of \S 3 covers.
Both families are indexed by $\ell\in\{0,1,2,3\}$, the coarsest meshes being
built on a $4\times4$ Cartesian mesh.  The errors $E_{e}$ and $E_{0}$ are
reported in Table~\ref{tab-conv1} and Fig.~\ref{fig-conv1}: the rates $k+1$
for $E_{e}$ and $k+2$ for $E_{0}$ predicted by
\Pro{\ref{Pro-energy-Err}}--\Pro{\ref{Pro-L2-Err}} are observed for
$k\in\{1,2,3\}$ on both families.  At comparable $h$, the polygonal meshes
carry about half as many cells as the triangular ones and yield errors larger
by a factor $1.3$ ($k=1$), $1.7$ ($k=2$) and $2.2$ ($k=3$) for $E_{e}$, and
$1.2$, $2.4$ and $3.7$ for $E_{0}$.

Finally, Table~\ref{tab-lock} reports the behaviour of the method in the
quasi-incompressible limit, obtained by letting
$\nu_1=\nu_2=\nu\to\frac12$ at fixed $E_1$, $E_2$, so that
$\lambda_1,\lambda_2\to+\infty$ while
$\lambda\norm{\nabla\cdot\udl{u}}_{\infty}$ stays bounded.  The errors
$E_{e}$ and $E_{0}$, as well as the observed rates, are seen to be essentially
independent of $\lambda$: from $\nu=0.49$ to $\nu=0.499999$, that is over four
orders of magnitude in $\lambda$, they vary by less than $0.8\,\%$ ($E_{e}$)
and $2.3\,\%$ ($E_{0}$); the larger variation observed between $\nu=0.3$ and
$\nu=0.49$ (a factor $1.4$, resp.\ $1.5$--$1.7$) is that of the exact solution
itself, as reflected by $\lambda\norm{\nabla\cdot\udl{u}}_{\infty}$.  The
method is therefore locking free, in agreement with the fact that the
right-hand sides of \reff{Pro-energy-Err-0} and \reff{Pro-L2-Err-0} involve
$\lambda$ only through the product
$\lambda\abs{\nabla\cdot\udl{u}}_{H^{k+1}(\cTh)}$.

\begin{rem}[On the $\lambda$-weighted terms of $\norme{\cdot}_h$]\label{Rem-lambda-norm}
The error $E_{h}$ measured in the norm $\norme{\cdot}_h$ of
\reff{Energy-norm-h} behaves differently: as shown in the last columns of
Table~\ref{tab-lock}, it converges at the optimal rate $h^{k+1}$ but with a
constant growing as $\lambda^{1/2}$.  
The displacement errors themselves are, as reported above, uniform in
$\lambda$.
\end{rem}

\begin{table}[htbp]
\centering
\small
\caption{Test case 1 (straight interface, $\alpha=\beta=5\cdot10^{-2}$). Errors \reff{errors} and observed convergence rates.}
\label{tab-conv1}
\begin{tabular}{rrrr|rr|rr}
\hline
$\mathrm{card}(\cTh)$ & $\mathrm{card}(\cF_h)$ & $N_{\mathrm{dof}}$ & $h$ & $E_{e}$ & rate & $E_{0}$ & rate\\
\hline
\multicolumn{8}{l}{\emph{triangular meshes}, $k=1$}\\
$32$ & $56$ & $160$ & $0.354$ & $4.29\cdot 10^{1}$ & -- & $2.60$ & --\\
$128$ & $208$ & $704$ & $0.177$ & $1.15\cdot 10^{1}$ & $1.90$ & $3.57\cdot 10^{-1}$ & $2.86$\\
$512$ & $800$ & $2{,}944$ & $0.088$ & $3.07$ & $1.91$ & $4.82\cdot 10^{-2}$ & $2.89$\\
$2{,}048$ & $3{,}136$ & $12{,}032$ & $0.044$ & $7.87\cdot 10^{-1}$ & $1.96$ & $6.25\cdot 10^{-3}$ & $2.95$\\
\hline
\multicolumn{8}{l}{\emph{triangular meshes}, $k=2$}\\
$32$ & $56$ & $240$ & $0.354$ & $8.89$ & -- & $2.94\cdot 10^{-1}$ & --\\
$128$ & $208$ & $1{,}056$ & $0.177$ & $1.36$ & $2.71$ & $2.43\cdot 10^{-2}$ & $3.59$\\
$512$ & $800$ & $4{,}416$ & $0.088$ & $1.74\cdot 10^{-1}$ & $2.97$ & $1.58\cdot 10^{-3}$ & $3.95$\\
$2{,}048$ & $3{,}136$ & $18{,}048$ & $0.044$ & $2.20\cdot 10^{-2}$ & $2.99$ & $1.00\cdot 10^{-4}$ & $3.98$\\
\hline
\multicolumn{8}{l}{\emph{triangular meshes}, $k=3$}\\
$32$ & $56$ & $320$ & $0.354$ & $2.02$ & -- & $4.50\cdot 10^{-2}$ & --\\
$128$ & $208$ & $1{,}408$ & $0.177$ & $1.19\cdot 10^{-1}$ & $4.08$ & $1.35\cdot 10^{-3}$ & $5.06$\\
$512$ & $800$ & $5{,}888$ & $0.088$ & $7.64\cdot 10^{-3}$ & $3.97$ & $4.38\cdot 10^{-5}$ & $4.95$\\
$2{,}048$ & $3{,}136$ & $24{,}064$ & $0.044$ & $4.82\cdot 10^{-4}$ & $3.99$ & $1.39\cdot 10^{-6}$ & $4.98$\\
\hline
\multicolumn{8}{l}{\emph{polygonal meshes}, $k=1$}\\
$18$ & $55$ & $152$ & $0.354$ & $5.54\cdot 10^{1}$ & -- & $3.60$ & --\\
$68$ & $205$ & $688$ & $0.177$ & $1.62\cdot 10^{1}$ & $1.78$ & $4.96\cdot 10^{-1}$ & $2.86$\\
$264$ & $793$ & $2{,}912$ & $0.088$ & $4.13$ & $1.97$ & $6.12\cdot 10^{-2}$ & $3.02$\\
$1{,}040$ & $3{,}121$ & $11{,}968$ & $0.044$ & $1.04$ & $1.99$ & $7.61\cdot 10^{-3}$ & $3.01$\\
\hline
\multicolumn{8}{l}{\emph{polygonal meshes}, $k=2$}\\
$18$ & $55$ & $228$ & $0.354$ & $1.55\cdot 10^{1}$ & -- & $7.45\cdot 10^{-1}$ & --\\
$68$ & $205$ & $1{,}032$ & $0.177$ & $2.19$ & $2.83$ & $5.58\cdot 10^{-2}$ & $3.74$\\
$264$ & $793$ & $4{,}368$ & $0.088$ & $2.89\cdot 10^{-1}$ & $2.92$ & $3.71\cdot 10^{-3}$ & $3.91$\\
$1{,}040$ & $3{,}121$ & $17{,}952$ & $0.044$ & $3.72\cdot 10^{-2}$ & $2.96$ & $2.38\cdot 10^{-4}$ & $3.96$\\
\hline
\multicolumn{8}{l}{\emph{polygonal meshes}, $k=3$}\\
$18$ & $55$ & $304$ & $0.354$ & $3.34$ & -- & $1.14\cdot 10^{-1}$ & --\\
$68$ & $205$ & $1{,}376$ & $0.177$ & $2.45\cdot 10^{-1}$ & $3.77$ & $4.48\cdot 10^{-3}$ & $4.66$\\
$264$ & $793$ & $5{,}824$ & $0.088$ & $1.66\cdot 10^{-2}$ & $3.88$ & $1.58\cdot 10^{-4}$ & $4.83$\\
$1{,}040$ & $3{,}121$ & $23{,}936$ & $0.044$ & $1.07\cdot 10^{-3}$ & $3.95$ & $5.21\cdot 10^{-6}$ & $4.93$\\
\hline
\end{tabular}
\end{table}
\begin{table}[htbp]
\centering
\small
\caption{Test case 1. Quasi-incompressible limit on the triangular mesh family ($\alpha=\beta=5\cdot10^{-2}$, $\lambda_1=\lambda_2/10$): errors on the third mesh of the family ($\mathrm{card}(\cTh)=512$, $h=8.84\cdot10^{-2}$) and rates observed on the last refinement. The errors $E_{e}$ and $E_{0}$ are essentially independent of $\lambda$, whereas $E_{h}$ grows as $\lambda^{1/2}$, cf.\ \Rem{\ref{Rem-lambda-norm}}.}
\label{tab-lock}
\begin{tabular}{rrr|rr|rr|rr}
\hline
$\nu$ & $\lambda_2$ & $\lambda\norm{\nabla\cdot\udl{u}}_{\infty}$ & $E_{e}$ & rate & $E_{0}$ & rate & $E_{h}$ & rate\\
\hline
\multicolumn{9}{l}{$k=1$}\\
$0.3$ & $5.77$ & $33.8$ & $2.87$ & $1.91$ & $4.49\cdot 10^{-2}$ & $2.89$ & $2.28$ & $1.91$\\
$0.49$ & $1.64\cdot 10^{2}$ & $67.8$ & $2.06$ & $1.95$ & $2.67\cdot 10^{-2}$ & $2.93$ & $1.92$ & $1.94$\\
$0.4999$ & $1.67\cdot 10^{4}$ & $70.1$ & $2.06$ & $1.95$ & $2.62\cdot 10^{-2}$ & $2.93$ & $1.37\cdot 10^{1}$ & $1.99$\\
$0.499999$ & $1.67\cdot 10^{6}$ & $70.1$ & $2.06$ & $1.95$ & $2.62\cdot 10^{-2}$ & $2.93$ & $1.36\cdot 10^{2}$ & $1.99$\\
\hline
\multicolumn{9}{l}{$k=2$}\\
$0.3$ & $5.77$ & $33.8$ & $1.65\cdot 10^{-1}$ & $2.97$ & $1.50\cdot 10^{-3}$ & $3.95$ & $9.57\cdot 10^{-2}$ & $2.95$\\
$0.49$ & $1.64\cdot 10^{2}$ & $67.8$ & $1.16\cdot 10^{-1}$ & $2.98$ & $9.44\cdot 10^{-4}$ & $3.96$ & $7.76\cdot 10^{-2}$ & $2.96$\\
$0.4999$ & $1.67\cdot 10^{4}$ & $70.1$ & $1.16\cdot 10^{-1}$ & $2.98$ & $9.24\cdot 10^{-4}$ & $3.96$ & $5.10\cdot 10^{-1}$ & $3.00$\\
$0.499999$ & $1.67\cdot 10^{6}$ & $70.1$ & $1.16\cdot 10^{-1}$ & $2.98$ & $9.24\cdot 10^{-4}$ & $3.96$ & $5.07$ & $3.00$\\
\hline
\multicolumn{9}{l}{$k=3$}\\
$0.3$ & $5.77$ & $33.8$ & $7.25\cdot 10^{-3}$ & $3.97$ & $4.20\cdot 10^{-5}$ & $4.95$ & $3.02\cdot 10^{-3}$ & $3.95$\\
$0.49$ & $1.64\cdot 10^{2}$ & $67.8$ & $5.06\cdot 10^{-3}$ & $3.99$ & $2.75\cdot 10^{-5}$ & $4.97$ & $2.34\cdot 10^{-3}$ & $3.96$\\
$0.4999$ & $1.67\cdot 10^{4}$ & $70.1$ & $5.03\cdot 10^{-3}$ & $3.99$ & $2.69\cdot 10^{-5}$ & $4.98$ & $1.44\cdot 10^{-2}$ & $4.00$\\
$0.499999$ & $1.67\cdot 10^{6}$ & $70.1$ & $5.03\cdot 10^{-3}$ & $3.99$ & $2.69\cdot 10^{-5}$ & $4.98$ & $1.43\cdot 10^{-1}$ & $4.00$\\
\hline
\end{tabular}
\end{table}

\subsection{Test case 2: circular inclusion with a compliant interface.}
\label{sec-num2}

We now turn to the inclusion problem of \cite[Sect.~5.1]{HH04}, for which the
solution is radially symmetric, and we extend its exact solution to a
compliant interface.  In contrast with \S\ref{sec-num1}, the interface is here
curved, and the two ways of accounting for it --- exactly fitted meshes with
curved faces, or straight-faced meshes that only approximate the geometry ---
the first takes us beyond the validity of the analysis of \S 4--\S 5; in the second the geometry error dominates for all $k \ge 1$.
Here we consider the first approach, this allows us to study the robustness of the method under non-admissible perturbations. 
 Let $\Gamma=\{r=a\}$ with $r=\abs{\udl{x}}$, let
$\Omega_1$ be the subdomain $\{r<a\}$ and $\Omega_2=\{a<r<b\}$, and let
$\udl{n}=\udl{e}_r$.  Writing $\udl{u}=u_r(r)\udl{e}_r$, $u_\theta=0$, and
$\udl{f}=\udl{0}$, the general solution of \reff{elasticity-interface} in
$\Omega_i$ is $u_r(r)=A_ir+B_i/r$, with
$\sigma_{rr}=2(\lambda_i+\mu_i)A_i-2\mu_iB_i/r^2$ and $\sigma_{r\theta}=0$.
Discarding the singular mode in the inclusion ($B_1=0$) and imposing
$\udl{u}=\udl{x}$ on $\{r=b\}$, the continuity of the normal traction and the
slip condition $\jm{\udl{u}}=-\bK\Vs(\udl{u})\udl{n}$ on $\Gamma$ yield the
linear system
\begin{equation}\label{radial}
  \left\{
  \ba{l}
    A_2+b^{-2}B_2=1,\\
    (\lambda_1+\mu_1)A_1-(\lambda_2+\mu_2)A_2+a^{-2}\mu_2B_2=0,\\
    \big(a+2\beta(\lambda_1+\mu_1)\big)A_1-aA_2-a^{-1}B_2=0,
  \ea
  \right.
\end{equation}
for the three coefficients $(A_1,A_2,B_2)$.  Since $\sigma_{r\theta}=0$, only
the normal compliancy $\beta$ enters \reff{radial}; the tangential compliancy
$\alpha$ is nevertheless active in the discrete problem, through $\bS_h$ and
$\bK$.  For $\beta=0$ one recovers the exact solution of
\cite[Sect.~5.1]{HH04}, namely $A_1=(1-b^2/a^2)c+b^2/a^2$, $A_2=c$ and
$B_2=(1-c)b^2$ with
$c=(\lambda_1+\mu_1+\mu_2)b^2/\big((\lambda_2+\mu_2)a^2+(\lambda_1+\mu_1)
(b^2-a^2)+\mu_2b^2\big)$.  Following \cite{HH04} we take $E_1=1$,
$\nu_1=0.25$ in the inclusion, $E_2=10$, $\nu_2=0.3$ in the surrounding
material, and $a=0.4$, $b=2$.

The computations are run on the quarter annulus
$\Omega=\{0.2<r<2,\ 0<\theta<\pi/2\}$, the Dirichlet data on $\p\Omega$ being
given by the exact solution, whose restriction to $\Omega$ solves the same
problem \reff{elasticity-interface}; excluding a neighbourhood of the origin
avoids the geometric degeneracy of a polar mesh there.  The mesh sequence,
depicted in Fig.~\ref{fig-meshes}(a), consists of $9\cdot 2^{\ell}$ radial
layers by $8\cdot 2^{\ell}$ angular sectors, $\ell\in\{0,1,2,3\}$, the radius
$r=a$ being always a mesh line; the cells are annular sectors, bounded by two
circular and two straight faces, and $\cTh$ is thus fitted to $\Gamma$ in the
sense of \S 3, the curved interface being resolved exactly.  On a circular
face, $\mP^k(F)$ is understood as the space spanned by the polynomials of
degree at most $k$ in the arclength.

\begin{rem}[On resolving $\Gamma$ exactly]\label{Rem-geom}
Approximating $\Gamma$ by the chords of the circle, as is done for convenience
in \cite[Sect.~5]{HH04}, perturbs the geometry by $O(h^2)$ and therefore caps
the accuracy of the method, whatever the polynomial degree.  This is
immaterial for the piecewise affine approximation used in \cite{HH04}, but it
would hide the high-order convergence predicted by
\Pro{\ref{Pro-energy-Err}}--\Pro{\ref{Pro-L2-Err}} as soon as $k\geq 1$. The above theory requires flat interfaces,
however we included this example to see how fragile the method is. It turns out the perturbation is small enough that optimal convergence is observed.
For a detailed analysis on how to extend the HHO method to the case of curved interfaces we refer to \cite{Yemm24}.

\end{rem}

Fig.~\ref{fig-solution} displays the magnitude of the discrete displacement
for $\alpha=\beta=0.5$ together with the radial profiles obtained for
$\alpha=\beta\in\{0,5\cdot10^{-2},0.5\}$ on the coarsest mesh: the slip at the
interface, which increases with the compliancy, is captured on the coarsest
mesh already.  The errors $E_{e}$ and $E_{0}$ are collected in
Table~\ref{tab-conv2} and displayed in Fig.~\ref{fig-conv2} for
$k\in\{1,2,3\}$, both for the perfectly bonded interface ($\alpha=\beta=0$)
and for a compliant interface ($\alpha=\beta=5\cdot10^{-2}$).  The observed
rates, $k+1$ for $E_{e}$ and $k+2$ for $E_{0}$, are in full agreement with
\Pro{\ref{Pro-energy-Err}}--\Pro{\ref{Pro-L2-Err}}, and the two sets of
errors are almost indistinguishable, which is a first indication of the
robustness of the method with respect to the interface data, in spite of the curved interfaces.  A slight
superconvergence (rate $\approx 4.1$ instead of $4$) is observed for $E_{0}$
when $k=2$.  We also note that, for $k=1$, the second-order convergence of the
$L^2$-error reported in \cite[Fig.~3]{HH04} for the piecewise affine unfitted
method is here improved to third order.

Table~\ref{tab-robust} addresses the robustness with respect to the compliancy
parameters.  Since the exact solution \reff{radial} itself depends on $\beta$,
we report there the relative errors
$\tilde{E}_{e}:=E_{e}/\norme{\updl{I}_h(\udl{u})}_{e}$ and
$\tilde{E}_{0}:=E_{0}/\norm{\Pi_h\udl{u}}_{\Omega}$, for $\alpha=\beta$
ranging over sixteen orders of magnitude and for the degenerate combinations
in which one of the two parameters vanishes, which are the cases that motivate
the perturbation $\bf P$ in the definition \reff{penalty-matrix} of $\bS_h$.
Recall from \Rem{\ref{Rem-limits}} that
$\udl{n}^T\bS_h\udl{n}=(\beta+h_T/(\mu+\lambda))^{-1}$ and
$\udl{t}_i^T\bS_h\udl{t}_i=(\alpha+h_T/\mu)^{-1}$, so that, since
$h_T/(\mu_i+\lambda_i)$ and $h_T/\mu_i$ range from $5.6\cdot10^{-3}$ to
$5.2\cdot10^{-1}$ on the cells adjacent to $\Gamma$ in the three meshes used
there, the whole range from the Nitsche-like penalty regime
$\alpha,\beta\ll h_T/\mu$ to the compliancy-dominated regime
$\alpha,\beta\gg h_T/\mu$ is covered.  Both the errors and the rates are
essentially independent of $(\alpha,\beta)$ -- the relative errors vary by at
most $30\%$ over the whole table -- as predicted by the analysis, in which the
constants do not depend on the compliancy.

\begin{rem}[The limit of a free interface]\label{Rem-free}
When $\alpha,\beta\to+\infty$ one has $\bI-\bK\bS_h\to \bf{0}$ and
$\bS_h\to\bf{0}$, so that the interface unknowns disappear from
\reff{gradient-reconstruct}, \reff{Local-aT} and \reff{stabilization-term},
and the discrete problem degenerates into two independent traction-free
problems set on $\Omega_1$ and $\Omega_2$, as anticipated in
\Rem{\ref{Rem-limits}}.  For very large compliancies the interface unknowns
become numerically inactive -- the corresponding diagonal entries of the
condensed matrix drop below $10^{-12}$ times the largest one -- and are then
simply removed from the linear system; the number of such degrees of freedom
is reported in the last column of Table~\ref{tab-robust}.  On the mesh
sequence used there this occurs for $\alpha=\beta\gtrsim10^{10}$, whereas the
free-interface regime is reached much earlier: for $k=1$ on the second mesh
one obtains $E_{e}=8.6015\cdot10^{-2}$, $E_{h}=7.9552\cdot10^{-2}$ and
$E_{0}=4.2861\cdot10^{-4}$ for every $\alpha=\beta$ between $10^{4}$ and
$10^{16}$, to five significant digits, the two subdomains being then
decoupled.  The row $\alpha=\beta=10^{12}$ of Table~\ref{tab-robust}, for
which all the interface unknowns are inactive, gives errors identical to those
of the row $\alpha=\beta=10^{4}$, for which none is: the elimination is
harmless.
\end{rem}

\begin{table}[htbp]
\centering
\small
\caption{Test case 2 (circular inclusion, polar meshes). Errors \reff{errors} and observed convergence rates.}
\label{tab-conv2}
\begin{tabular}{rrrr|rr|rr}
\hline
$\mathrm{card}(\cTh)$ & $\mathrm{card}(\cF_h)$ & $N_{\mathrm{dof}}$ & $h$ & $E_{e}$ & rate & $E_{0}$ & rate\\
\hline
\multicolumn{8}{l}{\emph{$\alpha=\beta=0$}, $k=1$}\\
$72$ & $161$ & $508$ & $0.422$ & $2.30\cdot 10^{-1}$ & -- & $2.68\cdot 10^{-3}$ & --\\
$288$ & $610$ & $2{,}168$ & $0.216$ & $6.83\cdot 10^{-2}$ & $1.81$ & $3.47\cdot 10^{-4}$ & $3.04$\\
$1{,}152$ & $2{,}372$ & $8{,}944$ & $0.109$ & $1.82\cdot 10^{-2}$ & $1.94$ & $4.41\cdot 10^{-5}$ & $3.02$\\
$4{,}608$ & $9{,}352$ & $36{,}320$ & $0.055$ & $4.66\cdot 10^{-3}$ & $1.98$ & $5.58\cdot 10^{-6}$ & $3.01$\\
\hline
\multicolumn{8}{l}{\emph{$\alpha=\beta=0$}, $k=2$}\\
$72$ & $161$ & $762$ & $0.422$ & $5.03\cdot 10^{-2}$ & -- & $4.08\cdot 10^{-4}$ & --\\
$288$ & $610$ & $3{,}252$ & $0.216$ & $7.73\cdot 10^{-3}$ & $2.79$ & $2.86\cdot 10^{-5}$ & $3.96$\\
$1{,}152$ & $2{,}372$ & $13{,}416$ & $0.109$ & $1.02\cdot 10^{-3}$ & $2.97$ & $1.75\cdot 10^{-6}$ & $4.09$\\
$4{,}608$ & $9{,}352$ & $54{,}480$ & $0.055$ & $1.28\cdot 10^{-4}$ & $3.01$ & $1.05\cdot 10^{-7}$ & $4.09$\\
\hline
\multicolumn{8}{l}{\emph{$\alpha=\beta=0$}, $k=3$}\\
$72$ & $161$ & $1{,}016$ & $0.422$ & $1.08\cdot 10^{-2}$ & -- & $9.17\cdot 10^{-5}$ & --\\
$288$ & $610$ & $4{,}336$ & $0.216$ & $9.74\cdot 10^{-4}$ & $3.59$ & $3.75\cdot 10^{-6}$ & $4.76$\\
$1{,}152$ & $2{,}372$ & $17{,}888$ & $0.109$ & $7.14\cdot 10^{-5}$ & $3.83$ & $1.28\cdot 10^{-7}$ & $4.95$\\
$4{,}608$ & $9{,}352$ & $72{,}640$ & $0.055$ & $4.78\cdot 10^{-6}$ & $3.93$ & $4.07\cdot 10^{-9}$ & $5.01$\\
\hline
\multicolumn{8}{l}{\emph{$\alpha=\beta=5\cdot10^{-2}$}, $k=1$}\\
$72$ & $161$ & $508$ & $0.422$ & $2.35\cdot 10^{-1}$ & -- & $2.71\cdot 10^{-3}$ & --\\
$288$ & $610$ & $2{,}168$ & $0.216$ & $6.99\cdot 10^{-2}$ & $1.80$ & $3.55\cdot 10^{-4}$ & $3.03$\\
$1{,}152$ & $2{,}372$ & $8{,}944$ & $0.109$ & $1.87\cdot 10^{-2}$ & $1.93$ & $4.54\cdot 10^{-5}$ & $3.01$\\
$4{,}608$ & $9{,}352$ & $36{,}320$ & $0.055$ & $4.82\cdot 10^{-3}$ & $1.97$ & $5.78\cdot 10^{-6}$ & $3.00$\\
\hline
\multicolumn{8}{l}{\emph{$\alpha=\beta=5\cdot10^{-2}$}, $k=2$}\\
$72$ & $161$ & $762$ & $0.422$ & $5.14\cdot 10^{-2}$ & -- & $4.27\cdot 10^{-4}$ & --\\
$288$ & $610$ & $3{,}252$ & $0.216$ & $7.88\cdot 10^{-3}$ & $2.79$ & $2.95\cdot 10^{-5}$ & $3.98$\\
$1{,}152$ & $2{,}372$ & $13{,}416$ & $0.109$ & $1.04\cdot 10^{-3}$ & $2.96$ & $1.77\cdot 10^{-6}$ & $4.12$\\
$4{,}608$ & $9{,}352$ & $54{,}480$ & $0.055$ & $1.32\cdot 10^{-4}$ & $3.00$ & $1.06\cdot 10^{-7}$ & $4.10$\\
\hline
\multicolumn{8}{l}{\emph{$\alpha=\beta=5\cdot10^{-2}$}, $k=3$}\\
$72$ & $161$ & $1{,}016$ & $0.422$ & $1.11\cdot 10^{-2}$ & -- & $9.46\cdot 10^{-5}$ & --\\
$288$ & $610$ & $4{,}336$ & $0.216$ & $9.98\cdot 10^{-4}$ & $3.58$ & $3.82\cdot 10^{-6}$ & $4.78$\\
$1{,}152$ & $2{,}372$ & $17{,}888$ & $0.109$ & $7.35\cdot 10^{-5}$ & $3.82$ & $1.29\cdot 10^{-7}$ & $4.97$\\
$4{,}608$ & $9{,}352$ & $72{,}640$ & $0.055$ & $4.95\cdot 10^{-6}$ & $3.92$ & $4.11\cdot 10^{-9}$ & $5.01$\\
\hline
\end{tabular}
\end{table}
\begin{table}[htbp]
\centering
\footnotesize
\caption{Test case 2. Robustness with respect to the compliancy parameters: relative errors on the third polar mesh ($\mathrm{card}(\cTh)=1{,}152$, $h=1.09\cdot10^{-1}$). Whatever $(\alpha,\beta)$ in this table, the rate observed on the last refinement lies in $[1.92,1.94]$, $[2.96,2.97]$, $[3.81,3.83]$ for $\tilde{E}_{e}$ and in $[2.98,3.03]$, $[4.09,4.13]$, $[4.92,4.97]$ for $\tilde{E}_{0}$, for $k=1,2,3$ respectively (the asymptotic values $k+1$ and $k+2$ are approached on the finer meshes of Table~\ref{tab-conv2}). The last column reports the number of interface degrees of freedom that are inactive (for $k=1$, $2$ and $3$ respectively), cf.\ \Rem{\ref{Rem-free}}: for $\alpha=\beta=10^{12}$ these are all of them, and the errors are unchanged.}
\label{tab-robust}
\begin{tabular}{rr|rr|rr|rr|r}
\hline
& & \multicolumn{2}{c|}{$k=1$} & \multicolumn{2}{c|}{$k=2$} & \multicolumn{2}{c|}{$k=3$} & \\
$\alpha$ & $\beta$ & $\tilde{E}_{e}$ & $\tilde{E}_{0}$ & $\tilde{E}_{e}$ & $\tilde{E}_{0}$ & $\tilde{E}_{e}$ & $\tilde{E}_{0}$ & \#\,inact.\\
\hline
$0$ & $0$ & $2.66\cdot 10^{-3}$ & $1.64\cdot 10^{-5}$ & $1.49\cdot 10^{-4}$ & $6.48\cdot 10^{-7}$ & $1.04\cdot 10^{-5}$ & $4.73\cdot 10^{-8}$ & $0$\\
$10^{-4}$ & $10^{-4}$ & $2.66\cdot 10^{-3}$ & $1.64\cdot 10^{-5}$ & $1.49\cdot 10^{-4}$ & $6.48\cdot 10^{-7}$ & $1.04\cdot 10^{-5}$ & $4.73\cdot 10^{-8}$ & $0$\\
$10^{-2}$ & $10^{-2}$ & $2.66\cdot 10^{-3}$ & $1.64\cdot 10^{-5}$ & $1.48\cdot 10^{-4}$ & $6.50\cdot 10^{-7}$ & $1.04\cdot 10^{-5}$ & $4.75\cdot 10^{-8}$ & $0$\\
$10^{0}$ & $10^{0}$ & $3.18\cdot 10^{-3}$ & $1.94\cdot 10^{-5}$ & $1.77\cdot 10^{-4}$ & $7.49\cdot 10^{-7}$ & $1.25\cdot 10^{-5}$ & $5.39\cdot 10^{-8}$ & $0$\\
$10^{4}$ & $10^{4}$ & $3.34\cdot 10^{-3}$ & $2.04\cdot 10^{-5}$ & $1.86\cdot 10^{-4}$ & $7.87\cdot 10^{-7}$ & $1.32\cdot 10^{-5}$ & $5.66\cdot 10^{-8}$ & $0$\\
$10^{12}$ & $10^{12}$ & $3.34\cdot 10^{-3}$ & $2.04\cdot 10^{-5}$ & $1.86\cdot 10^{-4}$ & $7.87\cdot 10^{-7}$ & $1.32\cdot 10^{-5}$ & $5.66\cdot 10^{-8}$ & $128$/$192$/$256$\\
$0$ & $10^{-2}$ & $2.67\cdot 10^{-3}$ & $1.65\cdot 10^{-5}$ & $1.49\cdot 10^{-4}$ & $6.58\cdot 10^{-7}$ & $1.05\cdot 10^{-5}$ & $4.82\cdot 10^{-8}$ & $0$\\
$10^{-2}$ & $0$ & $2.65\cdot 10^{-3}$ & $1.63\cdot 10^{-5}$ & $1.47\cdot 10^{-4}$ & $6.42\cdot 10^{-7}$ & $1.04\cdot 10^{-5}$ & $4.67\cdot 10^{-8}$ & $0$\\
$0$ & $10^{0}$ & $3.20\cdot 10^{-3}$ & $2.00\cdot 10^{-5}$ & $1.80\cdot 10^{-4}$ & $8.15\cdot 10^{-7}$ & $1.26\cdot 10^{-5}$ & $5.89\cdot 10^{-8}$ & $0$\\
$10^{0}$ & $0$ & $2.64\cdot 10^{-3}$ & $1.61\cdot 10^{-5}$ & $1.47\cdot 10^{-4}$ & $6.21\cdot 10^{-7}$ & $1.04\cdot 10^{-5}$ & $4.46\cdot 10^{-8}$ & $0$\\
\hline
\end{tabular}
\end{table}

\begin{figure}[htbp]
\centering
\includegraphics[width=\textwidth]{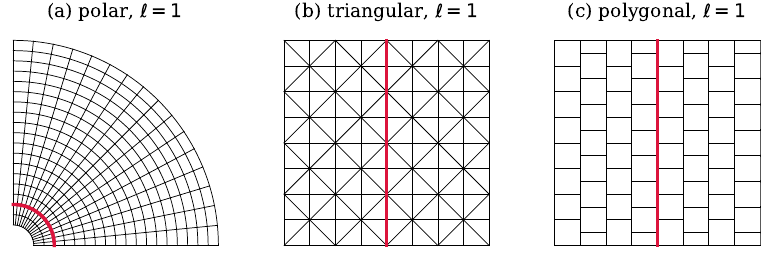}
\caption{The three mesh families ($\ell=1$), with the interface faces $\cIh$
highlighted: (a) polar mesh of the quarter annulus of \S\ref{sec-num2};
(b) triangular and (c) polygonal meshes of the unit square of
\S\ref{sec-num1}.}
\label{fig-meshes}
\end{figure}

\begin{figure}[htbp]
\centering
\includegraphics[width=\textwidth]{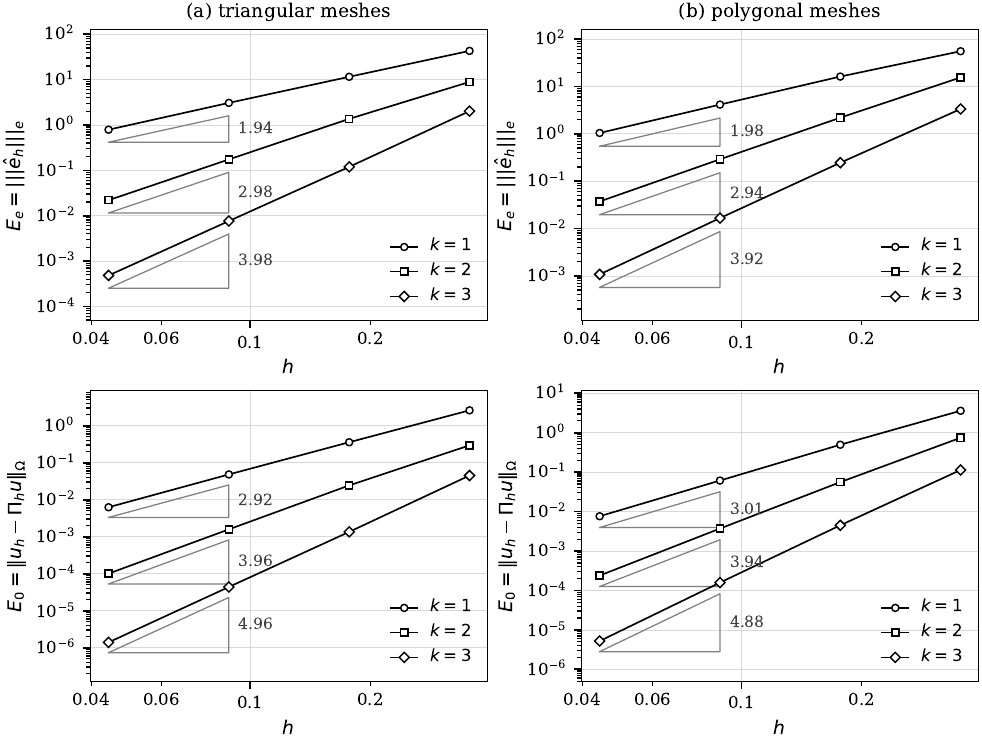}
\caption{Test case 1 (straight interface, $\alpha=\beta=5\cdot10^{-2}$).
Energy error $E_{e}$ (top) and $L^2$-error $E_{0}$ (bottom) versus $h$ on
(a) the triangular and (b) the polygonal mesh family.}
\label{fig-conv1}
\end{figure}

\begin{figure}[htbp]
\centering
\includegraphics[width=\textwidth]{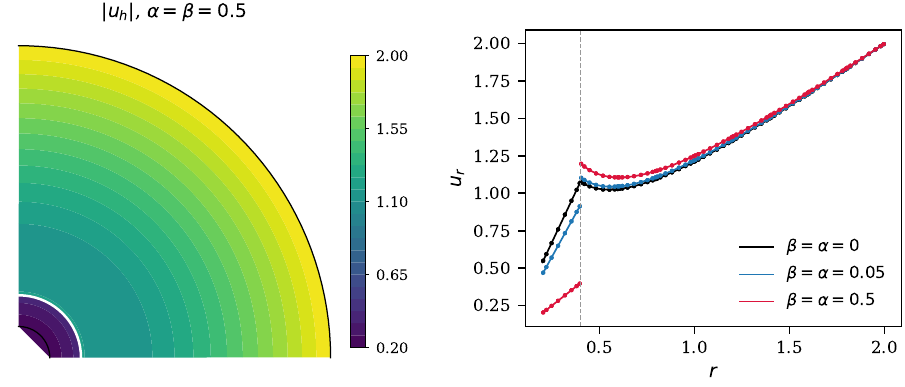}
\caption{Test case 2. Left: magnitude of the discrete displacement ($k=2$,
$\ell=2$, $\alpha=\beta=0.5$), the interface being drawn in white. Right:
radial profiles of $u_r$; solid lines are the exact solutions \reff{radial}
and the dots are the discrete solution on the coarsest mesh ($k=2$,
$\ell=0$), for $\alpha=\beta\in\{0,5\cdot10^{-2},0.5\}$. The slip
$\jm{\udl{u}}=-\bK\Vs(\udl{u})\udl{n}$ at $r=a$ grows with the compliancy.}
\label{fig-solution}
\end{figure}

\begin{figure}[htbp]
\centering
\includegraphics[width=\textwidth]{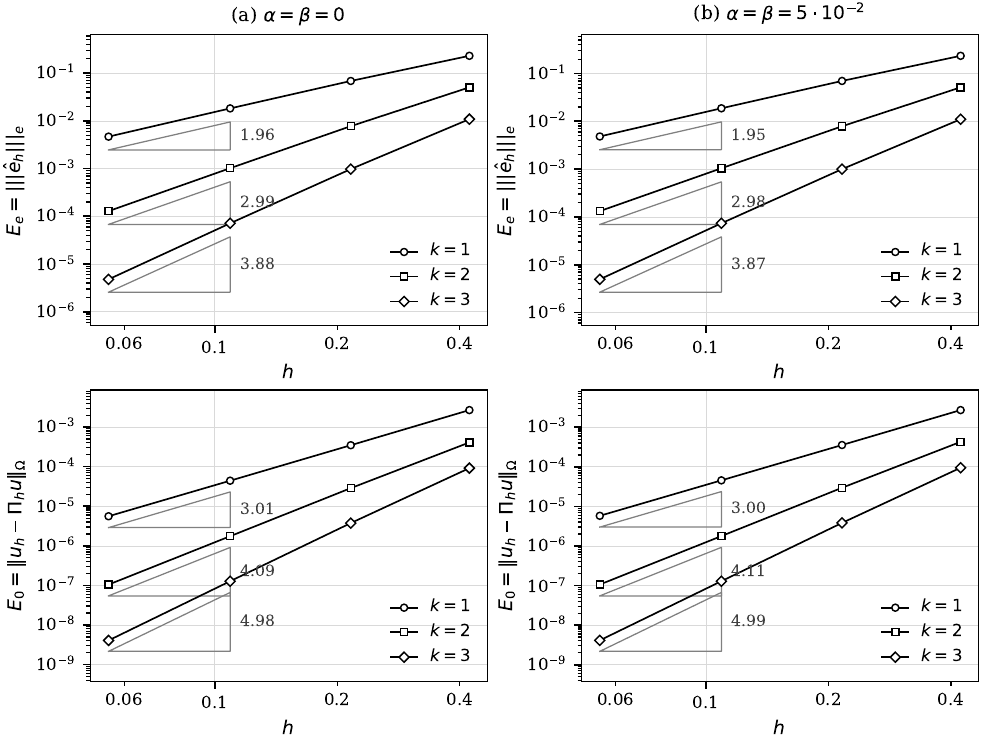}
\caption{Test case 2 (circular inclusion, polar meshes). Energy error $E_{e}$
(top) and $L^2$-error $E_{0}$ (bottom) versus $h$, for (a) a perfectly bonded
interface and (b) a compliant interface. The numbers next to the triangles are
the least-squares rates computed on the three finest meshes.}
\label{fig-conv2}
\end{figure}

\Section{Conclusions.}

We have designed and analysed a hybrid high-order method for the elasticity
problem with a linear slip interface, on meshes fitted to the interface.  The
method uses polynomials of degree $k\ge1$ on the mesh faces and of degree
$k+1$ in the mesh cells, it supports general polytopal cells and hanging
nodes, and the cell unknowns can be eliminated locally.  Its two specific
ingredients are the local symmetric strain reconstruction
\reff{gradient-reconstruct}, into which the interface condition is built
through the factor $\bI-\bK\bS_h$, and the interface stabilization of
\reff{stabilization-term}, weighted by the regularised interface stiffness
$\bS_h$ of \reff{penalty-matrix}.  As a consequence, one single formulation
covers the whole range of compliancies, and no unknown is attached to the
jump of the displacement.  We have proved that the discrete bilinear form is
coercive on the discrete space and that the errors converge as $h^{k+1}$ in
the energy norm and as $h^{k+2}$ in the $L^2$ norm, with constants that are
independent of the compliancy parameters $\alpha,\beta$ and of the Lam\'e
coefficient $\lambda$.  The numerical experiments of \S\ref{sec-numerics}
confirm these rates for $k\in\{1,2,3\}$ on three different mesh families, as
well as the robustness of the method with respect to the compliancy over
sixteen orders of magnitude, including the degenerate cases $\alpha=0$ or
$\beta=0$ and the free-interface limit, and in the quasi-incompressible limit.

\section*{Declarations}
\paragraph{Author contributions} All authors contributed equally to this work. All authors read and approved the final manuscript.

\paragraph{Funding} This work was supported by the National Natural Science Foundation of China grant 11301267, the Natural Science Foundation of Jiangsu Province grant BK20191386 and the Qing Lan Project of Jiangsu
Province. EB was partially supported by EPSRC grants EP/P01576X/1 and EP/V050400/1.

\paragraph{Data availability}
The code reproducing all numerical examples of this paper, together with a script that recomputes every table and figure and verifies them against the published values, is openly available at \cite{burman_huang_hho_slip_code}:\\[3mm]
{\tt{ https://github.com/burmanerik/hho-slip-interface}}\\[3mm]
 and archived at:\\[3mm]
{\tt{  https://doi.org/10.5281/zenodo.22809513.}}

\paragraph{Conflict of interest} The authors declare that they have no competing interests.

\paragraph{Acknowledgement}
During the preparation of this manuscript, the authors used Anthropic Claude Opus 5, to assist with drafting and revising text, improving mathematical exposition, checking notation and internal consistency, editing LaTeX, and developing computational code. These tools were not used to fabricate or directly alter research data or numerical results. All AI-assisted material, including mathematical statements and computational code, was critically reviewed and validated by the authors, who take full responsibility for the accuracy, originality, and integrity of the manuscript.

\bibliographystyle{plain}
\bibliography{ref1}
\end{document}